\documentclass[reqno,10pt, centertags]{amsart}
\usepackage{amssymb,upref,esint,color}
\usepackage{hyperref}
\usepackage{etoolbox}
\newcommand*{\mailto}[1]{\href{mailto:#1}{\nolinkurl{#1}}}

\makeatletter
\usepackage{verbatim}
\usepackage{amscd}
\usepackage{framed}
\numberwithin{equation}{section}
\input epsf

\newtheorem{theorem}{Theorem}[section]
\newtheorem{definition}[theorem]{Definition}
\newtheorem{lemma}[theorem]{Lemma}
\newtheorem{corollary}[theorem]{Corollary}
\newtheorem{proposition}[theorem]{Proposition}
\numberwithin{equation}{section}
\theoremstyle{definition}
\newtheorem{example}[theorem]{Example}

\begin{document}

\title[Complex analysis of symmetric operators. I]
{Complex analysis of symmetric operators. I }

\author[Yicao Wang]{Yicao Wang}
\address
{School of Mathematics, Hohai University, Nanjing 210098, China}

\baselineskip= 20pt
\begin{abstract}Based on the relationship of symmetric operators with Hermitian symmetric spaces, we introduce the notion of \emph{Weyl curve} for a symmetric operator $T$, which is the geometric abstraction and generalization of the well-known Weyl functions. We prove that there is a one-to-one correspondence between unitary equivalence classes of simple symmetric operators and congruence classes of Nevanlinna curves (the geometric analogue of operator-valued Nevanlinna functions). To prove this result, we introduce a \emph{canonical} functional model for $T$ in terms of its \emph{characteristic vector bundles}. In this geometric formalism, when the deficiency indices are $(n,n)$ ($n\leq +\infty$) we also introduce the notion of \emph{entire operators}, whose Weyl curves are entire in the Grassmannian $Gr(n,2n)$. If $n$ is finite, this makes it possible to introduce modern value distribution theory of entire curves into the picture and to demonstrate that the distribution of eigenvalues of a \emph{generic} abstract boundary value problem is the same thing as the value distribution of the Weyl curve with respect to the Cartier divisor induced by the corresponding boundary condition. Many other new concepts are introduced and many new results are obtained as well.
\end{abstract}
\maketitle


\tableofcontents

\section{Introduction}
Unbounded self-adjoint operators play a basic role in quantum mechanics. It was von Neumann who first distinguished self-adjoint operators from symmetric ones in his studies of the mathematical foundations of quantum mechanics in late 1920s. von Neumann then invented the self-adjoint extension theory of symmetric operators. He found that a symmetric operator has a self-adjoint extension if and only if the positive deficiency index equals to the negative deficiency index. von Neumann also succeeded in finding all self-adjoint extensions if they really exist. These are parameterized precisely by unitary maps from the positive deficiency subspace to the negative. The later years have witnessed applications of this theory in many other areas.

However, von Neumann's theory is an abstract framework and people are still interested in how to parameterize the self-adjoint extension in a convenient manner for various concrete applications. There are several schemes trying to achieve this, e.g., using Lagrangian subspaces, self-adjoint linear relations, unitary operators and even unbounded operators \cite{schmudgen2012unbounded, behrndt2020boundary, grubb1968characterization}.

In our viewpoint, even if one is only interested in self-adjoint extensions, non-self-adjoint extensions are essentially unavoidable in the big picture. For simplicity, right now let $T$ be a densely defined simple closed symmetric operator with deficiency indices $(n,n)$, where $n\in \mathbb{N}$. A closed extension of $T$ is actually determined uniquely by a closed subspace (we also view it as an abstract boundary condition, bearing differential operators in our mind) of the quotient $D(T^*)/D(T)$, where $D(T)$ is the domain of $T$ and $D(T^*)$ is that of its adjoint operator $T^*$. Consequently the most \emph{canonical} and \emph{natural} way to parameterize extensions of $T$ is to use the set of closed subspaces of $D(T^*)/D(T)$--the Grassmannian of closed subspaces in $D(T^*)/D(T)$. There are several topological components of this space, among which the most important is the Grassmannian of closed subspaces of dimension $n$. It can be holomorphically identified with $Gr(n, 2n)$--the Grassmannian of $n$-dimensional subspaces in $\mathbb{C}^{2n}$.

Another important concept in the theory of self-adjoint extensions is the Weyl function. It has its origin in H. Weyl's work on singular Sturm-Liouville problems on the half line in early 1910s and since then has played an essential role in the direct or inverse spectral theory of Sturm-Liouville problems \cite{weyl1910gewohnliche}. This construction was then generalized by M. G. Krein's school in Ukraine for general symmetric operators with equal deficiency indices in 1980s \cite{derkach1987weyl}. See the Introduction of \cite{behrndt2020boundary} for a more detailed historical account. However, in the general setting the Weyl function is an analytic operator-valued function on the upper and lower half planes $\mathbb{C}_+\cup \mathbb{C}_-$. Its definition is always related to a specific choice of data called a \emph{boundary triplet} introduced in 1970s \cite{VMBruk1976, kochubei1975extensions}, involving partly two specific self-adjoint extensions. In this sense, the Weyl function is never uniquely defined and often develops singularities at some points on the real line.

 From our viewpoint, a Weyl function is really not a function, but essentially a geometric object (we call it the \emph{Weyl curve}) defined \emph{intrinsically} by the symmetric operator $T$ itself, and has no direct relation to any specific self-adjoint extensions. We note that this curve is actually a holomorphic map $W_T(\lambda)$ from $\mathbb{C}_+\cup \mathbb{C}_-$ to the above $Gr(n, 2n)$ and in certain cases the curve can even be defined holomorphically on the whole $\mathbb{C}$, i.e., the Weyl curve is an entire curve in $Gr(n, 2n)$.

An important and fundamental observation in this paper is that the position of the image $W_T(\mathbb{C}_+)$ in $Gr(n, 2n)$ is very special. It lies in an open subset $M$ of $Gr(n, 2n)$, which is the famous non-compact irreducible Hermitian symmetric space of type $I_{n,n}$ while $Gr(n, 2n)$ is its compact dual. This observation explains the basic properties of Weyl functions often mentioned in the literature. $M$ can be realized as a bounded symmetric domain whose Bergman metric is hyperbolic while $\mathbb{C}_+$ itself (with its Poincar$\acute{e}$ metric) is a model of hyperbolic plane. This implies that \emph{the geometry of Hermitian symmetric spaces and holomorphic maps between them may serve as the very foundation of a complete treatment of symmetric operators}.

It is this belief that leads us to introduce complex analysis (and geometry) into the study of the extension theory of symmetric operators. The basic goal of this paper is to explore in depth von Neumann's theory in this direction. We succeed in generalizing the notion of Weyl functions to symmetric operators with arbitrary deficiency indices $(n_+,n_-)$ (but neither of $n_\pm$ is zero\footnote{If either $n_+$ or $n_-$ is zero, $T$ is called maximal. The structure of such operators is clear, see \cite[\S~104, Chap.~8]{akhiezer2013theory}}). As one may expect, this is related to the irreducible Hermitian symmetric space $N$ of type $I_{n_+,n_-}$ \footnote{Strictly speaking, if either of $n_\pm$ is $+\infty$, then we should call $N$ a Banach symmetric space, but we won't make this distinction.}. One of our main theorems is that, there is a one-to-one correspondence between the unitary equivalence classes of simple symmetric operators and the congruence classes of \emph{Nevanlinna curves} (the geometric abstraction of matrix-valued Nevanlinna functions) under the action of the automorphism group of $N$ (Thm.~\ref{thm12}). Somehow, this fact seems to be known in the community of spectral theorists for many years at least in the case $n_+=n_-$, though not in this geometric form. We think it is our introduction of the notion of Weyl curves that makes this basic result much more transparent.

To prove and understand the above theorem, we also introduce three characteristic vector bundles associated to a simple symmetric operator $T$. These vector bundles arise naturally and are a geometric consequence of von Neumann's basic observation that the deficiency index is locally constant on $\mathbb{C}_+\cup \mathbb{C}_-$. The holomorphicity of the Weyl curve then simply means these are holomorphic bundles. In a sense, Weyl curves are just the classifying maps of these vector bundles and the Hermitian symmetric space $N$ provides the classifying space, just as the common sense that general Grassmannians are classifying spaces for complex vector bundles.

One use of the characteristic bundles is to construct a functional model for $T$. A functional model of $T$ is to construct a Hilbert space of holomorphic functions on $\mathbb{C_+}\cup\mathbb{C_-}$ such that multiplication by the independent variable is symmetric and unitarily equivalent to $T$. In our model, the Hilbert space consists of certain holomorphic sections of our characteristic vector bundle of the second kind. The advantage of our functional model is that it is intrinsically defined in the sense that it only depends on $T$ itself, without resorting to any artificial choices. Other functional models in the literature can be viewed as different realizations of ours in terms of different trivializations of the characteristic vector bundle.

We pay special attention to the case with equal finite deficiency indices $(n,n)$ for its importance in spectral theory. If the Weyl curve is entire in $Gr(n,2n)$, we say the operator is entire. Our geometric viewpoint towards these operators motivates us to introduce the formalism of value distribution theory of entire curves into the picture. In this theory, people's main concern is the intersection property of an entire curve with Cartier divisors in the ambient projective algebraic manifold. Indeed a \emph{generic} abstract boundary condition produces a Cartier divisor in $Gr(n,2n)$ (called a Schubert hyperplane in algebraic geometry) and when the Weyl curve intersects the Schubert hyperplane, the corresponding abstract boundary value problem obtains an eigenvalue. Though some subtleties may arise, the basic philosophy is that \emph{the distribution of eigenvalues is the same thing as the value distribution of the Weyl curve}. This observation makes it possible to introduce many notions in value distribution theory into our theory. This picture gives us a unified formalism to treat generic extensions of $T$--whether they are self-adjoint or not. The main difficulty we have encountered here is how to extend the investigation to the case of $n=+\infty$. Our success in this respect is only partial and will be presented elsewhere.

The idea of using complex analysis in spectral theory goes back even to the early days of this exciting mathematical discipline, e.g., Poincar$\acute{e}$'s work on spectra of the Laplacian on a bounded domain in $\mathbb{R}^2$ or $\mathbb{R}^3$ with classical boundary conditions in 1890s. Its popularity in spectral theory now is more or less a common sense. However, our viewpoint is that, complex analysis and even complex geometry are an indispensable part of the theory--at least for the extension theory of symmetric operators, not just a useful tool. Note that the material of complex geometry used in this paper is almost standard\footnote{Almost all textbooks on complex geometry only deal with complex vector bundles of finite rank, but we do need those with infinite rank and luckily this won't cause serious problems.} and its details can be found in the great book \cite{griffiths2014principles}.

The paper is organized as follows. In \S~\ref{sec2}, we give a quick introduction to von Neumann's theory, which provides a starting point and is the germ of the subsequent developments. In \S~\ref{sec3}, we introduce the notion of strong symplectic Hilbert spaces and investigate basic properties of the relevant Hermitian symmetric space of type $I_{n_+,n_-}$. The material is not original in any sense, but presented in a manner facilitating the definition of Weyl curves in \S~\ref{sec4}. \S~\ref{sec5} introduces the three characteristic vector bundles mentioned above. The characteristic vector bundle of the third kind is not necessary for proving the basic correspondence between simple symmetric operators and Nevanlinna curves. It is there to answer the question: which vector bundles over $\mathbb{C}_+\cup \mathbb{C}_-$ are the characteristic vector bundles of simple symmetric operators? We haven't given all the details of the answer but relate it to the non-abelian Hodge theory. The correspondence between simple symmetric operators and Nevanlinna curves should be regarded as a classification theorem of simple symmetric operators and the set of congruence classes of Nevanlinna curves should be viewed as a moduli space and its structure shall reflect the sociology of simple symmetric operators. The following sections \S~\ref{sec7}, \S~\ref{sec8} are the first steps towards this direction. However, to give a definite theory goes beyond our ambition and we content ourselves with emphasizing the relevant concepts and giving some elementary observations. In \S~\ref{sec9} and the subsequent sections, we come back to the case $n_+=n_-$. The basic goal of \S~\ref{sec9} is to demonstrate how \emph{contractive} Weyl functions can be used to analyze the spectral theory of generic closed extensions of a simple symmetric operator, while in the literature this is usually done in terms of the Weyl functions. The section is also to pave the way for the following sections. \S~\ref{asing} gives an analytic interpretation of spectral kernel of a simple symmetric operator--points in spectral kernel are precisely those points on the real line at which the Weyl curve doesn't admit analytic continuations. Thus the simplest cases are simple symmetric operators with empty spectral kernel, whose Weyl curves are entire curves defined on the whole $\mathbb{C}$. These so-called entire operators are the main topic in \S~\ref{entire}, where we introduce modern value distribution theory into the analysis of abstract boundary value problems associated to an entire operator with a finite deficiency index. Many notions are introduced and several results are included in the subsections. The content in \S~\ref{real} is a supplement to the proceeding investigation, emphasizing the difference between symmetric operators in real Hilbert spaces and those in complex Hilbert spaces. The last section \S \ref{sturm} applies our theory to Sturm-Liouville problems with the goal of demonstrating the utility of our theory. This sheds new light on some old results and also provides some new ones.

\begin{center}
\textbf{Conventions and notations}
\end{center}

 In this paper, all (real or complex, and finite-dimensional or infinite-dimensional) Hilbert spaces $(H, (\cdot, \cdot))$ are separable. If $H$ is complex, we assume that the inner product $(\cdot, \cdot)$ is linear in the first variable and conjugate-linear in the second. The induced norm will be written as $\|\cdot\|$. If necessary, the inner product or norm is also denoted by $(\cdot, \cdot)_H$ or $\|\cdot\|_H$ to emphasize which underlying Hilbert space is referred to. We use $\oplus$ to denote topological direct sum while orthogonal direct sum will be denoted by $\oplus_\bot$. For a subspace $V$, its orthogonal complement is $V^\bot$. The zero vector space will be just denoted by $0$. For example, for subspaces $V_1, V_2\subset H$ such that $V_1\cap V_2=\{0\}$, we write $V_1\cap V_2=0$. $\mathbb{B}(H_1, H_2)$ denotes the Banach space of bounded operators between the Hilbert spaces $H_1$ and $H_2$. If $H_1=H_2=H$, we simply write this as $\mathbb{B}(H)$. The identity operator on $H$ will be just denoted by $Id$ if the underlying space is clear from the background. For a constant $c$, $c\times Id$ is often written simply as $c$ if no confusion arises. $D(A)$ is the domain of the operator $A$, and $\textup{ker}A$ (resp. $\textup{Ran}A$) is the kernel (resp. range) of $A$. The spectrum of a (bounded or unbounded) operator $A$ is denoted by $\sigma(A)$ and the resolvent set of $A$ by $\rho(A)$. The point (resp. continuous and residual) spectrum  of $A$ is denoted by $\sigma_p(A)$ (resp. $\sigma_c(A)$ and $\sigma_r(A)$). For operators $A$ and $B$, we use $A\subset B$ to mean
$D(A)\subset D(B)$ and $B|_{D(A)}=A$. In this setting, $B$ is called an extension of $A$.

If $\lambda\in \mathbb{C}$, then $\Re \lambda$ and $\Im \lambda$ are the real and the imaginary parts of $\lambda$ respectively. For a matrix $A$ (resp. a complex vector bundle $E$), $\textup{rk}A$ (resp. $\textup{rk}E$) is the rank of $A$ (resp. $E$). If $f(x)$ is a function, $O(f)$ represents a function $g(x)$ such that $|g/f|$ is bounded w.r.t. the underlying limiting process of the independent variable $x$. We shall use the words "holomorphic" and "analytic" interchangeably.

Additionally, some terminology and notation related to group actions are included in Appendix A.

\section{von Neumann's theory revisited}\label{sec2}
This section gives a very brief sketch of von Neumann's theory on self-adjoint extensions. Its basic purpose is to introduce our main object to be studied and to spot our starting point for the whole paper. For technical details of this section, see for example \cite[Chap.~13]{rudin1991functional} or \cite[Chap.~13]{schmudgen2012unbounded}.

Let $H$ be an infinite-dimensinal complex Hilbert space. Recall that a linear operator $A$ defined in $H$ is closed if the graph of $A$ is a closed subspace of $H\oplus_\bot H$. A densely defined symmetric operator $T$ in $H$ is a linear operator from a linear subspace $D(T)\subset H$ to $H$ such that $D(T)$ is dense in $H$ and $(Tx,y)=(x,Ty)$ for any $x,y\in D(T)$.

Let $T$ be a densely defined closed symmetric operator. The domain $D(T^*)$ of the adjoint $T^*$ consists of $y\in H$ such that
\[(Tx,y)=(x,z),\quad \forall x\in D(T)\]
for some (unique) $z\in H$, and $T^*y=z$. $T^*$ is always closed and $T\subset T^*$. If furthermore $T=T^*$, we say $T$ is self-adjoint. If $T\subset S \subset T^*$ and $S=S^*$, we say $S$ is a self-adjoint extension of $T$. The basic concern of the theory of self-adjoint extensions is to determine whether a self-adjoint extension of $T$ exists and to characterize them if they do exist.

Throughout the paper, a symmetric operator always means a closed one. von Neumann's fundamental observation is that $\textup{dim ker}(T^*-\lambda)$ as a function in $\lambda$ is constant on $\mathbb{C}_+$ and $\mathbb{C}_-$ respectively. It is well-known that
\[D(T^*)=D(T)\oplus N_+\oplus N_-,\]
where $N_\pm=\textup{ker}(T^*\mp i)$ are called deficiency subspaces of $T^*$. If $D(T^*)$ is equipped with the \emph{graph inner product} associated with $T^*$, i.e., $$(x,y)_{T^*}=(x,y)+(T^*x,T^*y)$$ for $x,y\in D(T^*)$, then $D(T^*)$ is a Hilbert space and $D(T)$ is a closed subspace. In particular, the above decomposition is actually orthogonal w.r.t. the graph inner product. The numbers $n_\pm:=\textup{dim}N_\pm$ ($\leq +\infty$) are called the positive and the negative deficiency index respectively. Self-adjoint extensions of $T$ exist if and only if $n_+=n_-$. The following theorem is due to von Neumann.
\begin{theorem}(\cite[Thm.~13.10]{schmudgen2012unbounded}) Let $T$ be a symmetric operator with $n_+=n_-\geq 1$. Then there is a one-to-one correspondence between the set $\mathcal{S}$ of self-adjoint extensions of $T$ and the set $\mathrm{U}(N_+,N_-)$ of unitary maps from $N_+$ to $N_-$.
\end{theorem}

We recall briefly how the correspondence (some authors call it the von Neumann map) is constructed. Note that we have the following orthogonal decompositions
\[H=\textup{Ran}(T-i)\oplus_\bot N_-=\textup{Ran}(T+i)\oplus_\bot N_+.\]
From this we can construct a partial isometry (the Cayley transform of $T$):
\[C:=(T-i)(T+i)^{-1}: \textup{Ran} (T+i)\longrightarrow \textup{Ran}(T-i).\]
Define $\tilde{C}: H\rightarrow H$ by $\tilde{C}|_{N_+^\bot}=C$ and $\tilde{C}|_{N_+}=0$. Given a unitary map $U: N_+\rightarrow N_-$, define $\tilde{U}:H \rightarrow H$ by $\tilde{U}|_{N_+}=U$, $\tilde{U}|_{N_+^\bot}=0$. Then
$C_U:=\tilde{C}+\bar{U}$ is a unitary operator on $H$. Solve the equation in $T_U$
\[C_U=(T_U-i)(T_U+i)^{-1}.\]
The resulting $T_U$ is the self-adjoint extension corresponding to $U$. In this way, we get the von Neumann map $v: U\mapsto T_U$.
$U$ can be reasonably called an abstract (self-adjoint) boundary condition and $v$ is precisely how the self-adjoint extensions depend on the corresponding self-adjoint boundary conditions.

The above von Neumann's theory is hardly a completed theory. The point is that the set $\mathrm{U}(N_+,N_-)$ has more natural structures than just being a set. If one such unitary map is fixed, then the set can be identified with $\mathbb{U}(N_+)$, the unitary group of the Hilbert space $N_+$.

On the other side, since we are considering unbounded self-adjoint operators depending on parameters, we'd better choose a topology on the set of (bounded or unbounded) self-adjoint operators. One choice is the gap topology.

\begin{definition} For two self-adjoint operators $A_1$, $A_2$ in a Hilbert space $H$, the gap metric between them is defined via
$$\gamma(A_1, A_2):=\|(A_1+i)^{-1}-(A_2+i)^{-1}\|.$$
The induced topology is called the gap topology.
\end{definition}
 This topology can be viewed from different angles: originally it is the "gap" between the graphs of $A_1$ and $A_2$. One can also map self-adjoint operators to unitary operators via Cayley transform and the metric is precisely the norm metric between the corresponding unitary operators. Indeed, in our context note that
\[C_U=(T_U+i-2i)(T_U+i)^{-1}=Id-2i(T_U+i)^{-1},\]
i.e.
\[(T_U+i)^{-1}=\frac{1}{2i}(Id-C_U).\]
Thus if $U_1, U_2\in \mathrm{U}(N_+, N_-)$, then
\begin{equation*}(T_{U_1}+i)^{-1}-(T_{U_2}+i)^{-1}=\frac{1}{2i}(C_{U_2}-C_{U_1})=\frac{1}{2i}(\tilde{U}_2-\tilde{U}_1).\end{equation*}
From this formula, $v$ is therefore continuous (we use the norm topology on $\mathrm{U}(N_+, N_-)$), and we say the self-adjoint extensions depend continuously on the corresponding boundary conditions. For more information on the above topological aspect, see \cite{booss2005unbounded}.
\section{Strong symplectic structures and related Hermitian symmetric spaces}\label{sec3}
The goal of this section is to present the basics of irreducible Hermitian symmetric spaces of type $I_{n_+,n_-}$. More or less, the material included here is classical and standard since E. Cartan's time. However, Hermitian symmetric spaces are generally introduced in the formalism of Lie theory and realized finally as spaces of complex matrices in the finite-dimensional case \cite{helgason1979differential, viviani2014tour}. We take an alternative approach starting with the notion of strong symplectic structures, and realize the relevant Hermitian symmetric spaces as spaces of operators (thus the infinite-dimensional case can be handled as well). We think this is a short-cut for applying Hermitian symmetric spaces to the extension theory of symmetric operators as it is actually the way how these spaces arise in the theory.

Let $H$ be a complex Hilbert space ($\textup{dim}H\leq +\infty$).
\begin{definition}
 A strong symplectic structure on $H$ is a continuous sesquilinear form $[\cdot, \cdot]: H\times H\rightarrow \mathbb{C}$ such that\\
 i) For $x,y\in H$, $[x, y]$ is linear in $x$ and conjugate-linear in $y$;\\
 ii) $[y,x]=-\overline{[x,y]}$ for all $x,y\in H$;\\
 iii)  $[\cdot, \cdot]$ is non-degenerate in the sense that if $[x,y]=0$ for all $y\in H$, then $x=0$;\\
 iv) the map $\mathfrak{i}:H\rightarrow H^*$ defined by $\mathfrak{i}(x)=[\cdot, x]$ is surjective, where $H^*$ is the dual space of continuous linear functionals on $H$.\\
 The triplet $(H, (\cdot, \cdot), [\cdot, \cdot])$ will be called a strong symplectic Hilbert space and if no ambiguity arises, we just say $H$ is a strong symplectic Hilbert space.
  \end{definition}
\emph{ Remark}. (1) In the literature, when defining a strong symplectic structure, one often needs \emph{not} to specify the choice of an inner product on $H$ and only an equivalence class of norms matters, because this is enough for equipping $H$ with the required topology.\footnote{Perhaps \cite{marsden1974reduction} was the first to define an infinite-dimensional strong symplectic structure in the real case, where no norm or inner product is specified.} However, we do include an inner product as part of the definition, because this is indeed the case we shall encounter. (2) If $\textup{dim} H<+\infty$, the condition iv) is superfluous and follows from iii). (3) A strong symplectic structure is in fact an indefinite inner product on $H$. Indeed, $-i[\cdot, \cdot]$ is Hermitian in the usual sense: $-i[y,x]=\overline{-i[x,y]}$. However, we prefer to use the terminology here to remind the reader of the similarity to real symplectic structures. Further information on this analogy can be found in \S~\ref{sec4} and \S~\ref{real}.

We are not interested in the case where $-i[\cdot, \cdot]$ is positive-definite or negative-definite and in the following we always assume that our symplectic structures are not of this kind. We also note that our definition and investigation to follow have some overlaps with \cite{everitt2004infinite}, which was also motivated by boundary value problems.
 \begin{example} \label{ex2}Given two complex Hilbert spaces $H_+$, $H_-$, of dimension $n_+$ and $n_-$ respectively. Then $H=H_+\oplus_\bot H_-$ can be equipped with a standard strong symplectic structure as follows: for $x=(x_1,x_2), y=(y_1,y_2)\in H_+\oplus_\bot H_-$, define
 \[[x,y]=i(x_1,y_1)_{H_+}-i(x_2,y_2)_{H_-}.\]
  \end{example}
  \begin{example}If $(H_i, (\cdot,\cdot)_{H_i}, [\cdot, \cdot]_i)$, $i=1,2$ are two strong symplectic Hilbert spaces, the direct sum of the two is the orthogonal direct sum $H_1\oplus_\bot H_2$, whose strong symplectic structure is defined as follows: For $x=(x_1,x_2), y=(y_1,y_2)\in H_1\oplus_\bot H_2$,
  \[[x,y]:=[x_1,y_1]_1+[x_2,y_2]_2.\]
  \end{example}
Let $H$ be a strong symplectic Hilbert space.
 \begin{definition}Let $M$ be a (closed or non-closed) subspace of $H$. The symplectic complement of $M$ is defined as the subspace
 \[M^{\bot_s}=\{x\in H|[x,y]=0,\, \forall y\in M\}.\]
\end{definition}
If $M_1$, $M_2$ are two subspaces, then it can be easily checked that \[(M_1+M_2)^{\bot_s}=M_1^{\bot_s}\cap M_2^{\bot_s},\quad (M_1\cap M_2)^{\bot_s}\supset\overline{M_1^{\bot_s}+ M_2^{\bot_s}}.\]
\begin{definition}A closed subspace $M\subset H$ is isotropic if $M\subset M^{\bot_s}$. An isotropic subspace $M$ is called maximal if there is no isotropic subspace of $H$ containing $M$ properly. $M$ is called co-isotropic if $M^{\bot_s}\subset M$. If $M=M^{\bot_s}$, we call $M$ a Lagrangian subspace.
\end{definition}

\begin{definition}\label{de1}A closed subspace $M\subset H$ is called completely positive- (resp. negative-) definite if $-i[\cdot, \cdot]$ is positive- (resp. negative-)definite when restricted on $M$ and there exists a positive constant $c_M$ such that $-i[x, x]\geq c_M \|x\|^2$ (resp. $i[x, x]\geq c_M \|x\|^2$) for all $x\in M$. Such an $M$ is called maximal if it cannot be contained properly in another completely positive- (resp. negative-)definite subspace.
\end{definition}
We can find a canonical maximal completely positive-definite subspace $M_+$ as follows: via Riesz representation theorem, there is a linear bounded operator $A$ on $H$ such that $[x,y]=(Ax,y)$ for all $x,y\in H$. The condition ii) implies that $A^*=-A$ and the conditions iii) and iv) imply $A$ is a topological linear isomorphism. Let $|A|$ be the square root of $A^*A=-A^2$. Since $A$ commutes with $A^2$ and hence with $A^*A$, $A$ commutes with $|A|$. Let $\mathcal{J}=|A|^{-1}A$. Then
\[\mathcal{J}^2=(|A|^{-1}A)^2=|A|^{-2}A^2=-Id.\]
In terms of $\mathcal{J}$ and $|A|$, one can obtain the basic identity $(M^{\bot_s})^{\bot_s}=\bar{M}$ for any subspace $M$, where $\bar{M}$ is the closure of $M$. Denote the closed subspaces $\textup{ker}(\mathcal{J}\mp i)$ by $M_\pm$ respectively.
\begin{proposition}\label{p10}We have the decomposition $H=M_+\oplus_\bot M_-$. $M_+$ (resp. $M_-$) is a maximal completely positive-definite (resp. negative-definite) subspace of $H$, and $M_-=M_+^{\bot_s}$.
\end{proposition}
\begin{proof}
The orthogonal decomposition is obvious. For $x\in H$, let $x=x_++x_-$ be the decomposition accordingly. Denote the restriction of $-i[\cdot, \cdot]$ on $M_+$ by $(\cdot, \cdot)_+$ and the restriction of $i[\cdot, \cdot]$ on $M_-$ by $(\cdot, \cdot)_-$. Note that for $x\in M_+$ and $y\in M_-$,
\[[x,y]=(Ax,y)=(|A|\mathcal{J}x,y)=i(|A|x,y),\]
and similarly,
\[[x,y]=-(x, Ay)=-i(x, |A|y),\]
implying that $[x,y]=0$. Therefore, for $x=x_++x_-$ and $y=y_++y_-$,
\[-i[x,y]=(x_+,y_+)_+-(x_-,y_-)_-.\]
Note that for $0\neq x\in M_+$,
\[-i[x,x]=(x,x)_+=-i(Ax,x)=(|A|x,x)>0.\]
This shows that $M_+$ is completely positive-definite since $|A|$ is a topological linear isomorphism.

If $M_+$ is not maximal, then there exists a completely positive-definite subspace $M'$ containing $M_+$ properly. Let $y=y_++y_-\in M'$ such that $y$ is orthogonal to $M$ in $M'$ (w.r.t. the inner product induced from $-i[\cdot, \cdot]$ on $M'$). Then immediately we have $y_+=0$ and $y=y_-$, contradicting the assumption that $M'$ is completely positive-definite. The conclusion for $M_-$ follows similarly.

From the above argument, it is obvious that $M_-\subset M_+^{\bot_s}$. If $x=x_++x_-\in M_+^{\bot_s}$, then by definition for any $y\in M_+$, we have
\[-i[y,x]=(y,x_+)_+=0,\]
implying that $x_+=0$ and $x=x_-\in M_-.$
 \end{proof}

There are many maximal completely positive-definite (resp. negative-definite) subspaces of $H$. Denote the set of all maximal completely positive-definite (resp. negative-definite) subspaces by $W_+(H)$ (resp. $W_-(H)$). $W_+(H)$ can be characterized in the following way. Note that the above $M_\pm$ are Hilbert spaces with inner product $(\cdot, \cdot)_\pm$ respectively.

\begin{theorem}\label{para}If $H=M_+\oplus_\bot M_-$, where $M_\pm$ are those in Prop.~\ref{p10}, then a subspace $M$ is a maximal completely positive-definite subspace of $H$ if and only if it is of the form $\{x+Bx|x\in M_+\}$ where $B\in \mathbb{B}(M_+,M_-)$ and $\|B\|<1$.
\end{theorem}
\begin{proof} Let $M$ be a maximal completely positive-definite subspace of $H$. If $x\in M$ and $x=x_++x_-$ according to the decomposition $H=M_+\oplus M_-$, then the map $\pi_+:M\rightarrow M_+, x\mapsto x_+$ is a linear topological isomorphism. Indeed, $\pi_+$ is obviously injective. Note that by definition for $x\in M$
\[-i[x,x]=\|x_+\|_+^2-\|x_-\|^2_-\geq c_M\|x\|^2,\]
implying that $\pi_+(M)$ is closed. If $\pi_+$ is not surjective, then there is a nonzero $y\in M_+$ such that $y$ is orthogonal to $\pi_+(M)$ in $M_+$. The linear span of $y$ and $M$ should again be a completely positive-definite subspace of $H$. This contradicts the maximality of $M$. One also obviously has that $\|\pi_+x\|\leq \|x\|$ for $x\in M$. The bounded inverse theorem then implies that $\pi_+$ is an isomorphism.

In this way, $x_-=\pi_+^{-1}(x_+)-x_+$ and we can write this as $x_-=Bx_+$ for $B\in \mathbb{B}(M_+, M_-)$. Note that for any $0\neq x_+\in M_+$
\begin{eqnarray*}-i[x_++Bx_+,x_++Bx_+]&=&(x_+,x_+)_+-(Bx_+,Bx_+)_-\\
&=&((Id-B^*B)x_+,x_+)_+\geq c_M (\|x_+\|_+^2+\|Bx_+\|_-^2).
\end{eqnarray*}
This implies that $B$ should lie in the open unit ball in $\mathbb{B}(M_+,M_-)$.

Conversely, given $B\in \mathbb{B}(M_+,M_-)$ such that $\|B\|< 1$, we can form the subspace $M=\{x+Bx|x\in M_+\}$. Note that
\begin{eqnarray*}-i[x+Bx,x+Bx]&=&(x,x)_+-(Bx,Bx)_-\\
&=&((Id-B^*B)x,x)_+\geq (1-\|B\|^2) (x,x)_+\\
&\geq& \frac{(1-\|B\|^2)}{2}\|x+Bx\|^2.\end{eqnarray*}
This implies that $M$ is completely positive-definite. If $M$ is not maximal, then there is a completely positive-definite subspace $M'$ containing $M$ properly. Let $y\in M'\setminus M$ and $y=y_++y_-$ according to the decomposition $H=M_+\oplus_\bot M_-$. Then \[y-(y_++By_+)=y_--By_+\in M'\cap M_-=0,\]
implying that $y=y_++By_+\in M$. A contradiction! This completes the proof.
\end{proof}

The theorem then gives rise to a complex (Banach) manifold structure on $W_+(H)$ and the open unit ball in $\mathbb{B}(M_+,M_-)$ serves as a \emph{global} coordinate chart. Similarly, all maximal completely negative-definite subspaces can be parameterized by elements in the open unit ball in $\mathbb{B}(M_-,M_+)$. Let $N_+$ be any maximal completely positive-definite subspace and define $N_-:=N_+^{\bot_s}$.
\begin{proposition}$N_-$ is a maximal completely negative-definite subspace of $H$ and $H=N_+\oplus N_-$.
\end{proposition}
\begin{proof}According to the above theorem, let $B\in \mathbb{B}(M_+,M_-)$ parameterize $N_+$. Then for any $y=y_++y_-\in N_-$, we have
\[0=-i[x,y]=(x_+,y_+)_+-(Bx_+, y_-)_-\]
for all $x_+\in M_+$. This implies that $y_+=B^*y_-$. Conversely, $B^*y_-+y_-$ certainly lies in $N_+^{\bot_s}$ for each $y_-\in M_-$. Thus $N_-=\{B^*y+y|y\in M_-\}$ and it is a maximal completely negative-definite subspace.

Obviously $N_+\cap N_-=0$. To see $N_+\oplus N_-=H$, we have to solve the equation
\[\left\{
\begin{array}{ll}
x_+=y_++B^*y_-,\\
x_-=y_-+By_+.
\end{array}
\right.\]
for any given $x=x_++x_-\in H=M_+\oplus_\bot M_-$. Note that the matrix of operators $\left(
                                                                 \begin{array}{cc}
                                                                   Id & B^* \\
                                                                   B & Id \\
                                                                 \end{array}
                                                               \right
)$ acting on $M_+\oplus M_-$ is invertible for $\|B\|=\|B^*\|< 1$. The solution $(y_+,y_-)$ exists and is unique.
\end{proof}
In Thm.~\ref{para}, we can replace $M_\pm$ by $N_\pm$ and a similar result holds except for that $H=N_+\oplus N_-$ is not orthogonal. Then $W_+(H)$ also has the open unit ball in $\mathbb{B}(N_+,N_-)$ as its coordinate chart. It's routine to check that  different charts are related to each other by biholomorphic maps (see the following Prop.~\ref{p8}) and the complex manifold structure on $W_+(H)$ is unique. In particular, let $(\cdot, \cdot)_+$ (resp. $(\cdot, \cdot)_-$) be the restriction of $-i[\cdot, \cdot]$ (resp. $i[\cdot, \cdot]$) on $N_+$ (resp. $N_-$) and $x=x_++x_-$  w.r.t. the decomposition $H=N_+\oplus N_-$. Then
\[-i[x,y]=(x_+,y_+)_+-(x_-,y_-)_-\]
for arbitrary $x,y\in H$.

 It is obvious from Thm.~\ref{para} that all maximal completely positive-definite (resp. negative-definite) subspaces have the same dimension $n_+$ (resp. $n_-$). We call the pair $(n_+, n_-)$ the signature of $[\cdot, \cdot]$. Actually the signature completely characterizes the strong symplectic structure $[\cdot, \cdot]$. We call two strong symplectic structures $[\cdot, \cdot]_1$ on $H_1$ and $[\cdot, \cdot]_2$ on $H_2$ are isomorphic if there is a linear topological isomorphism (not necessarily an isometry) $\Phi: H_1\rightarrow H_2$ such that \[[\Phi(x), \Phi(y)]_2=[x, y]_1.\] In particular, all automorphisms of a strong symplectic structure on $H$ form a group called pseudounitary group of signature $(n_+, n_-)$ and denoted by $\mathbb{U}(n_+, n_-)$ \footnote{The relevant strong symplectic Hilbert space will be clear from the background.}. $\mathbb{U}(n_+, n_-)$ is a subgroup of the general linear group $\mathbb{GL}(H)$, i.e., the group consisting of all invertible elements in $\mathbb{B}(H)$. $\mathbb{GL}(H)$ is open in $\mathbb{B}(H)$ and in this way is a Banach-Lie group. We always view $\mathbb{U}(n_+, n_-)$ as a real analytic Banach-Lie subgroup of $\mathbb{GL}(H)$.

\begin{proposition}All strong symplectic Hilbert spaces with signature $(n_+, n_-)$ are isomorphic.
\end{proposition}
\begin{proof}Given a strong symplectic Hilbert space $H$ with signature $(n_+,n_-)$, we choose $N_+\in W_+(H)$ and its symplectic complement $N_-$. Take an ordered orthonormal basis $\{e_j\}_{j=1}^{n_+}$ of $N_+$ w.r.t. the inner product $(\cdot, \cdot)_+$ and an ordered orthonormal basis $\{f_j\}_{j=1}^{n_-}$ of $N_-$ w.r.t. the inner product $(\cdot, \cdot)_-$. Then $N_+$ is identified with $\mathbb{C}^{n_+}$ and $N_-$ with $\mathbb{C}^{n_-}$ ($\mathbb{C}^{n_\pm}$ are equipped with the standard Hilbert space structure). Let $\mathbb{H}:=\mathbb{C}^{n_+}\oplus_\bot \mathbb{C}^{n_-}$ with the strong symplectic structure defined as in Example \ref{ex2}:
\[-i[x,y]=(x_+,y_+)_{\mathbb{C}^{n_+}}-(x_-,y_-)_{\mathbb{C}^{n_-}},\]
where $x=x_++x_-\in \mathbb{H}$. This identification is surely an isomorphism between the two strong symplectic structures.
\end{proof}

Given two isomorphic strong symplectic Hilbert spaces $H_1=N_+^1\oplus N_-^{1}$ and $H_2=N_+^2\oplus N_-^2$, an isomorphism $\Phi\in \mathbb{B}(H_1,H_2)$ can be written in the matrix form \[\left(
                                                                                                           \begin{array}{cc}
                                                                                                             g_{11} & g_{12} \\
                                                                                                             g_{21} & g_{22} \\
                                                                                                           \end{array}
                                                                                                         \right)
, g_{11}\in \mathbb{B}(N_+^1, N_+^2), g_{12}\in \mathbb{B}(N_-^1,N_+^2), g_{21}\in \mathbb{B}(N_+^1, N_-^2), g_{22}\in \mathbb{B}(N_-^1,N_-^2).\] These matrix entries should satisfy
\[g_{11}^*g_{11}-g_{21}^*g_{21}=Id,\quad g_{22}^*g_{22}-g_{12}^*g_{12}=Id,\quad g_{12}^*g_{11}-g_{22}^*g_{21}=0.\]
Certainly $\Phi$ transforms elements in $W_+(H_1)$ bijectively into elements in $W_+(H_2)$.
\begin{proposition}\label{p8}If $N_+\in W_+(H_1)$ is parameterized by $B\in \mathbb{B}(N_+^1, N_-^1)$ and $\Phi$ an isomorphism as above, then $\Phi(N_+)$ is parameterized by
\begin{equation}B'=(g_{21}+g_{22}B)(g_{11}+g_{12}B)^{-1}\in \mathbb{B}(N_+^2, N_-^2).\label{e6}\end{equation}
\end{proposition}
\begin{proof}For arbitrary $x_++Bx_+\in N_+$,
\[\Phi(x_++Bx_+)=\left(
                                                                                                           \begin{array}{cc}
                                                                                                             g_{11} & g_{12} \\
                                                                                                             g_{21} & g_{22} \\
                                                                                                           \end{array}
                                                                                                         \right)\left(
                                                                                                                  \begin{array}{c}
                                                                                                                    x_+ \\
                                                                                                                    Bx_+ \\
                                                                                                                  \end{array}
                                                                                                                \right)
                                                                                                        =\left(
                                                                                                           \begin{array}{c}
                                                                                                             (g_{11}+g_{12}B)x_+ \\
                                                                                                             (g_{21}+g_{22}B)x_+ \\
                                                                                                           \end{array}
                                                                                                         \right).                                                                                                 \]
Let $y_+:=(g_{11}+g_{12}B)x_+$. Then
\[\Phi(x_++Bx_+)=y_++(g_{21}+g_{22}B)(g_{11}+g_{12}B)^{-1}y_+.\]
The conclusion then follows.
\end{proof}
The result actually shows how different coordinate charts on $W_+(H)$ are related to each other. The formula also means that $\mathbb{U}(n_+,n_-)$ acts holomorphically on $W_+(H)$. Recall that a group $G$ acts on a set $M$ transitively if for any two elements $m_1, m_2\in M$, there is a $g\in G$ transforming $m_1$ into $m_2$. The isotropy group at $m\in M$ is the subgroup $G_m\subset G$ fixing $m$. Let $M_\pm$ be the subspaces as before.

\begin{proposition}
The action of $\mathbb{U}(n_+,n_-)$ on $W_+(H)$ is transitive and the isotropy subgroup at $M_+\in W_+(H)$ is $\mathbb{U}(M_+)\times \mathbb{U}(M_-)$.
\end{proposition}
\begin{proof}Let $B\in \mathbb{B}(M_+, M_-)$ parameterize $N_+$. It suffices to prove that there exists $\Phi\in \mathbb{U}(n_+,n_-)$ such that $\Phi(M_+)=N_+$. Define \[\Phi_0=\left(
                                                                 \begin{array}{cc}
                                                                   Id & B^* \\
                                                                   B & Id \\
                                                                 \end{array}
                                                               \right),\quad \Phi_1=\left(
                                                                                      \begin{array}{cc}
                                                                                        (Id-B^*B)^{-1/2} & 0 \\
                                                                                        0 & (Id-BB^*)^{-1/2} \\
                                                                                      \end{array}
                                                                                    \right)
\]
w.r.t. the decomposition $H=M_+\oplus_\bot M_-$. It can be checked directly that $\Phi:=\Phi_0\circ \Phi_1\in \mathbb{U}(n_+,n_-)$ and $\Phi(M_+)=N_+$.

Note that $M_+$ is parameterized by $0\in \mathbb{B}(M_+, M_-)$. If $\Phi\in\mathbb{U}(n_+,n_-)$ fixes $M_+$, then according to the above proposition, $g_{21}g_{11}^{-1}=0$ and hence $g_{21}=0$, $g_{12}=0$, $g_{11}^*g_{11}=Id$ and $g_{22}^*g_{22}=Id$. That's to say $\Phi=\left(
                                                                                                               \begin{array}{cc}
                                                                                                                 g_{11} & 0 \\
                                                                                                                 0 & g_{22} \\
                                                                                                               \end{array}
                                                                                                             \right)
$ where $g_{11}$ and $g_{22}$ are unitary operators on $M_+$ and $M_-$ respectively.
 \end{proof}

 Since $\mathbb{U}(M_\pm)$ are isomorphic to $\mathbb{U}(n_\pm)$ respectively, the above proposition means we can identify $W_+(H)$ with the homogeneous space $\mathbb{U}(n_+,n_-)/(\mathbb{U}(n_+)\times \mathbb{U}(n_-))$. Obviously, there is a normal subgroup of both $\mathbb{U}(n_+,n_-)$ and $\mathbb{U}(n_+)\times \mathbb{U}(n_-)$ consisting of elements like $c\times Id$ with $c\in \mathbb{U}(1)$. This normal subgroup acts trivially on $W_+(H)$. We denote the quotient groups by $\mathbb{PU}(n_+,n_-)$ and $\mathbb{P}(\mathbb{U}(n_+)\times \mathbb{U}(n_-))$. Still $\mathbb{PU}(n_+,n_-)$ acts on $W_+(H)$ effectively and transitively, with $\mathbb{P}(\mathbb{U}(n_+)\times \mathbb{U}(n_-))$ as the isotropy subgroup at $M_+$. That's to say $W_+(H)$ can also be identified with $\mathbb{PU}(n_+,n_-)/\mathbb{P}(\mathbb{U}(n_+)\times \mathbb{U}(n_-))$.

If $\textup{dim} H<+\infty$, the picture can be enlarged. Let $Gr(n_+, H)$ be the Grassmannian of subspaces of dimension $n_+$ in $H$. Then $Gr(n_+, H)$ is a complex projective manifold of dimension $n_+\cdot n_-$ while the previous $W_+(H)$ is an open subset of it. The projective general linear group $\mathbb{PGL}(H)$ acts transitively on $Gr(n_+, H)$. In this way $\mathbb{PU}(n_+,n_-)\subset \mathbb{PGL}(H)$ acts on $Gr(n_+, H)$ and $W_+(H)$ is just a $\mathbb{PU}(n_+,n_-)$-orbit. $W_+(H)$ is the famous non-compact irreducible Hermitian symmetric space of type $I_{n_+, n_-}$ in E. Cartan's classification of such spaces and $Gr(n_+, H)$ is its compact dual. The embedding $W_+(H)\hookrightarrow Gr(n_+, H)$ is a realization of the so-called Borel embedding and our Prop.~\ref{para} is actually a version of the so-called Harish-Chandra embedding, i.e., realizing $W_+(H)$ as a bounded symmetric domain. With its Bergman metric, $W_+(H)$ is a hyperbolic K$\ddot{a}$hler manifold. For more details about all these, see \cite{helgason1979differential,viviani2014tour,wolf1972fine}.

If $\textup{dim} H=+\infty$, many facts mentioned above still have their counterparts, but the situation is much more complicated for infinite-dimensional Grassmannians come into view. $W_+(H)$ still has a natural dual which can never be compact in any reasonable sense. See \cite{chu2020bounded} for further information on infinite-dimensional bounded symmetric spaces and \cite{abbondandolo2009infinite} for a topological investigation on various infinite-dimensional Grassmannians. In this paper we will only use the most obvious facts concerning $W_+(H)$ and its dual in the case $n_+=n_-=+\infty$ when necessary.

We are particularly interested in strong symplectic structures with signature $(n,n)$. Only in this setting, Lagrangian subspaces exist. For the remainder of this section, we always assume $n_+=n_-=n\leq +\infty$.
\begin{proposition}If the strong symplectic structure on $H$ has signature $(n, n)$, then all Lagrangian subspaces are parameterized by unitary maps from any $N_+\in W_+(H)$ to its symplectic complement $N_-$.
\end{proposition}
\begin{proof}This is easy and left to the interested reader.
\end{proof}

We denote the space of Lagrangian subspaces of $H$ by $\mathcal{L}(H)$. Since now $n_+=n_-=n$, we can identify $N_+$ with $N_-$ (using a unitary map between them) and then $\mathcal{L}(H)$ is parameterized by $\mathbb{U}(N_+)$. This gives rise to a real analytic Banach manifold structure on $\mathcal{L}(H)$. Different choices of $N_+$ only introduce different coordinate charts that are related to each other by analytic diffeomorphisms. If $H_1$ and $H_2$ are two strong symplectic Hilbert spaces and $\Phi$ is an isomorphism between them, then $\Phi$ transforms elements in $\mathcal{L}(H_1)$ bijectively into elements in $\mathcal{L}(H_2)$. In particular, $\mathbb{U}(n,n)$ acts on $\mathcal{L}(H)$.
\begin{proposition}Given $N_+^1\in W_+(H_1)$ and $N_+^2\in W_+(H_2)$. If $L\in \mathcal{L}(H_1)$ is parameterized by a unitary map $u$ from $N_+^1$ to $N_-^1$, then $\Phi(L)$ is parameterized by the unitary map
\[u'=(g_{21}+g_{22}u)(g_{11}+g_{12}u)^{-1}.\]
\end{proposition}
\begin{proof}The proof is similar to that of Prop.~\ref{p8}.
\end{proof}
There is another way to parameterize $W_+(H)$ in this case. If $N_-$ is identified with $N_+$ and $H$ with $N_+\oplus_\bot N_+$, then the strong symplectic structure becomes the standard one in Example \ref{ex2}. However, there is an alternative way to define a second strong symplectic structure on $N_+\oplus_\bot N_+$: If $x=(x_1,x_2)$ and $y=(y_1,y_2)$ are both in $N_+\oplus_\bot N_+$, then
\[[x,y]_{new}:=(x_2, y_1)_+-(x_1,y_2)_+\]
is also a strong symplectic structure with signature $(n,n)$. For $x=(x_1,x_2)$, we set $\beta(x)=(\frac{x_1-x_2}{i\sqrt{2}}, \frac{x_1+x_2}{\sqrt{2}})\in N_+\oplus_\bot N_+$. It's easy to check
\[[\beta(x),\beta(y)]_{new}=[x,y],\]
i.e., $\beta$ is a symplectic isomorphism between these two strong symplectic structures. Now if $N$ is a maximal completely positive-definite subspace w.r.t. $[\cdot, \cdot]$, then it is of the form $\{(x, Bx)|x\in N_+\}$ for some $B\in \mathbb{B}(N_+)$ with $\|B\|< 1$. We have
\begin{eqnarray*}\beta(N)&=&\{(\frac{x-Bx}{i\sqrt{2}}, \frac{x+Bx}{\sqrt{2}})\in N_+\oplus_\bot N_+|x\in N_+\}\\
&=&\{(x, i(Id+B)(Id-B)^{-1}x)\in N_+\oplus_\bot N_+|x\in N_+\}.\end{eqnarray*}
Thus the new operator $M=i(Id+B)(Id-B)^{-1}$ also parameterizes $N$. Note that the imaginary part of $M$
\[\frac{M-M^*}{2i}=(Id-B^*)^{-1}(Id-B^*B)(Id-B)^{-1}\]
is positive-definite and invertible. Conversely, any $M\in \mathbb{B}(N_+)$ with this property produces a $B=(M-i)(M+i)^{-1}$ with $\|B\|< 1$. In short, $B$ is the Cayley transform of $M$ while $M$ is the inverse transform of $B$.

The picture can be applied to parameterize Lagrangian subspaces as well. However, the inverse Cayley transform of a unitary operator $U$ is not necessarily an operator, but may have a multi-valued part. If $U$ parameterizes $L_U$, then
\[\beta(L_U)=\{(-i(Id-U)x, (Id+U)x)|x\in N_+\}.\]
Let $K$ be the orthogonal complement of $\textup{ker}(Id-U)$. Then $K$ is $U$-invariant and we set $U_K$ to be $U$ restricted on $K$. Let $A$ be the inverse Cayley transform of $U_K$. Then $A$ is a elf-adjoint operator in $K$. When $A$ is unbounded, its domain is $\textup{Ran}(Id-U)$. In terms of $A$,
\[\beta(L_U)=\{(x, Ax+y)|x\in D(A), y\in [D(A)]^\bot\}.\]
$\beta(L_U)$ is an example of (self-adjoint) linear relations, i.e., subspaces of $N_+\oplus_\bot N_+$.

\begin{proposition}$\mathbb{PU}(n,n)$ acts transitively on $\mathcal{L}(H)$.
\end{proposition}
\begin{proof}If $H$ is identified with $N_+\oplus_\bot N_+$ and $L\in \mathcal{L}(H)$ is parameterized by the unitary operator $U$. We only need to show there is a $\Phi\in \mathbb{U}(n,n)$ transforming $L$ into the Lagrangian subspace parameterized by the unitary operator $Id$. Taking $\Phi=\left(
                                                                                                       \begin{array}{cc}
                                                                                                         Id & 0 \\
                                                                                                         0 & U^{-1} \\
                                                                                                       \end{array}
                                                                                                     \right)
$ will do the job.
\end{proof}
\emph{Remark}. This just means $\mathcal{L}(H)$ is also a $\mathbb{PU}(n,n)$-orbit in $Gr(n, H)$.

For later convenience, we shall say something more about the complex manifold $Gr(n,H)$ and its relation to $W_\pm(H)$ and $\mathcal{L}(H)$. Obviously, $W_\pm(H)$ are two open subsets of $Gr(n, H)$ while $\mathcal{L}(H)$ is a closed subset of $Gr(n, H)$. This is in sharp contrast with the case $n_+\neq n_-$. In the latter case, $W_+(H)$ lies in $Gr(n_+, H)$, $W_-(H)$ lies in $Gr(n_-,H)$ and $\mathcal{L}(H)$ makes no sense. In particular, $Gr(n_+, H)$ and $Gr(n_-,H)$ are two different topological components of the Grassmannian of all closed subspaces of $H$. What's more important for us is the following.
\begin{proposition}\label{p11}Let $\overline{W_\pm(H)}$ be the topological closure of $W_\pm(H)$ in $Gr(n, H)$ respectively. Then $\overline{W_+(H)}\cap \overline{W_-(H)}=\mathcal{L}(H)$.
\end{proposition}
\begin{proof}Take $N_+\in W_+(H)$ and let $N_-\in W_-(H)$ be its symplectic complement. Then $H=N_+\oplus N_-$. This decomposition gives rise to two coordinate charts $U_\pm$ of $Gr(n,H)$:
\[U_+=\{S\in Gr(n, H)|S=\textup{graph}(B), B\in \mathbb{B}(N_+,N_-) \},\]
\[U_-=\{S\in Gr(n, H)|S=\textup{graph}(B), B\in \mathbb{B}(N_-,N_+) \},\]
and $W_\pm(H)\subset U_\pm$. Then $\mathcal{L}(H)\subset \overline{W_+(H)}\cap \overline{W_-(H)}$ is obvious. Now if $S\in \overline{W_+(H)}\cap \overline{W_-(H)}$, then it can be written in two ways, $\{(x, B_1x)\in N_+\oplus N_- |x\in N_+\}$ in $U_+$ and $\{(B_2y, y)\in N_+\oplus N_-|y\in N_-\}$ in $U_-$. The equality of the two implies that both $B_1$ and $B_2$ are invertible and $B_2=B_1^{-1}$. As a limit point of both $W_+(H)$ and $W_-(H)$ in $Gr(n, H)$, $S$ has to be both non-negative definite and non-positive definite. Therefore $Id-B_1^*B_1=0$ and consequently $B_1$ is a unitary map. This shows $S\in \mathcal{L}(H)$. The conclusion is thus proved.
\end{proof}
\emph{Remark}. If $n\in \mathbb{N}$, in the theory of Hermitian symmetric spaces, $\mathcal{L}(H)$ is the Shilov boundary of $W_+(H)$ in $Gr(n, H)$, which plays a fundamental role in the function theory on $W_+(H)$.

The following example is to demonstrate the above argument in the simplest case.
\begin{example} Let $H=\mathbb{C}\oplus_\bot \mathbb{C}=\mathbb{C}^2$ and the strong symplectic structure be given by
$-i[x,y]=x_1\overline{y_1}-x_2\overline{y_2}$,
where $x=(x_1,x_2)$ and $y=(y_1,y_2)$ are in $\mathbb{C}^2$. Then each maximal completely positive-definite subspace is of the form $\mathbb{C}\{(1,z)\}$ with $|z|<1$, each maximal completely negative-definite subspace is of the form $\mathbb{C}\{(z,1)\}$ for $|z|<1$, and each Lagrangian subspace is of the form $\mathbb{C}\{(1,z)\}$ with $|z|=1$. Note that in this case $W_+(H)\cup W_-(H)\cup \mathcal{L}(H)=Gr(1,2)=\mathbb{CP}^1$ and $\mathcal{L}(H)$ is the common boundary of $W_+(H)$ and $W_-(H)$ in $\mathbb{CP}^1$. It is easy to find that $\mathbb{PU}(1,1)$ is precisely the M$\ddot{o}$bius transformation group of the unit disc. In terms of the inverse Cayley transform, the unit disc is transformed into the upper half plane $\mathbb{C}_+$ and $\mathbb{PU}(1,1)$ is simply described as $\mathbb{PSL}(2,\mathbb{R})=\mathbb{SL}(2,\mathbb{R})/\{\pm Id\}$. If $g\in\mathbb{PSL}(2,\mathbb{R})$ is represented by $\pm\left(
                                                                                                                                                   \begin{array}{cc}
                                                                                                                                                     a & b \\
                                                                                                                                                     c & d \\
                                                                                                                                                   \end{array}
                                                                                                                                                 \right)
$ and $z\in \mathbb{C}_+$, then $g\cdot z=\frac{a z+b}{c z+d}$.
\end{example}

The following proposition will be used in \S~\ref{sec9}.
\begin{proposition}\label{p3}If $N$ is a maximal completely positive-definite subspace of $H$ and $L$ is Lagrangian, then $N$ and $L$ are transversal in $H$, i.e., $N\cap L=0$, and $N\oplus L=H$.
\end{proposition}
\begin{proof}$L\cap N=0$ is obvious. We can identify $H$ with $N\oplus_\bot N$. Then $N$ is identified with $\mathcal{N}:=\{(x,0)\in  N\oplus_\bot N|x\in N\}$. Let $U\in \mathbb{U}(N)$ be the operator parameterizing $L$, i.e.,
\[L=\{(x,Ux)\in N\oplus_\bot N|x\in N \}.\]
For any $(x_+, x_-)\in N\oplus_\bot N$, we have
\[(x_+,x_-)=(x_+-U^{-1}x_-,0)+(U^{-1}x_-,UU^{-1}x_-).\]
Thus $N\oplus_\bot N=N+L$.
\end{proof}

\section{Symmetric operator and its Weyl curve}\label{sec4}
In this paper, we are basically concerned with simple symmetric operators.
\begin{definition}A densely defined symmetric operator $T$ in a Hilbert space $H$ is called simple, if it is not a nontrivial orthogonal sum of a self-adjoint operator and a symmetric operator.
\end{definition}
The classification of self-adjoint operators up to unitary equivalence is well-known and described by the Hahn-Hellinger theory of spectral multiplicity. Hence to classify symmetric operators is essentially to classify the simple ones. It is known that simplicity of $T$ is equivalent to that $\textup{span}\cup_{\lambda\in \mathbb{C}_+\cup \mathbb{C}_-}\{\textup{ker}(T^*-\lambda)\}$ is dense in $H$.

Let $T$ be a simple symmetric operator in $H$ and $T^*$ its adjoint. Note that $D(T^*)$ is a Hilbert space with its graph inner product and $D(T)$ is a closed subspace of it. Then the quotient $\mathcal{B}_T:=D(T^*)/D(T)$ is also a Hilbert space. It is not hard to see that for $[x],[y]\in \mathcal{B}_T$
\[[[x],[y]]_T:=(T^*x,y)-(x, T^*y)\]
is well-defined and a strong symplectic structure on $\mathcal{B}_T$. Let $(n_+, n_-)$ be the signature of this structure. $n_\pm$ are just the deficiency indices. As before we only consider the case $n_+, n_-\geq 1$. It's a basic fact that closed extensions of $T$ correspond to closed subspaces of $\mathcal{B}_T$; in particular, self-adjoint extensions (if any) correspond to Lagrangian subspaces of $\mathcal{B}_T$ \cite[Prop.~14.7]{schmudgen2012unbounded}.

Let $W_+(\mathcal{B}_T)$ (resp. $W_-(\mathcal{B}_T)$) be the space of maximal completely positive-definite (resp. negative-definite) linear subspaces of $\mathcal{B}_T$ w.r.t. $[\cdot, \cdot]_T$. The following proposition is essentially \cite[Prop.~1.6.8]{behrndt2020boundary}, but we prefer a geometric reformulation.
\begin{proposition}\label{p6}
For $\lambda\in \mathbb{C}_+$ (resp. $\mathbb{C}_-$), the image $W(\lambda)$ of $\textup{ker}(T^*-\lambda)\subset D(T^*)$ in $\mathcal{B}_T$ is a maximal completely positive-definite (resp. negative-definite) subspace.
\end{proposition}
\begin{proof} If $\lambda\in \mathbb{C}_+$ and $0\neq x\in \textup{ker}(T^*-\lambda)$, then $\textup{ker}(T^*-\lambda)\cap D(T)=0$ due to the well-known fact that $T$ has no eigenvalues on $\mathbb{C}_+$. On the other side,
\[-i[[x],[x]]_T=-i[(T^*x,x)-(x,T^*x)]=-i(\lambda-\bar{\lambda})(x,x)=2\Im\lambda(x,x)>0.\]

According to the von Neumann decomposition \cite[Thm.~1.7.11]{behrndt2020boundary} \[D(T^*)=D(T)\oplus \textup{ker}(T^*-\lambda) \oplus \textup{ker}(T^*-\bar{\lambda}),\]
we have $\mathcal{B}_T=W(\lambda)\oplus W(\bar{\lambda})$. Note that $\textup{ker}(T^*-\lambda)$ is closed both in $D(T^*)$ and in $H$ and on it the two norms are equivalent. This shows that $W(\lambda)$ is completely positive-definite in $\mathcal{B}_T$. Similarly, $W(\bar{\lambda})$ is completely negative-definite. If $x\in \textup{ker}(T^*-\lambda)$ and $y\in \textup{ker}(T^*-\bar{\lambda})$, then
\[-i[[x],[y]]_T=-i[(T^*x,y)-(x,T^*y)]=-i(\lambda-\lambda)(x,y)=0,\]
implying that $W(\bar{\lambda})\subset (W(\lambda))^{\bot_s}$. These arguments are sufficient to show $W(\lambda)$ and $W(\bar{\lambda})$ are maximal.
\end{proof}
This basic proposition then gives rise to a map $$W_T:\mathbb{C}_+\cup \mathbb{C}_-\rightarrow W_+(\mathcal{B}_T)\cup W_-(\mathcal{B}_T).$$
 We call $W_T$ the two-branched Weyl map of $T$. When we want to emphasize the geometric content, we often call $W_T(\lambda)$ the Weyl curve of $T$. Since the component from $\mathbb{C}_-$ to $W_-(\mathcal{B}_T)$ is completely determined by the component from $\mathbb{C}_+$ to $W_+(\mathcal{B}_T)$, usually we can focus on the latter and call it the single-branched Weyl map (curve).
\begin{theorem}\label{thm11}The single-branched Weyl map $W_T(\lambda)$ is holomorphic.
\end{theorem}
This theorem is the most basic fact in the whole theory we are developing, because several definitions and constructions to follow actually are based on it.

Recall that a map between two complex manifolds is by definition holomorphic if and only if its local form in coordinate charts is given by holomorphic functions. To prove this theorem, we should introduce more concepts. If $\mathbb{H}$ is the standard strong symplectic Hilbert space in Example \ref{ex2} and $\Phi$ is a symplectic isomorphism between $\mathcal{B}_T$ and $\mathbb{H}$, then $\Phi$ can be extended to $D(T^*)$ simply by setting its value at $x\in D(T^*)$ to be $\Phi([x])$. By abuse of notation, we still use $\Phi$ to denote this extension. According to the decomposition $\mathbb{H}=H_+\oplus_\bot H_-$, $$\Phi(x)=(\Gamma_+x,\Gamma_-x)\in H_+\oplus_\bot H_-$$ and the abstract Green's second formula
\begin{equation}(T^*x, y)-(x,T^*y)=i(\Gamma_+x,\Gamma_+ y)_{H_+}-i(\Gamma_-x,\Gamma_-y)_{H_-}\end{equation}
holds for any $x,y\in D(T^*)$. By our discussion in the previous section, in terms of the above decomposition of $\mathbb{H}$, $\Phi(W_T(\lambda))$ is parameterized by $B(\lambda)\in \mathbb{B}(H_+,H_-)$ with $\|B(\lambda)\|<1$. It suffices to prove that $B(\lambda)$ is holomorphic on $\mathbb{C}_+$.

Consider the closed extension $T_+$ of $T$ whose domain is $\{x\in D(T^*)|\Gamma_+x=0\}$ and the closed extension $T_-$ whose domain is $\{x\in D(T^*)|\Gamma_-x=0\}$.
\begin{lemma}$T_+^*=T_-$.
\end{lemma}
\begin{proof}If $y\in D(T_-)$, then for any $x\in D(T_+)$
\[(T_+x, y)=i(\Gamma_+x,\Gamma_+y)-i(\Gamma_-x,\Gamma_-y)+(x, T^*y)=(x, T^*y)=(x, T_-y),\]
implying by definition that $y\in D(T_+^*)$ and $T_+^*y=T_-y$ and hence $T_-\subset T_+^*$. Conversely, if $y\in D(T^*_+)$, then for any $x\in D(T_+)$
\begin{eqnarray*}0&=&(T_+x,y)-(x, T_+^*y)=(T^*x,y)-(x, T^*y)\\
&=&i(\Gamma_+x,\Gamma_+y)_{H_+}-i(\Gamma_-x,\Gamma_-y)_{H_-}\\
&=&-i(\Gamma_-x,\Gamma_-y)_{H_-}.\end{eqnarray*}
It is easy to see $\Gamma_-(D(T_+))=H_-$ and thus $\Gamma_-y=0$, i.e., $y\in D(T_-)$. The proof is finished.
\end{proof}
\begin{proposition}\label{res}$\mathbb{C}_+\subset \rho(T_+)$ and $\mathbb{C}_-\subset \rho(T_-)$.
\end{proposition}
\begin{proof}For $\lambda\in \mathbb{C}_+$, we have to prove $T_+-\lambda$ has a bounded inverse. Let $\lambda=a+bi$ for $a,b\in \mathbb{R}$. Then for $x\in D(T_+)$
\begin{eqnarray*}\|(T_+-a-bi)x\|^2&=&((T_+-a-bi)x,(T_+-a-bi)x)\\
&=&((T_+-a)x, (T_+-a)x)+b^2\|x\|^2+ib[(T_+x,x)-(x,T_+x)]\\
&=&((T_+-a)x, T_+-a)x)+b^2\|x\|^2\\&-&b[(\Gamma_+x,\Gamma_+x)_{H_+}-(\Gamma_-x,\Gamma_-x)_{H_-}]\\
&=&((T_+-a)x, (T_+-a)x)+b^2\|x\|^2+b(\Gamma_-x,\Gamma_-x)_{H_-}\\
&\geq& b^2\|x\|^2.
\end{eqnarray*}
That's to say $\|(T_+-\lambda)x\|\geq b\|x\|$ for $x\in D(T_+)$. Since $T_+-\lambda$ is a closed operator, this implies that $\textup{Ran}(T_+-\lambda)$ is closed and $(T_+-\lambda)^{-1}$ is well-defined and bounded. Thus $\lambda\in \rho(T_+)\cup \sigma_r(T_+)$.  A similar argument applying to $T_-$ shows that $\bar{\lambda}\in \rho(T_-)\cup \sigma_r(T_-)$.

If $\lambda\in \sigma_r(T_+)$, then $\textup{Ran}(T_+-\lambda)$ is a proper closed subspace of $H$. Choose $y\neq 0$ in the orthogonal complement of $\textup{Ran}(T_+-\lambda)$ in $H$. Then for any $x\in D(T_+)$,
\[0=((T_+-\lambda)x, y)=(x, (T_--\bar{\lambda})y).\]
Since $D(T_+)$ is dense in $H$, we have $(T_--\bar{\lambda})y=0$, i.e., $\bar{\lambda}$ is an eigenvalue of $T_-$. A contradiction! Therefore, we must have $\lambda\in \rho(T_+)$. The conclusion for $T_-$ follows similarly.
\end{proof}
\begin{proposition}\label{p1} For any $\varphi\in H_+$ and $\lambda\in \mathbb{C}_+$, the abstract boundary value problem \[\left\{
\begin{array}{ll}
T^*x=\lambda x,\\
\Gamma_+ x=\varphi.
\end{array}
\right.\]
has a unique solution $x\in D(T^*)$.
\end{proposition}
\begin{proof}If both $x,y$ are solutions of the problem, then $\Gamma_+(x-y)=0$ and thus $x-y\in D(T_+)$. $x-y$ also fulfills the equation $T_+(x-y)=T^*(x-y)=\lambda(x-y)$. Thus $x-y=0$ for $\lambda$ is a regular value of $T_+$. The uniqueness is proved.

 This $x$ does exist. It is because $W(\lambda)\in W_+(\mathcal{B}_T)$ and $\Phi(W(\lambda))$ should again be maximally completely positive-definite in $\mathbb{H}$ and the graph of a $B(\lambda)\in \mathbb{B}(H_+,H_-)$.
\end{proof}

Thus from this proposition, we can define an operator-valued function $\gamma_+:\mathbb{C}_+\rightarrow \mathbb{B}(H_+,H)$ by letting $\gamma_+(\lambda)\varphi$ be the unique solution of the equation in Prop.~\ref{p1}. Then in the coordinate chart induced from $\Phi$, the parameter corresponding to $W_T(\lambda)$ is $B(\lambda)=\Gamma_-\gamma_+(\lambda)$. The following proposition and its corollary then complete the proof of Thm.~\ref{thm11}.
\begin{proposition}\label{p2}For $\lambda, \lambda'\in \mathbb{C}_+$,
\[\gamma_+(\lambda)-\gamma_+(\lambda')=(\lambda-\lambda')(T_+-\lambda)^{-1}\gamma_+(\lambda').\]
\end{proposition}
\begin{proof}Since by definition \[\Gamma_+(\gamma_+(\lambda)\varphi-\gamma_+(\lambda')\varphi)=\varphi-\varphi=0,\]  $\gamma_+(\lambda)\varphi-\gamma_+(\lambda')\varphi\in D(T_+)$. Note that
\begin{eqnarray*}T_+(\gamma_+(\lambda)\varphi-\gamma_+(\lambda')\varphi)&=&T^*(\gamma_+(\lambda)\varphi-\gamma_+(\lambda')\varphi)\\
&=&\lambda\gamma_+(\lambda)\varphi-\lambda'\gamma_+(\lambda')\varphi\\
&=& \lambda(\gamma_+(\lambda)\varphi-\gamma_+(\lambda')\varphi)+(\lambda-\lambda')\gamma_+(\lambda')\varphi.
\end{eqnarray*}
Thus,
\[\gamma_+(\lambda)\varphi-\gamma_+(\lambda')\varphi=(\lambda-\lambda')(T_+-\lambda)^{-1}\gamma_+(\lambda')\varphi\]
as required.
\end{proof}
We obviously have

\begin{corollary}\label{c1}For $\lambda, \lambda'\in \mathbb{C}_+$,
\[B(\lambda)-B(\lambda')=(\lambda-\lambda')\Gamma_-(T_+-\lambda)^{-1}\gamma_+(\lambda').\]
\end{corollary}
It should be mentioned that in terms of $T_-$, a map $\gamma_-(\lambda)\in \mathbb{B}(H_-,H)$ can be defined similarly for $\lambda\in \mathbb{C}_-$ and $\Gamma_+\gamma_-(\lambda)=(B(\bar{\lambda}))^*$. This fact will be freely used in later sections.

We call the analytic operator-valued function $B(\lambda)$ the \emph{contractive Weyl function} associated to the isomorphism $\Phi$. Here "contractive" simply means $\|B(\lambda)\|< 1$ for $\lambda\in \mathbb{C}_+$. It is clear that $B(\lambda)$ is only a representation of $W_T(\lambda)$ in terms of a specific choice of the "coordinate system" $\Phi$. If $\mathbb{H}'=H_+'\oplus_\bot H_-'$ is another standard strong symplectic Hilbert space and $\mathbb{H}$, $\mathbb{H}'$ are related by a symplectic isomorphism $\Phi=\left(
                                                                                                                                                        \begin{array}{cc}
                                                                                                                                                          g_{11} & g_{12}\\
                                                                                                                                                          g_{21} & g_{22}\\
                                                                                                                                                        \end{array}
                                                                                                                                                      \right)
$, then according to Prop.~\ref{p8}, in terms of $\mathbb{H}'$, the new contractive Weyl function is
\[B'(\lambda)=(g_{21}+g_{22}B(\lambda))(g_{11}+g_{12}B(\lambda))^{-1}.\]

If $n_+=n_-=n\leq +\infty$ and $G$ is a Hilbert space of dimension $n$, due to the discussion in \S~\ref{sec3} we can use the two standard strong symplectic structures on $G\oplus_\bot G$ to describe the Weyl curve: if $B(\lambda)\in \mathbb{B}(G)$ is the contractive Weyl function in terms of an isomorphism from $\mathcal{B}_T$ to $(G\oplus_\bot G, [\cdot, \cdot])$, then in terms of $(G\oplus_\bot G, [\cdot, \cdot]_{new})$, the curve can also be described as the inverse Cayley transform of $B(\lambda)$ i.e.,
\[M_+(\lambda)=i(Id+B(\lambda))(Id-B(\lambda))^{-1}\quad \textup{for}\quad \lambda\in \mathbb{C}_+.\]
It is also easy to find that for $\lambda\in \mathbb{C}_-$, the curve is described by $M_-(\lambda)=M_+(\bar{\lambda})^*$. The analytic operator-valued function
\begin{equation*}M(\lambda):=\left \{
\begin{array}{ll}M_+(\lambda),\quad \lambda\in \mathbb{C}_+\\ M_-(\lambda),\quad \lambda\in \mathbb{C}_-.
\end{array}
 \right.
 \end{equation*}
is called the Weyl function associated to the corresponding isomorphism. However, traditionally this Weyl function arises in another way.
The following definition is standard and only applies to the case $n_+=n_-$.
\begin{definition}Given a densely defined symmetric operator $T$ in $H$, a boundary triplet for $T$ is a triplet $(G, \Gamma_0, \Gamma_1)$ consisting of a Hilbert space $G$ and two linear maps $\Gamma_0, \Gamma_1: D(T^*)\rightarrow G$ such that\\
(1) The map $(\Gamma_0, \Gamma_1):D(T^*)\rightarrow G\times G$ is surjective;\\
(2) For any $x, y\in D(T^*)$, the following abstract Green's second formula holds:
\[(T^*x, y)-(x,T^*y)=(\Gamma_1x, \Gamma_0y)_G-(\Gamma_0 x, \Gamma_1y)_G.\]
\end{definition}
In our viewpoint, a boundary triplet is nothing else but essentially an isomorphism between the strong symplectic Hilbert spaces $\mathcal{B}_T$ and $(G\oplus_\bot G, [\cdot, \cdot]_{new})$. We can set $\Gamma_+=(\Gamma_1+i\Gamma_0)/\sqrt{2}$ and $\Gamma_-=(\Gamma_1-i\Gamma_0)/\sqrt{2}$. Therefore, Green's formula now reads
\[(T^*x, y)-(x,T^*y)=i(\Gamma_+x, \Gamma_+y)_G-i(\Gamma_- x, \Gamma_-y)_G.\]
This in essence gives an isomorphism between $\mathcal{B}_T$ and $(G\oplus_\bot G, [\cdot, \cdot])$. By abuse of language, we will call both $(G, \Gamma_0, \Gamma_1)$ and $(G, \Gamma_\pm)$ a boundary triplet. The context will imply which we exactly mean.

It is known that self-adjoint extensions of $T$ are parameterized by points in $\mathcal{L}(\mathcal{B}_T)$. So, given a boundary triplet $(G, \Gamma_\pm)$, all self-adjoint extensions are characterized by elements in $\mathbb{U}(G)$, i.e., the group of unitary operators in $G$. Actually, the domain of the self-adjoint extension $T_{U}$ associated with $U\in \mathbb{U}(G)$ is
\begin{equation}D(T_U)=\{x\in D(T^*)|\Gamma_-x=U\Gamma_+x\}.\label{main1}\end{equation}
In particular, $U=Id$ represents the self-adjoint extension corresponding to the abstract boundary condition $\Gamma_0x=0$ while $U=-Id$ represents that corresponding to the boundary condition $\Gamma_1x=0$. For convenience, if the boundary triplet is understood, we denote these two by $T_0$ and $T_1$ respectively.

Fix a boundary triplet $(G, \Gamma_0, \Gamma_1)$ for $T$. For $\varphi\in G$ and $\lambda\in \rho(T_0)$, the abstract boundary value problem
\[\left\{
\begin{array}{ll}
T^*x=\lambda x,\\
\Gamma_0x=\varphi.
\end{array}
\right.\]
has a unique solution $\gamma(\lambda) \varphi$. This map $\gamma(\lambda):G\rightarrow \textup{ker}(T^*-\lambda)\subset H$ is called the $\gamma$-field associated to the boundary triplet and the map $M(\lambda):=\Gamma_1\circ\gamma(\lambda)\in \mathbb{B}(G)$ is called the Weyl function associated to the boundary triplet. Obviously, this Weyl function coincides with our previous one.

There is also a more conceptual and intuitive interpretation of the rise of Weyl curves. We explain it as follows. The direct sum $H\oplus_\bot H$ is equipped with the second standard strong symplectic structure: for $x=(x_1,x_2)\in H\oplus_\bot H$, $y=(y_1,y_2)\in H\oplus_\bot H$, \[[x, y]_{new}=(x_2,y_1)_H-(x_1,y_2)_H.\]
 Then a simple symmetric operator $T$ induces an isotropic subspace (the graph of $T$) \[A_T:=\{(x, Tx)\in H\oplus_\bot H|x\in D(T)\}\] and its symplectic complement is a co-isotropic subspace (the graph of $T^*$)  \[A_T^{\bot_s}:=\{(x, T^*x)\in H\oplus_\bot H|x\in D(T^*)\}.\]
 The quotient $A_T^{\bot_s}/A_T$ of these two is essentially $D(T^*)/D(T)$ and acquires a reduced strong symplectic structure from $H\oplus_\bot H$.

 Recall that there is a similar reduction procedure in real symplectic geometry, see for example \cite{marsden1974reduction} \cite[Lemma~2.1.7]{mcduff2017introduction} and our \S~\ref{real}. It is also well-known that several structures (e.g., Lagrangian subspaces, compatible complex structures) on the original linear symplectic space have their reduced version on the quotient linear symplectic space.

In our setting, the story is indeed similar. In $H\oplus_\bot H$, there is a holomorphic curve $\mathcal{W}(\lambda)$ of subspaces: $\lambda\in \mathbb{C}\mapsto \{(x, \lambda x)\in H\oplus_\bot H|x\in H\}$. Note that for $\lambda\in \mathbb{C}_+$, $\mathcal{W}(\lambda)$ is maximally completely positive-definite in $H\oplus_\bot H$ and that for $\lambda\in \mathbb{C}_-$, $\mathcal{W}(\lambda)$ is maximally completely negative-definite. As for $\lambda\in \mathbb{R}$, $\mathcal{W}(\lambda)$ is Lagrangian. The reduced version of $\mathcal{W}(\lambda)$ is $(\mathcal{W}(\lambda)\cap A_T^{\bot_s}+A_T)/A_T$. Since for $\lambda\in \mathbb{C}_+\cup \mathbb{C}_-$
 \[\mathcal{W}(\lambda)\cap A_T^{\bot_s}=\{(x, \lambda x)|x\in D(T^*),\, T^*x=\lambda x\}\]
 which is essentially $\textup{ker}(T^*-\lambda)$,
 so $W_T(\lambda)$ is precisely the reduced version of $\mathcal{W}(\lambda)$ and the maximal complete positive-definiteness (and negative-definiteness) survives after the reduction procedure. In this sense, we call $\mathcal{W}(\lambda)$ the \emph{universal} Weyl curve associated to $H$, and any Weyl curve can be derived from it through this reduction procedure. One should note that singularities along the real line $\mathbb{R}\subset \mathbb{C}$ may arise in this procedure and $W_T(\lambda)$ needs not to be holomorphic at $\lambda\in \mathbb{R}$.

 The advantage of this briefly sketched formalism is that it also applies to simple symmetric operators which are not densely defined. In this case, $T^*$ doesn't exist as an operator, but $A_T^{\bot_s}$ as the symplectic complement of $A_T$ still makes sense and is actually a linear relation. In this paper, we also deal with simple symmetric operators that are not densely defined, but we often pretend that $T$ is densely defined and $T^*$ is well-defined. This is only for convenience of a direct presentation and the gap can be easily filled in by using terminology and notation developed for linear relations in \cite[Chap.~1]{behrndt2020boundary}. We have collected from \cite{behrndt2020boundary} some elementary facts on linear relations in Appendix B, which provides the necessary material to fill in the gap.

\section{Symmetric operator and its characteristic vector bundles}\label{sec5}
We want to compare different simple symmetric operators with the same deficiency indices $(n_+,n_-)$. So it is important to compare their Weyl curves. However, conceptually they live in different spaces. The solution is to use the same standard strong symplectic Hilbert space $\mathbb{H}$ and to consider the corresponding curves $B(\lambda)$ in $W_+(\mathbb{H})$ as Weyl curves. But even if $\mathbb{H}$ is fixed, the choice of $\Phi$ may still differ by the action of elements in $\mathbb{PU}(n_+,n_-)$. This ambiguity leads us to the following definitions.
\begin{definition}Provided that the standard strong symplectic Hilbert space $\mathbb{H}$ with signature $(n_+,n_-)$ is fixed. A holomorphic map $N(\lambda)$ from $\mathbb{C}_+$ to $W_+(\mathbb{H})$ is called a Nevanlinna curve with genus $(n_+,n_-)$.
\end{definition}

\begin{definition} If two Nevanlinna curves in $W_+(\mathbb{H})$ differ only by the action of an element in $\mathbb{PU}(n_+,n_-)$, we say they are congruent. This is an equivalence relation and each equivalence class is called a congruence class. We call the congruence class determined by $W_T(\lambda)$ the Weyl class of $T$, and denote it by $[W_T(\lambda)]$.
\end{definition}
Recall that two symmetric operators $T$ in $H$ and $T'$ in $H'$ are unitarily equivalent if there is a unitary map $U$ from $H$ to $H'$ such that $UD(T)=D(T')$ and $T'Ux=UTx$ for any $x\in D(T)$. If two simple symmetric operators $T_1$ and $T_2$ with deficiency indices $(n_+, n_-)$ are in the same unitary equivalence class, then it can be easily found that $[W_{T_1}(\lambda)]=[W_{T_2}(\lambda)]$. This shows that the Weyl class $[W_T(\lambda)]$ is a unitary invariant of $T$.  One of our goals in this paper is to show
\begin{theorem}\label{thm12}
 $[W_T(\lambda)]$ is a \emph{complete} unitary invariant of a (not necessarily densely defined) simple symmetric operator $T$, i.e., $[W_T(\lambda)]$ uniquely determines the unitary equivalence class of $T$. Moreover, there is a one-to-one correspondence between unitary equivalence classes of simple symmetric operators with deficiency indices $(n_+,n_-)$ and congruence classes of Nevanlinna curves with genus $(n_+,n_-)$.
\end{theorem}

In fact, this basic theorem has already been established before, at least in the case $n_+=n_-$, see for example \cite[Chap.~4]{behrndt2020boundary}. However, it was never stated in this geometric fashion and the existing proofs are not \emph{conceptually} transparent. In this section, we shall prove this theorem in its total generality and in an \emph{intrinsic} manner. We insist that this theorem should be viewed as the most basic result towards the classification of simple symmetric operators and its importance lies partly in the fact that it connects spectral theory with complex analysis and complex geometry in a fundamental way.

To prove Thm.~\ref{thm12}, we should introduce more concepts, which obviously have their own interest.
\subsection{Characteristic vector bundle of the first kind}\label{ssec1}
Recall that a complex vector bundle $E$ over a complex manifold $M$ is holomorphic, if the local transition functions of $E$ can be chosen to be holomorphic functions defined on overlaps of local coordinate charts in $M$. The basics of vector bundles and connections can be found in \cite[Chap.~0]{griffiths2014principles}.

Now it may be not so strange that the theory of holomorphic vector bundles enters into the theory of symmetric operators naturally. Given a simple symmetric operator $T$ in $H$, at each point $\lambda\in \mathbb{C_+}\cup \mathbb{C}_-$ ($\mathbb{C_+}\cup \mathbb{C}_-$ is basically a complex manifold of dimension 1, with a specified coordinate chart), we can attach the linear vector space $E_\lambda:=\textup{ker}(T^*-\lambda)$, whose rank is $n_+$ or $n_-$. The inner product on $H$ can be restricted on each fiber $E_\lambda$. Thus we actually have at hand a Hermitian vector bundle $E(T)$ over $\mathbb{C_+}\cup \mathbb{C}_-$. Given a standard strong symplectic Hilbert space $\mathbb{H}$ which is isomorphic to $\mathcal{B}_T$ via $\Phi$. Then $\gamma_\pm(\lambda)$ give rise to a bundle isomorphism between $E$ and the trivial vector bundle $(\mathbb{C}_+\times H_+)\cup (\mathbb{C}_-\times H_-)$. The trivial holomorphic structure on the latter can thus be transported to $E(T)$. A different choice of $\Phi$ won't change the holomorphic structure of $E$--This is precisely what the transformation formula (\ref{e6}) tells us. Consequently, the holomorphic structure on $E(T)$ is actually \emph{intrinsic} and determined by $T$ itself.

\begin{definition}The holomorphic Hermitian vector bundle $E(T)$ over $\mathbb{C}_+\cup \mathbb{C}_-$ is called the characteristic vector bundle of the first kind for $T$.
\end{definition}
Two holomorphic Hermitian vector bundles $E$ and $E'$ over the same base manifold $M$ are called isomorphic if there is a holomorphic bundle isomorphism which is an isometry between the fibers. Recall that a holomorphic bundle isomorphism between $E$ and $E'$ is a holomorphic map $\zeta: E\longrightarrow E'$ covering the identity map on the base manifold and being an isomorphism between the fibers.

\begin{theorem}\label{thm2}Two simple symmetric operators $T$ and $T'$ are unitarily equivalent, if and only if their characteristic bundle $E(T)$ and $E(T')$ are isomorphic as holomorphic Hermitian vector bundles.
\end{theorem}
\begin{proof}The necessity is obvious. The sufficiency part will be clear after the discussion in the next subsection.
\end{proof}

A basic quantity of a holomorphic Hermitian vector bundle is its curvature. Recall that a connection in a vector bundle $E$ over a smooth manifold $M$ is a $\mathbb{C}$-bilinear map
\[\nabla : \mathcal{X}(M)\times \Gamma(E)\rightarrow \Gamma(E),\]
where $\mathcal{X}(M)$ is the space of smooth vector fields on $M$ (or smooth sections of the tangent bundle $TM$) and $\Gamma(E)$ is the space of smooth sections of $E$. $\nabla$ should satisfy additionally the following conditions: For $f\in C^\infty(M)$, $X\in \mathcal{X}(M)$ and $s\in \Gamma(E)$,
\[\nabla_X(fs)=f\nabla_Xs+(Xf)s\]
and
\[\nabla_{fX}s=f\nabla_Xs,\]
where $Xf$ is the Lie derivative of $f$ along $X$. The curvature of $\nabla$ is a 2-form valued in the bundle $\textup{Hom}(E)$, and defined by
\[R(X,Y)s=\nabla_X\nabla_Ys-\nabla_Y\nabla_Xs-\nabla_{[X,Y]_L}s\]
for $X,Y\in \mathcal{X}(M)$ and $s\in \Gamma(E)$, where $[X,Y]_L$ is the Lie bracket of $X$ and $Y$. The curvature is of great importance because it carries topological information of the bundle $E$. This can be extracted through the famous Chern-Weil theory.

Let $E$ be equipped with a Hermitian metric $h$. A connection $\nabla$ is said to be compatible with $h$, if for $s_1,s_2\in \Gamma(E)$
\[dh(s_1,s_2)=h(\nabla s_1, s_2)+h(s_1,\nabla s_2),\]
where $d$ is the exterior differential operator. For a holomorphic vector bundle $E$ over a complex manifold $M$, a connection $\nabla$ can be decomposed into $\partial+\bar{\partial}$ according to the fact that a 1-form can be decomposed into its (1,0)-part and (0,1)-part. $\nabla$ is said to be compatible with the holomorphic structure if the $(0,1)$-part of $\nabla$ coincides with the Cauchy-Riemann operator $\bar{\partial}$ associated to the holomorphic structure.\footnote{Locally the kernel of $\bar{\partial}$ just tells us what should be considered as holomorphic sections of $E$.} A basic fact is that given a Hermitian structure $h$, there is a unique connection that is compatible with both $h$ and the holomorphic structure on $E$. It is called the Chern connection.

Let us sketch briefly how the Chern connection can be computed using localizations. In the usual literature this is generally done by using local frames and matrix calculus. We won't take this approach because the rank of the bundle interesting us could be $+\infty$. Let $\mathcal{O}$ be an open subset of $M$ and $\Psi$ a local holomorphic isomorphism between $E|_{\mathcal{O}}$ and $\mathcal{O}\times K$, where $K$ is a Hilbert space (the fiber of this trivial bundle). Then for any $x\in K$, $\Psi^{-1}(x)$ is a holomorphic section of $E|_{\mathcal{O}}$ and
\[\nabla \Psi^{-1}(x)=\Psi^{-1}(\Theta x),\]
where $\Theta$ is a (1,0)-form valued in $\mathbb{B}(K)$, called the connection operator. Let the Hermitian metric be given in this trivialization by
\[h_m(\Psi^{-1}(x),\Psi^{-1}(y))=(\mathrm{H}(m)x,y)_K,\]
where $\mathrm{H}(m)$ is a family of positive-definite bounded self-adjoint operators parameterized by $m\in \mathcal{O}$. Note that from the above formula we have $\mathrm{H}=(\Psi^{-1})^*\Psi^{-1}$. Now due to the compatibility of $\nabla$ with the Hermitian metric, we have
\begin{eqnarray*}\partial (h(\Psi^{-1}(x),\Psi^{-1}(y)))&=&(\partial \mathrm{H}x,y)_K=(\Psi^{-1}(\Theta x),\Psi^{-1}(y))\\
&=&((\Psi^{-1})^*\Psi^{-1}(\Theta x),y)_K=(\mathrm{H}\Theta x,y)_K,\end{eqnarray*}
where $\partial$ is the $(1,0)$-part of the exterior differential operator $d$. Since $x,y\in K$ are arbitrary, we have the formula $\Theta=\mathrm{H}^{-1}\partial\mathrm{H}$,
and the curvature operator is $R=\bar{\partial}(\mathrm{H}^{-1}\partial\mathrm{H})$. If $E$ is of rank 1, i.e., a line bundle, then $\frac{i}{2\pi}R=\frac{i}{2\pi}\bar{\partial}\partial\ln \mathrm{H}$ is a globally well-defined real and closed 2-form on $M$, called the first Chern form of $(E, h)$ and usually denoted by $c_1(E,h)$. A different choice of $h$ won't change the cohomology class of the first Chern form.

  In our holomorphic Hermitian vector bundle $E(T)$, we then have the canonical Chern connection $\nabla$. Since $\nabla$ is determined completely by $T$, we also call the curvature of $\nabla$ the curvature of $T$, and denote it by $R_T$. The basic goal of this subsection is to compute $R_T$ in terms of the contractive Weyl function $B(\lambda)$ and explore its meaning. We shall mainly focus on the component of $E(T)$ over $\mathbb{C_+}$, since the discussion for the component over $\mathbb{C_-}$ follows similarly.

We shall view $\mathbb{C_+}$ as a model of hyperbolic plane equipped with the standard Poincar$\acute{e}$ metric $ds^2=\frac{(du)^2+(dv)^2}{v^2}$ where $\lambda=u+iv$ for $u,v\in \mathbb{R}$. The corresponding area element is $\sigma=du\wedge dv/v^2$. The advantage of this metric is that it's invariant under the action of the M$\ddot{o}$bius transformation group $\mathbb{PSL}(2, \mathbb{R})$. This invariance is useful for other purposes in this paper. Now $R_T$ should be of the form $R_T=\frac{i}{2}\sigma \times r_T$ where $r_T$ is a smooth section of $\textup{Hom}(E(T))$. In the following, $R_T$ and $r_T$ are actually their operator form in terms of a given trivialization induced from a certain symplectic isomorphism $\Phi$.
\begin{theorem}\label{thm5}Fix a symplectic isomorphism $\Phi$. If $B(\lambda)$ is the corresponding contractive Weyl function of $T$, then the curvature operator is
\[R_T=\bar{\partial}[\Im{\lambda}(Id-B^*B)^{-1}\partial(\frac{Id-B^*B}{\Im{\lambda}})].\]
Denote $Id-B^*B$ by $\mathrm{K}$ and $Id-BB^*$ by $\tilde{\mathrm{K}}$. Then
\[r_T=Id-4(\Im{\lambda})^2\mathrm{K}^{-1}B'^*\tilde{\mathrm{K}}^{-1}B'.\]
\end{theorem}
\begin{proof}
If $\varphi, \psi\in H_+$, then $\gamma_+(\cdot)\varphi$, $\gamma_+(\cdot)\psi$ are global holomorphic sections of $E(T)$. At $\lambda\in \mathbb{C}_+$, we have
\[(T^*\gamma_+(\lambda)\varphi,\gamma_+(\lambda)\psi)_H-(\gamma_+(\lambda)\varphi,T^*\gamma_+(\lambda)\psi)_H=2i\Im{\lambda}(\gamma_+(\lambda)\varphi,\gamma_+(\lambda)\psi)_H.\]
On the other side, by the abstract Green's second formula
\begin{eqnarray*}&\quad&(T^*\gamma_+(\lambda)\varphi,\gamma_+(\lambda)\psi)_H-(\gamma_+(\lambda)\varphi,T^*\gamma_+(\lambda)\psi)_H\\
&=&i(\Gamma_+\gamma_+(\lambda)\varphi,\Gamma_+\gamma_+(\lambda)\psi)_{H_+}-i(\Gamma_-\gamma_+(\lambda)\varphi,\Gamma_-\gamma_+(\lambda)\psi)_{H_-}\\
&=&i(\varphi,\psi)_{H_+}-i(B(\lambda)\varphi, B(\lambda)\psi)_{H_-}\\
&=& i((Id-B^*B)\varphi, \psi)_{H_+}.\end{eqnarray*}
Combining all these, we obtain
\[(\gamma_+(\lambda)\varphi,\gamma_+(\lambda)\psi)_H=\frac{1}{2\Im{\lambda}}((Id-B^*B)\varphi, \psi)_{H_+}.\]
Thus in this setting
\[\mathrm{H}=\frac{Id-B^*B}{2\Im{\lambda}},\]
and the first formula follows.

Now we have
\begin{eqnarray*}R_T&=&\bar{\partial}[(\lambda-\bar{\lambda})\mathrm{K}^{-1}\partial\frac{\mathrm{K}}{\lambda-\bar{\lambda}}]=\bar{\partial}[(\lambda-\bar{\lambda})\mathrm{K}^{-1}\frac{(\lambda-\bar{\lambda})\partial \mathrm{K}-\mathrm{K}d\lambda}{(\lambda-\bar{\lambda})^2}]\\
&=&\bar{\partial}( \mathrm{K}^{-1}\partial \mathrm{K}-\frac{d\lambda}{\lambda-\bar{\lambda}}Id)\\
&=&-\mathrm{K}^{-1}\bar{\partial}\mathrm{K}\mathrm{K}^{-1}\partial \mathrm{K}+\mathrm{K}^{-1}\bar{\partial}\partial \mathrm{K}-\frac{d\bar{\lambda}d\lambda}{(\lambda-\bar{\lambda})^2}Id\\
&=&-(\mathrm{K}^{-1}\frac{\partial \mathrm{K}}{\partial \bar{\lambda}}\mathrm{K}^{-1}\frac{\partial \mathrm{K}}{\partial \lambda}(\lambda-\bar{\lambda})^2-\mathrm{K}^{-1}\frac{\partial^2 \mathrm{K}}{\partial\bar{\lambda}\partial \lambda}(\lambda-\bar{\lambda})^2+Id)\times \frac{d\bar{\lambda}d\lambda}{(\lambda-\bar{\lambda})^2}.
\end{eqnarray*}
Note that
\begin{eqnarray*}
\frac{\partial \mathrm{K}}{\partial \bar{\lambda}}\mathrm{K}^{-1}\frac{\partial \mathrm{K}}{\partial \lambda}-\frac{\partial^2 \mathrm{K}}{\partial\bar{\lambda}\partial \lambda}&=
&B'^*B\mathrm{K}^{-1}B^*B'+B'^*B'\\
&=&B'^*(B\mathrm{K}^{-1}B^*+Id)B',\end{eqnarray*}
and that $B\mathrm{K}^{-1}B^*+Id=(Id-BB^*)^{-1}$.
The second formula then follows.
\end{proof}

From the above computation, we immediately have
\begin{proposition}If two simple symmetric operators $T_1$ and $T_2$ have the same Weyl class, then they are unitarily equivalent.
\end{proposition}
\begin{proof}Since $T_1$ and $T_2$ have the same Weyl class, suitable symplectic isomorphisms $\Phi_1$ and $\Phi_2$ can be chosen such that they have the same contractive Weyl function $B(\lambda)$. In terms of the corresponding trivializations the Hermitian metrics over $E(T_1)$ and $E(T_2)$ are at the same time given by $\frac{Id-B(\lambda)^*B(\lambda)}{2\Im \lambda}$ on $\mathbb{C}_+$ and $-\frac{Id-B(\bar{\lambda})B(\bar{\lambda})^*}{2\Im \lambda}$ on $\mathbb{C}_-$. This means these two trivializations give rise to a holomorphic isometry between $E(T_1)$ and $E(T_2)$. Thus by Thm.~\ref{thm2}, $T_1$ and $T_2$ are unitarily equivalent.
 \end{proof}

As for the meaning of $r_T$, we note that $\mathrm{H}=\frac{\mathrm{K}}{2\Im{\lambda}}$ and with the new (not holomorphic) trivialization $\gamma_+(\lambda)\mathrm{H}^{-1/2}:H_+\rightarrow E_\lambda$ we have the new curvature operator
 \[\mathrm{H}^{1/2}r_T \mathrm{H}^{-1/2}=Id-4(\Im{\lambda})^2\mathrm{K}^{-1/2}B'^*\tilde{\mathrm{K}}^{-1}B'\mathrm{K}^{-1/2},\]
 which is a bounded self-adjoint operator. Then we obtain the following basic operator inequality.
 \begin{theorem}In terms of notation as above, for $\lambda\in \mathbb{C}_+$, we have
 \[Id-4(\Im{\lambda})^2\mathrm{K}^{-1/2}B'^*\tilde{\mathrm{K}}^{-1}B'\mathrm{K}^{-1/2}\geq 0,\]
 i.e., the self-adjoint operator on the left side is nonnegative.
 \end{theorem}
 \begin{proof}Though the inequality is far from obvious at first glance, the geometric proof is a simple application of a general principle concerning the curvatures of holomorphic vector bundles. It can be stated as follows \cite[pp.79]{griffiths2014principles}: "\emph{Curvature decreases in holomorphic subbundles and increases in holomorphic quotient bundles}."

 In our setting, the bundle $E(T)$ is actually a holomorphic subbundle of the trivial holomorphic Hermitian vector bundle $\mathbb{C}_+\times H$ (at $\lambda$, $E_\lambda(T)$ is $\textup{ker}(T^*-\lambda)\subset H$). The trivial connection on $\mathbb{C}_+\times H$ is flat (i.e., with zero curvature). The above principle applies and shows the curvature on $E(T)$ is non-positive. Recall that for us the curvature $R_T$ is non-positive (in the sense of Griffiths) if and only if $iR_T/\sigma$ is non-positive as a Hermitian endomorphism of $E(T)$. The statement then becomes our operator inequality.
 \end{proof}
 \emph{Remark}. The above discussion also applies to the component of $E(T)$ over $\mathbb{C}_-$, simply by replacing $B(\lambda)$ with $(B(\bar{\lambda}))^*$. We leave the details to the interested reader.

 Since actually any Nevanlinna curve can be realized as the Weyl curve of a certain simple symmetric operator, the above inequality holds for any contractive operator-valued analytic function on $\mathbb{C}_+$ (some authors call it a Schur function). Without our geometric scheme, it is hard to imagine why such an inequality should hold. It is of interest to give a direct analytic proof of this inequality without resorting to symmetric operators and the geometry of $E(T)$.

To see the meaning of the inequality, let's assume the deficiency indices of $T$ to be $(1,1)$. Then $B(\lambda)$ is a holomorphic function from $\mathbb{C}_+$  to the open unit disc in $\mathbb{C}$ and the inequality takes the following form:
\[\frac{2\Im{\lambda}|B'|}{1-|B|^2}\leq 1.\]
Set $\lambda=c(w)=i(1+w)(1-w)^{-1}$ (hence $\Im{\lambda}=\frac{1-|w|^2}{|1-w|^2}$). Then $f(w):=B(c(w))$ is a holomorphic map from the unit disc to itself. In terms of this $f$, the inequality means
\[\frac{|f'|(1-|w|^2)}{1-|f|^2}\leq 1,\]
 which is nothing else but the famous Schwarz-Pick Lemma! So our inequality is actually a generalization of this lemma, simply by replacing the target disc with a higher dimensional bounded symmetric domain. Since the Schwarz-Pick Lemma reflects the hyperbolicity of the disc, it can be imagined that the hyperbolic geometry of bounded symmetric domains ought to play a role in the investigation of symmetric operators.
 \begin{theorem}\label{thm6}If $T$ is of deficiency indices $(n,n)$, and a boundary triplet is given, then in terms of the associated Weyl function $M(\lambda)$,
 \[r_T=Id-Z^{-1}M'(\bar{\lambda})Z^{-1}M'(\lambda),\]
 where $Z=\frac{\Im M(\lambda)}{\Im \lambda}$.
 \end{theorem}
 \begin{proof}This time the metric operator is given by $\frac{\Im M(\lambda)}{\Im \lambda}$ on $\mathbb{C}_+$. The details are left to the interested reader.
 \end{proof}
 \begin{example}\label{ex5}We consider two special Nevanlinna curves with constant curvature. Let $G$ be a Hilbert space and $Id$ the identity operator over $G$. The first is $M_1(\lambda)=\lambda Id$. Obviously the corresponding symmetric operator $T$ in Thm.~\ref{thm12} has $r_T\equiv0$ and we say $T$ is flat. The second is $M_2(\lambda)$ such that it is identical to $iId$ on $\mathbb{C}_+$. For this case, $r_T\equiv Id$. Then the corresponding operator $T$ is the most "curved". To some extent, the two examples are at the opposite extreme from each other.
 \end{example}
 \begin{example}\label{ex3}Let $G$ and $Id$ be as in the above example. Define $M_3(\lambda)=(\lambda+ai)Id$ for $\lambda\in \mathbb{C}_+$ where $a>0$. Then for the corresponding symmetric operator $T_a$,
 \[r_{T_a}=(1-(\frac{\Im \lambda}{\Im \lambda +a})^2)Id.\]
 In particular, as $\Im \lambda$ approaches $+\infty$, $r_T$ goes to zero, while as $\lambda$ approaches the real axis, $r_{T_a}$ goes to $Id$.
 \end{example}
 The following example is less artificial, see also Example \ref{ex1}.
 \begin{example}Consider $B(\lambda)=e^{i\lambda}$ and let $v=\Im \lambda$. Then for the corresponding operator $T$ (with deficiency indices (1,1)),
 \[r_T=1-\frac{4v^2e^{-2v}}{(1-e^{-2v})^2}.\]
 \end{example}
 Simple symmetric operators $T$ with constant curvature $Id$ are very special in the following sense.
 \begin{proposition}\label{p12}A simple symmetric operator $T$ has the constant curvature $Id$ if and only if its single-branched Weyl map is constant.
 \end{proposition}
 \begin{proof}Choose a trivialization $\Phi$. If $W_T(\lambda)$ is constant, then its contractive Weyl function $B(\lambda)$ is constant and thus by Thm.~\ref{thm5}, $r_T\equiv Id$. Conversely, if $r_T\equiv Id$, then $\mathrm{K}^{-1/2}B'^*\tilde{\mathrm{K}}^{-1}B'\mathrm{K}^{-1/2}\equiv0$. Note that this is of the form $A^*A=0$. Hence $\tilde{\mathrm{K}}^{-1/2}B'\mathrm{K}^{-1/2}\equiv0$ and consequently $B'(\lambda)\equiv 0$. The conclusion then follows.
 \end{proof}

 When $n_+=n_-=1$, by Schwarz-Pick Lemma, if $r_T=0$ at some point $\lambda\in \mathbb{C}_+$, then $M(\lambda)=\lambda$ or is a M$\ddot{o}$bius transform of this $M(\lambda)$. We don't know if this result holds generally as the above proposition. If $r_T\equiv\gamma Id$ for some real number $\gamma\in [0,1]$, we shall say $T$ has constant curvature. It's interesting to learn more about symmetric operators of constant curvature.

 Since the curvature operator in terms of a different trivialization $\Phi'$ of $E(T)$ is similar to $r_T$, just like in Chern-Weil theory, we can define the determinant $\det(t-r_T)$ if $n_+< \infty$. This is a polynomial of degree $n_+$ in the indeterminate $t$. If
\[\det(t-r_T)=t^{n_+}-c_1(T)t^{n_+-1}+\cdots+(-1)^jc_j(T)t^{n_+-j}+\cdots+(-1)^{n_+}c_{n_+}(T),\]
$c_j(T)$ is called the $j$-th Chern function of $T$ over $\mathbb{C}_+$. These functions $c_j$ are real analytic bounded non-negative functions over $\mathbb{C}_+$ and are all unitary invariants of $T$. The construction applies as well to the component of $E(T)$ over $\mathbb{C}_-$, but essentially the same functions would be obtained because of the simple fact that $A^*A$ and $AA^*$ have the same nonzero eigenvalues.
\begin{example}\label{ex4}In Example \ref{ex3}, if $\textup{dim} G=n\in \mathbb{N}$, then \[c_1(T_a)=n(1-(\frac{\Im \lambda}{\Im \lambda +a})^2).\]
In particular, all these $T_a$'s are mutually unitarily inequivalent.
\end{example}

All these Chern functions may be of great importance, as can be conjectured from the following theorem.
\begin{theorem}\label{thm14}If $T$ is a simple symmetric operator with deficiency indices $(1,1)$, then the first Chern function $c_1(T)$ of $T$ over $\mathbb{C}_+$ determines the unitary equivalence class of $T$.
\end{theorem}
\begin{proof}The proof in essence is borrowed from \cite{cowen1978complex} which initiated the complex analysis of the now so-called Cowen-Douglas operators.\footnote{Part of our constructions in this section is indeed stimulated by the theory of Cowen-Douglas operators.} Assume the symplectic isomorphisms $\Phi_1$ for $T_1$ and $\Phi_2$ for $T_2$ are chosen and $B_1$, $B_2$ are the associated contractive Weyl functions. The corresponding trivializations will be denoted by $\gamma_+^1$ and $\gamma_+^2$ respectively. In this case, $R_T=\bar{\partial}\partial\ln(\frac{1-|B|^2}{\Im \lambda })$. If $c_1(T_1)=c_1(T_2)$, then
\[\bar{\partial}\partial\ln\frac{1-|B_1|^2}{1-|B_2|^2}=0,\]
implying that $\ln\frac{1-|B_1|^2}{1-|B_2|^2}$ is a harmonic function on $\mathbb{C}_+$. It should be of the form $f+\bar{f}$ for a holomorphic function $f$ on $\mathbb{C}_+$. Thus
\[1-|B_1|^2=e^f\times e^{\bar{f}}(1-|B_2|^2).\]
This means precisely \[(\gamma_+^1(\lambda)1, \gamma_+^1(\lambda)1)_{E_\lambda(T_1)}=(e^f\gamma_+^2(\lambda)1, e^f\gamma_+^2(\lambda)1)_{E_\lambda(T_2)},\]
where $(\cdot, \cdot)_{E_\lambda(T)}$ is the Hermitian metric on the fibre $E_\lambda(T)$. Then the bundle map $c\times \gamma_+^1(\lambda)1\mapsto c\times e^f\gamma_+^2(\lambda)1$ for any $c\in \mathbb{C}$ provides a holomorphic isometry for the two line bundles. Note that the first Chern function over $\mathbb{C}_-$ is essentilly the same as that over $\mathbb{C}_+$. The holomorphic isometry can be extended to the whole of $E(T_1)$ canonically. Thus by Thm.~\ref{thm2}, $T_1$ and $T_2$ are unitarily equivalent.
\end{proof}
\begin{example}If the deficiency indices are $(1,1)$, then
\[c_1(T)=1-\frac{4(\Im \lambda)^2|B'|^2}{(1-|B|^2)^2}=1-\frac{|M'|^2}{(\Im M/\Im \lambda)^2}.\]
It's a remarkable fact that the function $\frac{|M'|}{(\Im M/\Im \lambda)}$ doesn't depend on the specific choice of $M(\lambda)$ and completely determines its corresponding unitary equivalence class of simple symmetric operators.
\end{example}
Generally it cannot be expected that the curvature determines the unitary equivalence class of $T$ completely, but it is still interesting to know to what extent $T$ can be determined by $r_T$.
\subsection{Characteristic vector bundle of the second kind}\label{ssec2}
To prove Thm.~\ref{thm12}, for each Nevanlinna curve $N(\lambda)$ one should construct a simple symmetric operator $T$ whose Weyl curve $W_T(\lambda)$ is congruent to $N(\lambda)$ and show that  congruent Nevanlinna curves produce in this manner unitarily equivalent symmetric operators.

This has relation to the theory of functional models of symmetric operators. Indeed, if $T$ has deficiency indices $(n,n)$ and $M(\lambda)$ is the Weyl function associated to a certain boundary triplet, then one can construct a reproducing kernel Hilbert space of vector-valued holomorphic functions over $\mathbb{C}_+\cup \mathbb{C}_-$ with kernel $\frac{M(\lambda)-M(\bar{\mu})}{\lambda-\bar{\mu}}$ and multiplication by the independent variable $\lambda$ is a symmetric operator unitarily equivalent to $T$. This multiplication operator is called a model of $T$ in the literature. This construction of a reproducing kernel Hilbert space also applies to any (uniformly strict) Nevanlinna operator-valued function $N(\lambda)$ and produces a symmetric operator with $N(\lambda)$ as its Weyl function. For details of this construction, see \cite[Chap.~4]{behrndt2020boundary}.

This traditional approach is not satisfying for our purpose. In our formalism, the role of the Weyl function $M(\lambda)$ is only secondary and definitely coordinate-dependent. Consequently, the intrinsic meaning of the reproducing kernel Hilbert space in the functional model is not explicit: the kernel $\frac{M(\lambda)-M(\bar{\mu})}{\lambda-\bar{\mu}}$ has no intrinsic meaning because the relevant Hermitian symmetric space has no intrinsic linear structure and the kernel makes no sense for the coordinate-independent $W_T(\lambda)$. Conversely, we should construct a symmetric operator starting with a Nevanlinna curve rather than a Nevanlinna function.

To settle the problem we have to introduce a second characteristic vector bundle $F(T)$ to adapt the formalism of reproducing kernels to this setting. Then the intrinsic Hilbert space will be a space of holomorphic sections of $F(T)$. This also has the effect of explaining directly why there should be a functional model for $T$, without explicitly mentioning its Weyl function at all.

Let us recall the basics of reproducing kernel in a holomorphic vector bundle. This is the vector bundle version of the ordinary theory of reproducing kernel Hilbert spaces. It was invented in \cite{bertram1998reproducing} and our presentation follows \cite{koranyi2011classification} closely.

 Assume that $E\rightarrow M$ is a holomorphic Hermitian vector bundle over a complex manifold\footnote{$\textup{rk}E$ may be infinite, but the dimension of $M$ is always assumed finite.} $M$ and $\mathcal{H}$ a Hilbert space of holomorphic sections of $E$. Let $E^\dag$ be the conjugate-linear dual\footnote{We make the natural identification $(E^\dag)^\dag=E$.} of $E$. The pairing between $\varphi\in E_x^\dag$ and $\psi\in E_x$ will be denoted by $\varphi(\psi)=((\varphi, \psi))$. Obviously $E^\dag$ is anti-holomorphic, i.e., its transition functions are anti-holomorphic. We suppose the evaluation map $ev_x: \mathcal{H}\rightarrow E_x$, $s\mapsto s(x)$ is continuous for any $x\in M$. By Riesz representation theorem, we have a $\mathbb{C}$-linear map $ev_x^\dag: E_x^\dag\rightarrow \mathcal{H}$ defined by (reproducing property)
 \[(ev_x^\dag(\varphi),v)_\mathcal{H}=((\varphi, ev_x(v)))=((\varphi, v(x)))\]
 for any $\varphi\in E_x^\dag$ and $v\in \mathcal{H}$.
 We set $K(x,y)=ev_xev_y^\dag$, which is a linear map from $E^\dag_y$ to $E_x$ and called the reproducing kernel of $\mathcal{H}$. $K(x,y)$ is holomorphic in $x$, anti-holomorphic in $y$ and $K(x,y)^\dag=K(y,x)$. Also $K(x,y)$ satisfies the positivity condition: for any points $x_j\in M$, $j=1,\cdots, p$ and $\varphi_j\in E_{x_j}^\dag$, we have \[\sum_{i,j}((\varphi_i,K(x_i,x_j)\varphi_j))\geq 0,\] which is trivially
 \[(\sum_iev^\dag_{x_i}\varphi_i,\sum_jev^\dag_{x_j}\varphi_j)_\mathcal{H}\geq 0.\]
 Furthermore, the kernel $K$ is locally uniformly bounded in the sense that $K(x,x)$ is uniformly bounded over any compact subset of $M$.

 Conversely, if a kernel $K$ satisfying these properties is given, it is always the reproducing kernel of a Hilbert space $\mathcal{H}$ of holomorphic sections of $E$. $\mathcal{H}$ can be constructed as follows. Fix $y\in M$ and $\varphi\in E^\dag_y$, and let $x$ run over $M$. Then $K(x, y)\varphi$ is a holomorphic section of $E$. Define
 $$\breve{\mathcal{H}}:=\textup{span}\{K(\cdot, y)\varphi|y\in M, \varphi\in E^\dag_y\}.$$
  On $\breve{\mathcal{H}}$, we can give a (positive-definite) inner product defined by
 \[(K(\cdot, y_1)\varphi_1,K(\cdot, y_2)\varphi_2)_{\breve{\mathcal{H}}}=((K(y_2,y_1)\varphi_1,\varphi_2))\]
 for $\varphi_i\in E^\dag_{y_i}$, $i=1,2$. Completing $\breve{\mathcal{H}}$ w.r.t. this inner product then gives rise to the required $\mathcal{H}$.

 With these stated, we can now introduce our new vector bundle induced from a simple symmetric operator $T$ in a Hilbert space $H$.
 \begin{definition}Let $T$ be a simple symmetric operator in a Hilbert space $H$. We construct a vector bundle $F^\dag(T)$ over $\mathbb{C}_+\cup \mathbb{C}_-$ as follows: at each $\lambda\in \mathbb{C}_+\cup\mathbb{C}_-$, the fiber $F^\dag_\lambda(T)$ is $\textup{ker}(T^*-\bar{\lambda})$ (rather than $\textup{ker}(T^*-\lambda)$). $F^\dag(T)$ is anti-holomorphic w.r.t. the standard complex structure on $\mathbb{C}_+\cup\mathbb{C}_-$. The conjugate-linear dual $F(T)$ of $F^\dag(T)$ is holomorphic and called the characteristic vector bundle for $T$ of the second kind.
  \end{definition}
  The vector bundle $F^\dag(T)$ is obtained from $E(T)$ in the previous subsection by simply exchanging the fibers at $\lambda$ and $\bar{\lambda}$. In this sense they are essentially the same. Though the definition of $F(T)$ seems strange at first glance, it will be clear soon that this is the proper definition to produce results as expected.

  Now we have a natural continuous embedding $\iota_\lambda:F^\dag_\lambda(T)\rightarrow H$ and consequently a continuous dual linear map $\iota_\lambda^\dag: H\rightarrow F_\lambda(T)$. Therefore, we have a kernel \[K(\lambda, \mu):=\iota_\lambda^\dag\iota_\mu: F^\dag_\mu(T)\rightarrow F_\lambda(T).\] In particular, the Hermitian structure on $F_\lambda^\dag(T)$ is given by $K(\lambda, \lambda)$. It is routine to check that all the above kernel properties are satisfied. So there is a Hilbert space $\mathfrak{H}$ of holomorphic sections of $F(T)$ with $K$ as its reproducing kernel. Simplicity of $T$ implies that the linear span of all these $F^\dag_\lambda(T)$ is dense in $H$.

Elements in $\mathfrak{H}$ can be described as follows: given $x\in H$, an element $\hat{x}(\lambda)$ in $F_\lambda(T)$ arises via
  $$((\hat{x}(\lambda), w)):=(x, w)_H=(x, \iota_\lambda w)_H=((\iota_\lambda^\dag x, w))$$
   for any $w\in F^\dag_\lambda(T)=\textup{ker}(T^*-\bar{\lambda})$. This means $\hat{x}(\cdot)=\iota^\dag_{\cdot}x$ is actually a holomorphic section of $F(T)$. Due to the simplicity of $T$, the linear map $x\mapsto \hat{x}$ is injective. We can transport the Hilbert space structure of $H$ onto the image $\hat{H}$, which is precisely our $\mathfrak{H}$. This can be seen by checking that the reproducing kernel of $\hat{H}$ is just $K$. We leave the details to the interested reader.

  The multiplication operator $\mathcal{X}$ by the independent variable $\lambda$ sends a holomorphic section of $F(T)$ to another. We can define an operator $\mathfrak{T}$ in $\mathfrak{H}$ by restricting $\mathcal{X}$ to $D(\mathfrak{T})=\{f\in \mathfrak{H}|\mathcal{X} f \in \mathfrak{H}\}$.
  \begin{theorem}$\mathfrak{T}$ is a symmetric operator in $\mathfrak{H}$ and unitarily equivalent to $T$. Consequently, $\mathfrak{T}^*$ is also unitarily equivalent to $T^*$.
  \end{theorem}
  \begin{proof}By definition, the "Fourier transform" $x\mapsto \hat{x}$ is a unitary map from $H$ to $\mathfrak{H}$. If $x\in D(T)$ and $w\in F^\dag_\lambda(T)$, then
  \begin{eqnarray*}((\widehat{Tx}(\lambda), w))&=&((\iota_\lambda^\dag(Tx),w))=(Tx,\iota_\lambda w)_H=(x,T^*w)_H\\
  &=&(x, \bar{\lambda}w)_H=(\lambda x, \iota_\lambda w)_H=((\lambda\iota_\lambda^\dag x, w))\\
  &=&((\lambda \hat{x}(\lambda),w))=(((\mathcal{X}\hat{x})(\lambda),w)).
  \end{eqnarray*}
  This shows that $\widehat{Tx}=\mathcal{X}\hat{x}$. Since $Tx\in H$, we ought to have $\mathcal{X}\hat{x}\in \mathfrak{H}$. Thus $\hat{x}\in D(\mathfrak{T})$ and $\widehat{Tx}=\mathfrak{T}\hat{x}$.

  Conversely, if $\hat{x}\in D(\mathfrak{T})$, then for any $w\in F^\dag_\lambda(T)$, we have
  \begin{eqnarray*}(x,T^*w)_H&=&(x, \bar{\lambda}w)_H=(\lambda x,w)_H=(\lambda x,\iota_\lambda w)_H\\
  &=&((\lambda\iota_\lambda^\dag x, w))=((\mathcal{X}\hat{ x})(\lambda), w))=(((\mathfrak{T}\hat{x})(\lambda), w))\\
  &=&((ev_\lambda(\mathfrak{T}\hat{x}), w))=(\mathfrak{T}\hat{x}, ev_\lambda^\dag w)_{\mathfrak{H}}.
  \end{eqnarray*}
  Thus by Schwarz's inequality,
  \[|(x,T^*w)_H|\leq \|\mathfrak{T}\hat{x}\|\times \|ev_\lambda^\dag w\|.\]
  Note that
  \begin{eqnarray*}\|ev_\lambda^\dag w\|^2&=&(ev_\lambda^\dag w, ev_\lambda^\dag w)_\mathfrak{H}=((ev_\lambda ev_\lambda^\dag w, w))=((K(\lambda,\lambda)w,w))\\
  &=&((\iota_\lambda^\dag\iota_\lambda w,w))=(\iota_\lambda w,\iota_\lambda w)_H=(w,w)_H=\|w\|^2.\end{eqnarray*}
  Consequently, $|(x,T^*w)_H|\leq \|\mathfrak{T}\hat{x}\|\times \|w\|$. If $w_i\in F^\dag_{\lambda_i}(T)$, $i=1,2$, we can prove along the same line that
  \[|(x,T^*(w_1+w_2))_H|\leq \|\mathfrak{T}\hat{x}\|\times \|w_1+w_2\|.\]
  Since the linear span of $\cup_\lambda\textup{ker}(T^*-\lambda)$ is dense in $H$, $(T^*\cdot,x)_H$ can be extended as a continuous linear functional on $H$, i.e., $x\in D((T^*)^*)=D(T)$. Thus our conclusion follows.
  \end{proof}
  From this theorem, the sufficiency part of Thm.~\ref{thm2} follows immediately:
  \begin{proof}Proof of Thm.~\ref{thm2} (continued). If $E(T_1)$ and $E(T_2)$ are holomorphically isometric (say, the isometry is $\Psi$), so are $F(T_1)$ and $F(T_2)$, and the isometry is $(\Psi(\bar{\cdot}))^\dag$ defined by $$(((\Psi(\bar{\lambda}))^\dag x, \varphi))=((x, \Psi(\bar{\lambda})\varphi))$$ for $x\in F_\lambda(T_2)$ and $\varphi\in F_\lambda^\dag(T_1)$. Let $K_1$ and $K_2$ be the reproducing kernels in $F(T_1)$ and $F(T_2)$ respectively. Then that $\Psi$ is an isometry means for any $\lambda\in \mathbb{C}_+\cup \mathbb{C}_-$, $$K_1(\lambda,\lambda)=(\Psi(\bar{\lambda}))^\dag K_2(\lambda,\lambda)\Psi(\bar{\lambda}).$$
  The Identity Theorem implies
  \[K_1(\lambda,\mu)=(\Psi(\bar{\lambda}))^\dag K_2(\lambda,\mu)\Psi(\bar{\mu})\]
 for $\lambda, \mu$ in $\mathbb{C}_+$ or $\mathbb{C}_-$. Consequently $\Psi$ induces a unitary map $U$ from $\mathfrak{H}_1$ to $\mathfrak{H}_2$: $(Uf)(\lambda)=[\Psi(\bar{\lambda})^\dag]^{-1}f(\lambda)$ for any $f\in \mathfrak{H}_1$. Let $\mathfrak{T}_i$, $i=1,2$ be the corresponding model operators. Then we have $U\mathfrak{T}_1=\mathfrak{T}_2U$ simply because $\Psi$ is linear fiberwise. This shows $\mathfrak{T}_1$ and $\mathfrak{T}_2$ are unitarily equivalent and by the above theorem, so are $T_1$ and $T_2$.
   \end{proof}

The above theorem has established the basic fact that an \emph{intrinsic} functional model for $T$ can be constructed without resorting to its Weyl function or any other artificial choices. However, if a sympelectic isomorphism $\Phi$ between $\mathcal{B}_T$ and a certain standard strong symplectic Hilbert space $\mathbb{H}$ is chosen, then the bundles $F(T)$ and $F^\dag(T)$ are trivialized and the Hilbert space $\mathfrak{H}$ turns out to be a space of (vector-valued) holomorphic functions over $\mathbb{C}_+\cup \mathbb{C}_-$. When $n_+=n_-$, this, as one may hope, shall result in the traditional functional model in \cite{behrndt2020boundary}. It is indeed the case. We shall check this in the following. The relevant computations also motivate our later development.

If a trivialization $\Phi$ is chosen, we now turn to computing the local form of the reproducing kernel $K(\lambda,\mu)$.

(1) Let's assume $\lambda, \mu\in \mathbb{C}_+$ first. Then $\gamma_-(\bar{\lambda}):H_-\rightarrow F^\dag_\lambda $ produces an anti-holomorphic trivialization for $F^\dag(T)$ over $\mathbb{C}_+$. By duality we have $(\gamma_-(\bar{\lambda}))^\dag: F_\lambda\rightarrow H_-$. Let $\psi_-(\lambda):=[(\gamma_-(\bar{\lambda}))^\dag]^{-1}$. Then $\psi_-(\lambda)$ induces a holomorphic trivialization of $F(T)$ over $\mathbb{C}_+$. If $\phi,\varphi\in H_-$, then
\begin{eqnarray*}(\psi_-(\lambda)^{-1}K(\lambda,\mu)\gamma_-(\bar{\mu})\varphi, \phi)_{H_-}&=&((\iota_\lambda^\dag \iota_\mu\gamma_-(\bar{\mu})\varphi, \gamma_-(\bar{\lambda})\phi))=(\iota_\mu\gamma_-(\bar{\mu})\varphi,\iota_\lambda\gamma_-(\bar{\lambda})\phi)_{H}\\
&=&(\gamma_-(\bar{\mu})\varphi,\gamma_-(\bar{\lambda})\phi)_H.\end{eqnarray*}
Note that
\[(T^*\gamma_-(\bar{\mu})\varphi,\gamma_-(\bar{\lambda})\phi)_H-(\gamma_-(\bar{\mu})\varphi,T^*\gamma_-(\bar{\lambda})\phi)_H=(\bar{\mu}-\lambda)(\gamma_-(\bar{\mu})\varphi,\gamma_-(\bar{\lambda})\phi)_H.\]
Due to the abstract Green's formula, the left hand side of the above formula is
\begin{eqnarray*}&\quad&i(\Gamma_+\gamma_-(\bar{\mu})\varphi,\Gamma_+\gamma_-(\bar{\lambda})\phi)_{H_+}-i(\Gamma_-\gamma_-(\bar{\mu})\varphi,\Gamma_-\gamma_-(\bar{\lambda})\phi)_{H_-}\\
&=&i(B(\mu)^*\varphi,B(\lambda)^*\phi)_{H_+}-i(\varphi, \phi)_{H_-}\\
&=&i((B(\lambda)B(\mu)^*-Id)\varphi, \phi)_{H_-}.\end{eqnarray*}
We have
\[(\gamma_-(\bar{\mu})\varphi,\gamma_-(\bar{\lambda})\phi)_H=(\frac{i(Id-B(\lambda)B(\mu)^*)}{\lambda-\bar{\mu}}\varphi,\phi)_{H_-},\]
and finally obtain
\[\psi_-(\lambda)^{-1}K(\lambda,\mu)\gamma_-(\bar{\mu})=\frac{i(Id-B(\lambda)B(\mu)^*)}{\lambda-\bar{\mu}}.\]

(2) Similarly, the inverse of $\gamma_+(\bar{\lambda})$ for $\lambda\in \mathbb{C}_-$ induces a localization $\psi_+(\lambda)$ of $F(T)$ over $\mathbb{C}_-$ and for $\lambda, \mu\in \mathbb{C}_-$
\[\psi_+(\bar{\lambda})^{-1}K(\lambda, \mu)\gamma_+(\bar{\mu})=-i\frac{Id-B(\bar{\lambda})^*B(\bar{\mu})}{\lambda-\bar{\mu}}.\]

(3) If $\mu\in \mathbb{C}_+$, and $\lambda\in \mathbb{C}_-$ ($\bar{\lambda}\neq \mu$), then for $\varphi\in H_-$ and $\phi\in H_+$,
\[(\psi_+(\lambda)^{-1}K(\lambda,\mu)\gamma_-(\bar{\mu})\varphi, \phi)_{H_+}=(\gamma_-(\bar{\mu})\varphi,\gamma_+(\bar{\lambda})\phi)_H.\]
Again by the abstract Green's formula, we have
\[(\gamma_-(\bar{\mu})\varphi,\gamma_+(\bar{\lambda})\phi)_H=i(\frac{B(\bar{\lambda})^*-B(\mu)^*}{\lambda-\bar{\mu}}\varphi, \phi)_{H_+},\]
implying
\[\psi_+(\lambda)^{-1}K(\lambda,\mu)\gamma_-(\bar{\mu})=i\frac{B(\bar{\lambda})^*-B(\mu)^*}{\lambda-\bar{\mu}}.\]
For $\bar{\lambda}= \mu$, by continuity, we simply have
\[\psi_+(\lambda)^{-1}K(\lambda,\bar{\lambda})\gamma_-(\bar{\lambda})=i[B'(\bar{\lambda})]^*.\]

(4) Similarly, if $\mu\in \mathbb{C}_-$, and $\lambda\in \mathbb{C}_+$ ($\bar{\lambda}\neq \mu$), then
\[\psi_-(\lambda)^{-1}K(\lambda,\mu)\gamma_+(\bar{\mu})=-i\frac{B(\lambda)-B(\bar{\mu})}{\lambda-\bar{\mu}},\]
and for $\bar{\lambda}=\mu$, by continuity we have
\[\psi_-(\lambda)^{-1}K(\lambda,\bar{\lambda})\gamma_+(\bar{\lambda})=-iB'(\lambda).\]

 To see our approach does coincide with the traditional one included in \cite[Chap.~4]{behrndt2020boundary} for the case $n_+=n_-$, we should carry out the above computations in terms of the Weyl function $M(\lambda)$ rather than $B(\lambda)$.

 This time, we can use the $\gamma$-field $\gamma(\bar{\lambda})$ associated to a boundary triplet $(G, \Gamma_0, \Gamma_1)$ to trivialize $F^\dag(T)$. Denote $[(\gamma(\bar{\lambda}))^\dag]^{-1}$ by $\psi(\lambda)$. Then $\psi$ produces a holomorphic trivialization of $F(T)$. Note that
 \[(\psi(\lambda)^{-1}K(\lambda,\mu)\gamma(\bar{\mu})\varphi, \phi)_G=((\iota_\lambda^\dag \iota_\mu\gamma(\bar{\mu})\varphi, \gamma(\bar{\lambda})\phi))=(\gamma(\bar{\mu})\varphi, \gamma(\bar{\lambda})\phi)_H
 \]
 and that
 \[(T^*\gamma(\bar{\mu})\varphi, \gamma(\bar{\lambda})\phi)_H-(\gamma(\bar{\mu})\varphi, T^*\gamma(\bar{\lambda})\phi)_H=(\bar{\mu}-\lambda)(\gamma(\bar{\mu})\varphi, \gamma(\bar{\lambda})\phi)_H.\]
 By the abstract Green's formula, the left side of the above equation is
\begin{eqnarray*}(\Gamma_1\gamma(\bar{\mu})\varphi, \Gamma_0\gamma(\bar{\lambda})\phi)_G-(\Gamma_0\gamma(\bar{\mu})\varphi,\Gamma_1\gamma(\bar{\lambda})\phi )_G&=&(M(\bar{\mu})\varphi,\phi)_G-(\varphi, M(\bar{\lambda})\phi)_G\\
&=&((M(\bar{\mu})-M(\lambda))\varphi, \phi)_G.\end{eqnarray*}

Combining all these together, we have
\[\psi(\lambda)^{-1}K(\lambda,\mu)\gamma(\bar{\mu})=\frac{M(\lambda)-M(\bar{\mu})}{\lambda-\bar{\mu}},\]
which is precisely of the form in \cite[Chap.~4]{behrndt2020boundary}.

To motivate our later construction, we would like to determine an explicit formula for $\mathfrak{T}^*$ in terms of the contractive Weyl function $B(\lambda)$. Though $\mathfrak{T}^*$ may not be an operator, we won't worry about this issue for the necessary modification is almost clear when this does happen.
\begin{proposition}If $\hat{x}\in D(\mathfrak{T}^*)$, then for $\lambda\in \mathbb{C}_+$,
\[\psi_-(\lambda)^{-1}(\mathfrak{T}^*\hat{x})(\lambda)=\lambda\psi_-(\lambda)^{-1}\hat{x}(\lambda)+i(B(\lambda)\Gamma_+-\Gamma_-)x,\]
and for $\lambda\in \mathbb{C}_-$,
\[\psi_+(\lambda)^{-1}(\mathfrak{T}^*\hat{x})(\lambda)=\lambda\psi_+(\lambda)^{-1}\hat{x}(\lambda)+i(\Gamma_+-B(\bar{\lambda})^*\Gamma_-)x.\]
\end{proposition}
\begin{proof}If $\lambda\in \mathbb{C}_+$, then for any $w\in H_-$,
\begin{eqnarray*}&\quad& (\psi_-(\lambda)^{-1}(\mathfrak{T}^*\hat{x})(\lambda), w)_{H_-}=(((\mathfrak{T}^*\hat{x})(\lambda),\gamma_-(\bar{\lambda})w))\\&=&((\iota_\lambda^\dag T^*x, \gamma_-(\bar{\lambda})w))
=(T^*x, \iota_\lambda\gamma_-(\bar{\lambda})w)_H=(T^*x, \gamma_-(\bar{\lambda})w)_H\\
&=&(x,T^*\gamma_-(\bar{\lambda})w)_H+i(\Gamma_+x,\Gamma_+\gamma_-(\bar{\lambda})w)_{H_+}-i(\Gamma_-x,\Gamma_-\gamma_-(\bar{\lambda})w)_{H_-}\\
&=&(x,\bar{\lambda}\gamma_-(\bar{\lambda})w)_H+i(\Gamma_+x,B(\lambda)^*w)_{H_+}-i(\Gamma_-x,w)_{H_-}\\
&=&(\lambda x, \gamma_-(\bar{\lambda})w)_H+i((B(\lambda)\Gamma_+-\Gamma_-)x,w)_{H_-}.
\end{eqnarray*}
Note that
\begin{eqnarray*}(\lambda x, \gamma_-(\bar{\lambda})w)_H&=&(\lambda x, \iota_\lambda\gamma_-(\bar{\lambda})w)_H=((\lambda\iota_\lambda^\dag x, \gamma_-(\bar{\lambda})w))\\&=&((\lambda\hat{x}(\lambda),\gamma_-(\bar{\lambda})w))=(\lambda\psi_-(\lambda)^{-1}\hat{x}(\lambda),w)_{H_-}.\end{eqnarray*}
The first formula then follows. The second formula can be similarly obtained.
\end{proof}
\emph{Remark}. If $T$ has deficiency indices $(n,n)$, we can also use the Weyl function $M(\lambda)$ and the $\gamma$-field $\gamma$ to carry out a similar computation. Then the formula takes the following form
\[\psi(\lambda)^{-1}(\mathfrak{T}^*\hat{x})(\lambda)=\lambda\psi(\lambda)^{-1}\hat{x}(\lambda)+(\Gamma_1-M(\lambda)\Gamma_0)x\]
for any $\hat{x}\in D(\mathfrak{T}^*)$.

Let's now come to the problem of constructing a simple symmetric operator for a given Nevanlinna curve. We assume $N(\lambda)$ is a Nevanlinna curve for the standard strong symplectic Hilbert space $\mathbb{H}=H_+\oplus_\bot H_-$. We extend $N(\lambda)$ to $\mathbb{C}_+\cup \mathbb{C}_-$ by setting $N(\lambda)=(N(\bar{\lambda}))^{\bot_s}$ for $\lambda\in \mathbb{C}_-$. $W_+(\mathbb{H})\cup W_-(\mathbb{H})$ has a tautological vector bundle for which the fiber at $N_\pm\in W_\pm(\mathbb{H})$ is $N_\pm$ itself. We pull back this bundle to $\mathbb{C}_+\cup \mathbb{C}_-$ through $N(\lambda)$ and exchange the fibers at $\lambda$ and $\bar{\lambda}$. This results in an anti-holomorphic bundle $F^\dag$ over $\mathbb{C}_+\cup \mathbb{C}_-$. Let $F$ be the conjugate-linear dual of $F^\dag$. $F^\dag$ has a natural Hermitan metric induced from the strong symplectic structure on $\mathbb{H}$ by restriction. Then $F$ is equipped with the dual metric.

We need a reproducing kernel $K(\lambda, \mu)$ to produce a Hilbert space of holomorphic sections of $F$. Motivated by the investigation of the characteristic bundle $F(T)$ for $T$, this can be done as follows. Note that at each $\lambda\in \mathbb{C}_+\cup \mathbb{C}_-$ there is the decomposition $\mathbb{H}=F^\dag_\lambda\oplus F^\dag_{\bar{\lambda}}$. We denote the projection along $F^\dag_{\bar{\lambda}}$ onto $F^\dag_\lambda$ by $\mathcal{P}_\lambda$. For $\lambda, \mu\in \mathbb{C}_+\cup \mathbb{C}_-$, we can construct a map $\mathcal{Q}(\lambda, \mu)$ from $F^\dag_\mu$ to $F_\lambda$: if $\varphi\in F^\dag_\mu$,
then $\mathcal{P}_\lambda\varphi\in F^\dag_\lambda$ and $(\mathcal{P}_\lambda\varphi,\cdot)$ (the Hermitian structure on the fiber is used) is a conjugate-linear functional on $F^\dag_\lambda$. Consequently, $(\mathcal{P}_\lambda\varphi,\cdot)\in F_\lambda$ and we simply set $\mathcal{Q}(\lambda, \mu)\varphi=(\mathcal{P}_\lambda\varphi,\cdot)$. Now for $\lambda\neq \bar{\mu}$ we define
\[K(\lambda,\mu)=i\frac{\mathcal{Q}(\lambda, \mu)}{\lambda-\bar{\mu}}.\]
As for $\lambda=\bar{\mu}$, $\mathcal{Q}(\bar{\mu},\mu)$ is certainly zero. It turns out that $\mathcal{Q}(\lambda, \mu)$ as an analytic function in $\lambda$ has a zero at $\bar{\mu}$, and we can simply set $K(\bar{\mu},\mu)$ to be the limit of the above formula as $\lambda$ approaches $\bar{\mu}$.
\begin{theorem}\label{thm3}The above $K(\lambda,\mu)$ is the reproducing kernel of a Hilbert space $\mathfrak{H}$ of holomorphic sections of $F$. Set $\mathfrak{X}$ to be multiplication by the independent variable $\lambda$ on holomorphic sections of $F$ and define $\mathcal{D}=\{s\in \mathfrak{H}|\mathfrak{X}s\in \mathfrak{H}\}$. Then the restriction $\mathfrak{T}$ of $\mathfrak{X}$ on $\mathcal{D}$ is a simple symmetric operator, whose Weyl class is just the congruence class $[N(\lambda)]$.
\end{theorem}
The proof will be accomplished in several stages. It is a generalization of the proof of Thm.~4.2.4 in \cite{behrndt2020boundary}--we have translated it into geometric language and used the contractive Weyl function instead of the Weyl function, which only makes sense for the case $n_+=n_-$. Due to the conceptual transparency of our construction, we don't need the integral representation of an operator-valued Nevanlinna function, which definitely depends on viewing a Nevanlinna curve as a function.

We shall first give the local form of $K(\lambda,\mu)$ and then the required properties for $K(\lambda,\mu)$ to be a reproducing kernel can be checked immediately. The result takes the same form as before, and so we only do the computation for $\lambda,\mu\in \mathbb{C}_+$. The rest can be given in the same manner.

We shall identify the fibers $F_{\mp i}^\dag$ with $H_\pm$. Then $N(\lambda)$ for $\lambda\in \mathbb{C}_+$ is described by $B(\lambda)\in \mathbb{B}(H_+,H_-)$; in particular, $B(i)=0$. This gives a trivialization of $F^\dag$: $F_\lambda^\dag$ for $\lambda\in \mathbb{C}_+$ is identified with $$\{(B(\lambda)^*x,x)\in H_+\oplus_\bot H_-|x\in H_-\},$$ while $F_\lambda^\dag$ for $\lambda\in \mathbb{C}_-$ is identified with $$\{(x,B(\bar{\lambda})x)\in H_+\oplus_\bot H_-|x\in H_+\}.$$ In this sense, we set $\gamma_+(\bar{\lambda})x=(x,B(\bar{\lambda})x)$ for $x\in H_+$, $\lambda\in \mathbb{C}_-$, and set $\gamma_-(\bar{\lambda})x=(B(\lambda)^*x,x)$ for $x\in H_-$, $\lambda\in \mathbb{C}_+$.

For any $x\in H_-$, we have $\gamma_-(\bar{\mu})x=(B(\mu)^*x,x)$. Due to the decomposition $\mathbb{H}=F^\dag_\lambda\oplus F^\dag_{\bar{\lambda}}$, there are unique $y_1\in H_-$ and $y_2\in H_+$ such that
\[(B(\mu)^*x,x)=(B(\lambda)^*y_1,y_1)+(y_2,B(\lambda)y_2),\]
i.e.,
\[B(\mu)^*x=B(\lambda)^*y_1+y_2,\quad x=y_1+B(\lambda)y_2,\]
from which we can find
\[y_1=(Id-B(\lambda)B(\lambda)^*)^{-1}(Id-B(\lambda)B(\mu)^*)x.\]
Note that $\mathcal{P}_\lambda\gamma_-(\bar{\mu})x=\gamma_-(\bar{\lambda})y_1$. Therefore, for any $w\in H_-$, we have
\[(\mathcal{P}_\lambda\gamma_-(\bar{\mu})x,\gamma_-(\bar{\lambda})w)=((Id-B(\lambda)B(\lambda)^*)y_1,w)_{H_-}=((Id-B(\lambda)B(\mu)^*)x,w)_{H_-}.\]
By definition,
\[((\mathcal{Q}(\lambda,\mu)\gamma_-(\bar{\mu})x,\gamma_-(\bar{\lambda})w))=(\mathcal{P}_\lambda\gamma_-(\bar{\mu})x,\gamma_-(\bar{\lambda})w)=((Id-B(\lambda)B(\mu)^*)x,w)_{H_-}.\]
This shows that
\[\psi_-(\lambda)^{-1}\mathcal{Q}(\lambda,\mu)\gamma_-(\bar{\mu})x=(Id-B(\lambda)B(\mu)^*)x,\]
and consequently,
\[\psi_-(\lambda)^{-1}K(\lambda,\mu)\gamma_-(\bar{\mu})=\frac{i(Id-B(\lambda)B(\mu)^*)}{\lambda-\bar{\mu}},\]
where $\psi_-(\lambda)^{-1}=(\gamma_-(\bar{\lambda}))^\dag$.

We have to prove multiplication by the independent variable in $\mathfrak{H}$ is a simple symmetric operator $A$ and the bundles $F^\dag$ and $F^\dag(A)$ are anti-holomorphically isometric. The strategy is to encode the data into a Lagrangian subspace of a bigger strong symplectic Hilbert space. For convenience, we shall use the above trivialization to identify $\mathfrak{H}$ with a Hilbert space of holomorphic functions on $\mathbb{C}_+\cup \mathbb{C_-}$. For $f\in \mathfrak{H}$, $f(\lambda)\in H_-$ if $\lambda\in \mathbb{C}_+$, and $f(\lambda)\in H_+$ if $\lambda\in \mathbb{C}_-$. With this identification, for $\lambda, \mu\in \mathbb{C}_+$
\[K(\lambda,\mu)=\frac{i(Id-B(\lambda)B(\mu)^*)}{\lambda-\bar{\mu}};\]
for $\lambda, \mu\in \mathbb{C}_-$,
\[K(\lambda,\mu)=-i\frac{Id-B(\bar{\lambda})^*B(\bar{\mu})}{\lambda-\bar{\mu}};\]
for $\mu\in \mathbb{C}_+$, and $\lambda\in \mathbb{C}_-$,
\[K(\lambda,\mu)=i\frac{B(\bar{\lambda})^*-B(\mu)^*}{\lambda-\bar{\mu}};\]
for $\mu\in \mathbb{C}_-$, and $\lambda\in \mathbb{C}_+$,
\[K(\lambda,\mu)=-i\frac{B(\lambda)-B(\bar{\mu})}{\lambda-\bar{\mu}}.\]
If $\lambda\in \mathbb{C}_+$, $K(\lambda, \bar{\lambda})=-iB'(\lambda)$, and if $\lambda\in \mathbb{C}_-$, $K(\lambda, \bar{\lambda})=iB'(\bar{\lambda})^*$.

We denote by $\check{\mathbb{H}}$ the same space $\mathbb{H}$ but with the strong symplectic structure $-[\cdot,\cdot]$. Consider the direct sum $(\mathcal{H},[\cdot, \cdot])$ of the strong symplectic Hilbert spaces $\mathfrak{H}\oplus_\bot\mathfrak{H}$ and $\check{\mathbb{H}}$. Note that $\mathfrak{H}\oplus_\bot\mathfrak{H}$ is equipped with the second standard strong symplectic structure.
\begin{proposition}\label{p7}Consider the subspace $\mathbb{L}$ of $\mathcal{H}$ defined by
$$\{(f,g;\varphi_+,\varphi_-)\in \mathcal{H}| \textup{for} \, \lambda\in \mathbb{C}_+,\, g(\lambda)=\lambda f(\lambda)+i(B(\lambda)\varphi_+-\varphi_-),$$
$$\, \textup{for} \, \lambda\in \mathbb{C}_-,\, g(\lambda)=\lambda f(\lambda)+i(\varphi_+-B(\bar{\lambda})^*\varphi_-)\}.$$
Then $\mathbb{L}$ is Lagrangian in $\mathcal{H}$.
\end{proposition}
\begin{proof}We say a subspace of $\mathcal{H}$ is essentially Lagrangian, if its closure is Lagrangian. Since the evaluation map at each $\lambda\in \mathbb{C}_+\cup \mathbb{C}_-$ is continuous, $\mathbb{L}$ is clearly closed. Now to prove the claim, it suffices to find an essentially Lagrangian subspace densely contained in $\mathbb{L}$.

Consider the two sets in $\mathcal{H}$:
\[\mathbb{B}_1=\{(K(\cdot, \bar{\mu})\varphi, \mu K(\cdot, \bar{\mu})\varphi;\varphi, B(\mu)\varphi)\in \mathcal{H}|\mu\in \mathbb{C}_+, \varphi\in H_+\}\]
and
\[\mathbb{B}_2=\{(K(\cdot, \bar{\mu})\varphi, \mu K(\cdot, \bar{\mu})\varphi;B(\bar{\mu})^*\varphi, \varphi)\in \mathcal{H}|\mu\in \mathbb{C}_-, \varphi\in H_-\}.\]
We have $\mathbb{B}_1\subset \mathbb{L}$. Indeed, let $\mu\in \mathbb{C}_+$ and $\varphi\in H_+$, and set $s(\cdot)=K(\cdot, \bar{\mu})\varphi$. Then for $\lambda\in \mathbb{C_+}$ we have
\[(\lambda-\mu)s(\lambda)=-i(B(\lambda)\varphi-B(\mu)\varphi),\]
i.e.,
\[\mu s(\lambda)=\lambda s(\lambda)+i(B(\lambda)\varphi-B(\mu)\varphi).\]
For $\lambda\in \mathbb{C_-}$ we have
\[(\lambda-\mu)s(\lambda)=-i(\varphi-B(\bar{\lambda})^*B(\mu)\varphi),\]
i.e.,
\[\mu s(\lambda)=\lambda s(\lambda)+i(\varphi-B(\bar{\lambda})^*B(\mu)\varphi).\]
Similarly, we have $\mathbb{B}_2\subset \mathbb{L}$. We set $\mathbb{B}:=\textup{span}\{\mathbb{B}_1\cup \mathbb{B}_2\}$ and have $\mathbb{B}\subset \mathbb{L}$. It can be justified by direct computations that $\mathbb{B}$ is isotropic in $\mathcal{H}$. We only prove $[\cdot, \cdot]$ restricted on $\mathbb{B}_1$ vanishes and leave the remaining computations to the interested reader: if
\[a_i:=(K(\cdot, \bar{\mu}_i)\varphi_i, \mu_i K(\cdot, \bar{\mu}_i)\varphi_i;\varphi_i, B(\mu_i)\varphi_i),\quad i=1,2\]
are in $\mathbb{B}_1$, then
\begin{eqnarray*}[a_1,a_2]&=&(\mu_1 K(\cdot, \bar{\mu}_1)\varphi_1, K(\cdot, \bar{\mu}_2)\varphi_2)_\mathfrak{H}-(K(\cdot, \bar{\mu}_1)\varphi_1,\mu_2 K(\cdot, \bar{\mu}_2)\varphi_2)_\mathfrak{H}\\
&-&i(\varphi_1,\varphi_2)_{H_+}+i(B(\mu_1)\varphi_1,B(\mu_2)\varphi_2)_{H_-}\\
&=&(\mu_1-\bar{\mu}_2)(K(\bar{\mu}_2, \bar{\mu}_1)\varphi_1,\varphi_2)_{H_+}-i(\varphi_1,\varphi_2)_{H_+}+i(B(\mu_2)^*B(\mu_1)\varphi_1,\varphi_2)_{H_+}.
\end{eqnarray*}
Note that in the present setting,
\[(\mu_1-\bar{\mu}_2)K(\bar{\mu}_2, \bar{\mu}_1)=i(Id-B(\mu_2)^*B(\mu_1)).\]
We see immediately that $[a_1,a_2]=0$.

For $(f,\tilde{f})\in \mathfrak{H}\oplus_\bot\mathfrak{H}$, we set $\beta_\pm((f,\tilde{f}))=(\tilde{f}\pm if)/\sqrt{2}\in \mathfrak{H}$. Then obviously
\begin{eqnarray*}[(f, \tilde{f};\varphi_+,\varphi_-), (g, \tilde{g};\phi_+,\phi_-)]&=&i(\beta_+(f, \tilde{f}),\beta_+(g, \tilde{g}))_{\mathfrak{H}}+i(\varphi_-,\phi_-)_{H_-}\\&-&i(\beta_-(f, \tilde{f}),\beta_-(g, \tilde{g}))_{\mathfrak{H}}-i(\varphi_+,\phi_+)_{H_+}.\end{eqnarray*}
Clearly, $\beta_\pm$ transforms the second standard strong symplectic structure on $\mathfrak{H}\oplus_\bot\mathfrak{H}$ into the first one. Let
\[\mathbb{B}_\pm=\{(\beta_\pm((f,\tilde{f})), \varphi_\mp)\in \mathfrak{H}\oplus_\bot H_\mp|(f,\tilde{f};\varphi_+, \varphi_-)\in \mathbb{B}\}.\]
Since $\mathbb{B}$ is isotropic, it obviously defines a map $\tau$ from $\mathbb{B}_+$ to $\mathbb{B}_-$, preserving the inner product. To prove $\mathbb{B}$ is essentially Lagrangian, we have to prove $\mathbb{B}_\pm$ are dense in $\mathfrak{H}\oplus_\bot H_\mp$ respectively and hence $\tau$ can be extended to a unitary map from $\mathfrak{H}\oplus_\bot H_-$ to $\mathfrak{H}\oplus_\bot H_+$. If $(f, \psi)\in \mathfrak{H}\oplus_\bot H_-$ is orthogonal to $\mathbb{B}_+$, then for any $\mu\in \mathbb{C}_+$ and $\varphi\in H_+$, we have
\[(f, \frac{(\mu+i)(K(\cdot, \bar{\mu})\varphi)}{\sqrt{2}})_{\mathfrak{H}}+(\psi, B(\mu)\varphi)_{H_-}=0,\]
and for any $\mu\in \mathbb{C}_-$ and $\varphi\in H_-$, we have
\[(f, \frac{(\mu+i)(K(\cdot, \bar{\mu})\varphi)}{\sqrt{2}})_{\mathfrak{H}}+(\psi, \varphi)_{H_-}=0.\]
By the reproducing property, these equations imply
\[\frac{\bar{\mu}-i}{\sqrt{2}}f(\bar{\mu})+B(\mu)^*\psi=0\]
for $\mu\in \mathbb{C}_+$ and
\[\frac{\bar{\mu}-i}{\sqrt{2}}f(\bar{\mu})+\psi=0\]
for $\mu\in \mathbb{C}_-$. Setting $\mu=-i$ in the second equation, we see $\psi=0$ and $f\equiv 0$ follows immediately. The denseness of $\mathbb{B}_-$ in $\mathfrak{H}\oplus_\bot H_+$ holds similarly.

The remainder is to prove $\overline{\mathbb{B}}=\mathbb{L}$. Since $\mathbb{B}\subset \mathbb{L}$, we only have to prove $\mathbb{L}\subset \overline{\mathbb{B}}=\mathbb{B}^{\bot_s}$. For $(f,\tilde{f}; \varphi_+, \varphi_-)\in \mathbb{L}$, $\mu\in \mathbb{C}_+$ and $\varphi\in H_+$,
\begin{eqnarray*}&\quad&(\tilde{f},  K(\cdot, \bar{\mu})\varphi)_\mathfrak{H}-(f, \mu K(\cdot, \bar{\mu})\varphi)_\mathfrak{H}-i(\varphi_+, \varphi)_{H_+}+i(\varphi_-, B(\mu)\varphi)_{H_-}\\
&=&(\tilde{f}(\bar{\mu})-\bar{\mu}f(\bar{\mu}),\varphi)_{H_+}+i(-\varphi_++B(\mu)^*\varphi_-, \varphi)_{H_+}.
\end{eqnarray*}
The above line vanishes by definition of $\mathbb{L}$ and thus we have $\mathbb{L}\subset \mathbb{B}_1^{\bot_s}$. Similarly, we can prove $\mathbb{L}\subset \mathbb{B}_2^{\bot_s}$. Consequently, we have \[\mathbb{L}\subset \mathbb{B}_1^{\bot_s}\cap \mathbb{B}_2^{\bot_s}=(\textup{span}\{\mathbb{B}_1\cup \mathbb{B}_2\})^{\bot_s}=\mathbb{B}^{\bot_s}.\]
\end{proof}
The proof of Prop.~\ref{p7} has the following byproduct, which is not obvious from first glance.
\begin{corollary}If $n_+\neq n_-$, then $\textup{dim}\mathfrak{H}=+\infty$.
\end{corollary}
\begin{proof}This is simply because the map $\tau$ is actually a unitary map from $\mathfrak{H}\oplus_\bot H_-$ to $\mathfrak{H}\oplus_\bot H_+$.
\end{proof}
The corollary is not trivial since the statement is not necessarily true when $n_+=n_-$.
\begin{example}Let $H$ be finite-dimensional and $D(T)=0\subset H$. We say $T$ is the zero symmetric operator in $H$. In this case, the Weyl curve is the universal Weyl curve associated to $H$.
\end{example}

\begin{proposition}Let $\mathfrak{T}$ be the operator defined in the statement of Thm.~\ref{thm3}. In terms of the above identification again, the adjoint of $\mathfrak{T}$ (as a linear relation)
is
$$\mathfrak{T}^*=\{(f,g)\in \mathfrak{H}\oplus_\bot\mathfrak{H}|\exists \varphi_\pm\in H_\pm, \textup{for}\, \lambda\in \mathbb{C}_+, g(\lambda)=\lambda f(\lambda)+i(B(\lambda)\varphi_+-\varphi_-),$$
$$\, \textup{for} \, \lambda\in \mathbb{C}_-, \, g(\lambda)=\lambda f(\lambda)+i(\varphi_+-B(\bar{\lambda})^*\varphi_-)\}.$$
\end{proposition}
\begin{proof}Note that $(f,g)\in A_\mathfrak{T}$ (the graph of $\mathfrak{T}$) if and only if $(f,g;0, 0)\in \mathbb{L}\cap (\mathfrak{H}\oplus_\bot\mathfrak{H}+0)\subset \mathcal{H}$. This clearly implies that $\mathfrak{T}$ is a closed symmetric operator. Denote the linear relation in the statement of the proposition by $T_B$. Then obviously,
\[(f,g)\in T_B^{\bot_s}\Leftrightarrow (f,g,0,0)\in \mathbb{L}^{\bot_s}=\mathbb{L}\Leftrightarrow (f,g)\in A_\mathfrak{T}.\]
Thus $ T_B^{\bot_s}=A_\mathfrak{T}$ and consequently $\overline{T_B}=\mathfrak{T}^*$. Thus to prove the conclusion, it suffices to prove that $T_B$ is closed. If $(f_n, g_n)\in T_B$ is a Cauchy sequence in $\mathfrak{H}\oplus_\bot\mathfrak{H}$ with its limit $(f,g)$, then there is a sequence $(\varphi_{+n}, \varphi_{-n})\in H_+\oplus_\bot H_-$ such that for $\mu\in \mathbb{C}_+$,
\[g_n(\mu)=\mu f_n(\mu)+i(B(\mu)\varphi_{+n}-\varphi_{-n}),\]
and for $\mu\in \mathbb{C}_-$,
\[g_n(\mu)=\mu f_n(\mu)+i(\varphi_{+n}-B(\bar{\mu})^*\varphi_{-n}).\]
Since the evaluation map at each $\lambda\in \mathbb{C}_+\cup \mathbb{C}_-$ is continuous, we know that for $\mu=\pm i$
\[i(B(i)\varphi_{+n}-\varphi_{-n})\rightarrow g(i)-i f(i), \quad \textup{as}\quad n\rightarrow \infty,\]
and
\[i(\varphi_{+n}-B(i)^*\varphi_{-n})\rightarrow g(-i)+i f(-i), \quad \textup{as}\quad n\rightarrow \infty.\]
Since $B(i)=0$, from these we find that as $n\rightarrow \infty$, $\varphi_{+n}\rightarrow \varphi_+$ in $H_+$ and $\varphi_{-n}\rightarrow \varphi_-$ for some $\varphi_\pm\in H_\pm$. Therefore, we have for $\mu\in \mathbb{C}_+$,
\[g(\mu)=\mu f(\mu)+i(B(\mu)\varphi_+-\varphi_-),\]
and for $\mu\in \mathbb{C}_-$,
\[g(\mu)=\mu f(\mu)+i(\varphi_+-B(\bar{\mu})^*\varphi_-).\]
This shows $(f,g)\in T_B$ and $T_B$ is thus closed.
\end{proof}
\begin{proposition}The symmetric operator $\mathfrak{T}$ defined in the statement of Thm.~\ref{thm3} is simple and its Weyl class is $[N(\lambda)]$.
\end{proposition}
\begin{proof}We continue using the chosen trivialization to facilitate our calculation. Let's compute $\textup{ker}(\mathfrak{T}^*-\mu)$ for $\mu\in \mathbb{C}_+\cup \mathbb{C}_-$. If $(f,g)\in \textup{ker}(\mathfrak{T}^*-\mu)$, then $g(\lambda)-\mu f(\lambda)=0$ for any $\lambda\in \mathbb{C}_+\cup \mathbb{C}_-$. If $\mu\in \mathbb{C}_+$, then
\[g(\lambda)-\mu f(\lambda)=\lambda f(\lambda)+i(B(\lambda)\varphi_+-\varphi_-)-\mu f(\lambda)=0, \quad \forall \lambda\in \mathbb{C}_+.\]
In particular setting $\lambda=\mu$, we get
\begin{equation}\varphi_-=B(\mu)\varphi_+.\label{eq3}\end{equation} Consequently,
\[f(\lambda)=-i\frac{(B(\lambda)-B(\mu))\varphi_+}{\lambda-\mu}.\]
With these in hand, we see
\[0=g(\lambda)-\mu f(\lambda)=\lambda f(\lambda)+i(\varphi_+-B(\bar{\lambda})^*B(\mu)\varphi_+)-\mu f(\lambda),\quad \forall \lambda\in \mathbb{C}_-.\]
Therefore,
\[f(\lambda)=-i\frac{(Id-B(\bar{\lambda})^*B(\mu))\varphi_+}{\lambda-\mu}.\]
Comparing these expressions with the local form of the reproducing kernel $K$, we find $f(\cdot)=K(\cdot, \bar{\mu})\varphi_+$. Conversely, it's also very easy to see $f\in \mathfrak{H}$ of the form $K(\cdot, \bar{\mu})\varphi_+$ lies in $\textup{ker}(\mathfrak{T}^*-\mu)$. Therefore,
\[\textup{ker}(\mathfrak{T}^*-\mu)=\{K(\cdot, \bar{\mu})\varphi|\varphi\in H_+\}.\]
Similarly, this identity also holds for $\mu\in \mathbb{C}_-$ except that $H_+$ should be replaced with $H_-$. Since $\textup{span}\{K(\cdot, \bar{\mu})\varphi|\mu\in \mathbb{C}_\pm, \varphi\in H_\pm\}$ is dense in $\mathfrak{H}$ by definition, this shows that $\mathfrak{T}$ is simple.

To prove the Weyl class of $\mathfrak{T}$ is $[N(\lambda)]$, we have to construct an isomorphism $\Phi$ between the strong symplectic Hilbert spaces $D(\mathfrak{T}^*)/D(\mathfrak{T})$ and $\mathbb{H}=H_+\oplus_\bot H_-$. This is very easy: by definition, if $(f,g)\in \mathfrak{T}^*$, there are $\varphi_\pm\in H_\pm$ such that
\[g(\lambda)=\lambda f(\lambda)+i(B(\lambda)\varphi_+-\varphi_-),\quad \forall \lambda\in \mathbb{C}_+\]
and
\[g(\lambda)=\lambda f(\lambda)+i(\varphi_+-B(\bar{\lambda})^*\varphi_-),\quad \forall \lambda\in \mathbb{C}_-.\]
We simply set $\Gamma_\pm(f,g)=\varphi_\pm$. These maps are well-defined and $\textup{ker}\Gamma_+\cap \textup{ker}\Gamma_-=A_\mathfrak{T}$. That $\mathbb{L}$ is Lagrangian implies that Green's second formula holds, i.e., $\Phi$ preserves strong symplectic structures. Additionally, $\Phi$ is onto because for any fixed $\mu\in \mathbb{C}_+$, $\Phi(\textup{ker}(\mathfrak{T}^*-\mu))$ and $\Phi(\textup{ker}(\mathfrak{T}^*-\bar{\mu}))$ are transversal in $\mathbb{H}$, while $\textup{ker}(\mathfrak{T}^*-\mu)+\textup{ker}(\mathfrak{T}^*-\bar{\mu})\subset D(\mathfrak{T}^*)$. With these, the formula (\ref{eq3}) already shows $B(\lambda)$ is the associated contractive Weyl function.
\end{proof}
If $[N_1(\lambda)]=[N_2(\lambda)]$, we can choose trivializations such that $N_1(\lambda)$ and $N_2(\lambda)$ share the same contractive operator-valued function $B(\lambda)$ and then the corresponding  multiplication operators $\mathfrak{T}_1$ and $\mathfrak{T}_2$ are obviously unitarily equivalent.

Let $T$ (resp. $T'$) be a simple symmetric operator with deficiency indices $(n,n)$ in $H$ (resp. $H'$), and $(G, \Gamma_0, \Gamma_1)$ (resp. $(G, \Gamma'_0, \Gamma'_1)$) a boundary triplet associated to it. If the corresponding Weyl functions coincide, then of course $T$ and $T'$ are unitarily equivalent, but even more information can be derived.
\begin{theorem}If $T$ and $T'$ are simple symmetric operators with deficiency indices $(n,n)$ and they have the same Weyl function $M(\lambda)$ in the above sense, then there is a unitary map $U$ from $H$ to $H'$ such that
\[UT=T'U,\quad UT_0=T_0'U,\quad UT_1=T_1'U.\]
\end{theorem}
\emph{Remark}. The formula $UT=T'U$ should be interpreted as follows: $UD(T)=D(T')$ and $UTx=T'Ux$ for any $x\in D(T)$. The other two formulae should be interpreted in a similar manner.
\begin{proof}Let $\gamma(\lambda)$ (resp. $\gamma'(\lambda)$) be the associated $\gamma$-field. Since they share the same Weyl function $M(\lambda)$, for any $\varphi\in G$, $\|\gamma(\lambda)\varphi\|_{H_1}=\|\gamma'(\lambda)\varphi\|_{H_2}$. This means the map \[\gamma(\lambda)\varphi\in \textup{ker}(T^*-\lambda)\rightarrow \gamma'(\lambda)\varphi\in \textup{ker}(T'^*-\lambda)\]
is unitary from $\textup{ker}(T^*-\lambda)$ to $\textup{ker}(T'^*-\lambda)$. Due to simplicity of $T$ and $T'$, this map can be extended to a unitary map $U$ from $H_1$ to $H_2$ \cite[Thm.~4.2.6]{behrndt2020boundary}. By definition, we have $U\gamma(\lambda)=\gamma'(\lambda)$. Just as in Prop.~\ref{p2}, we can easily obtain for $\lambda, \mu\in \mathbb{C}_+\cup \mathbb{C}_-$
\[\gamma(\lambda)=(Id+(\lambda-\mu)(T_0-\lambda)^{-1})\gamma(\mu).\]
A similar formula for $\gamma'(\lambda)$ holds. Thus for arbitrary $\lambda, \mu\in \mathbb{C}_+\cup \mathbb{C}_-$,
\[U(Id+(\lambda-\mu)(T_0-\lambda)^{-1})\gamma(\mu)=(Id+(\lambda-\mu)(T'_0-\lambda)^{-1})U\gamma(\mu).\]
Again due to simplicity of $T$ and $T'$, this implies for any $\lambda\in \mathbb{C}_+\cup \mathbb{C}_-$
\[U(T_0-\lambda)^{-1}=(T'_0-\lambda)^{-1}U,\]
and thus $UT_0=T_0'U$ \cite[Lemma.~1.3.8]{behrndt2020boundary}.

We can define $\tilde{\gamma}(\lambda)\varphi=x$ to be the unique solution of the abstract boundary value problem
\[\left\{
\begin{array}{ll}
T^*x=\lambda x,\\
\Gamma_1x=\varphi.
\end{array}
\right.\]
For $T'$, a similar $\tilde{\gamma'}(\lambda)$ can be defined. By definition, for any $\lambda\in \mathbb{C}_+\cup \mathbb{C}_-$ and $\varphi\in G$,
\[\tilde{\gamma}(\lambda)M(\lambda)\varphi=\gamma(\lambda)\varphi,\quad \tilde{\gamma'}(\lambda)M(\lambda)\varphi=\gamma'(\lambda)\varphi.\]
Then $U\tilde{\gamma}(\lambda)M(\lambda)=\tilde{\gamma'}(\lambda)M(\lambda)$,
and consequently $U\tilde{\gamma}(\lambda)=\tilde{\gamma'}(\lambda)$. Proceeding as before, we can prove $UT_1=T_1'U$.
\end{proof}
The theorem means the Weyl function determines the unitary equivalence class $[(T, T_0, T_1)]$ of the triple $(T, T_0, T_1)$. In some applications, it's possible to find a unique representative in
$[(T, T_0, T_1)]$, which takes a specific form. See \S~\ref{sturm} for an example.

\subsection{Characteristic vector bundle of the third kind }\label{sec6}
The investigation of the characteristic bundle $E(T)$ poses the problem to determine when a holomorphic Hermitian vector bundle over $\mathbb{C}_+\cup \mathbb{C}_-$ is of the form $E(T)$ for a certain simple symmetric operator $T$. It is this problem that leads us to the relation of simple symmetric operators with a special kind of Higgs bundles in the non-abelian Hodge theory.

Recall that a Higgs bundle over a Riemann surface $M$ is a holomorphic vector bundle $E$ over $M$ together with a holomorphic section (called Higgs field) $\Psi$ of $\textup{Hom}(E)\otimes \kappa$ where $\kappa$ is the holomorphic cotangent bundle (i.e., canonical line bundle) of $M$. Higgs bundles arose originally in N. Hitchin's work \cite{hitchin1987self} in 1980s on dimensional reduction of Yang-Mills equations in real dimension 4 and since then have played increasingly important roles in several mathematical disciplines. In particular, Higgs bundles have been used to study character varieties of representations of the fundamental group $\pi_1(M)$. Basically, in the literature it is often assumed that $M$ is compact and $E$ is of finite rank.

A Higgs bundle $(E, \Psi)$ over $M$ is called a $\mathbb{U}(p,q)$-Higgs bundle, if $E=E_1\oplus E_2$ where $E_1, E_2\subset E$ are holomorphic subbundles of rank $p$ and $q$ respectively, and $\Psi=\left(
                                                                                                                                                             \begin{array}{cc}
                                                                                                                                                               0 & \beta_1 \\
                                                                                                                                                               \beta_2 & 0 \\
                                                                                                                                                             \end{array}
                                                                                                                                                           \right)
$ w.r.t. the decomposition, i.e., $\beta_1$ (resp. $\beta_2$) is a holomorphic section of $\textup{Hom}(E_2,E_1)\otimes \kappa$ (resp. $\textup{Hom}(E_1,E_2)\otimes \kappa$). If the decomposition $E=E_1\oplus E_2$ is orthogonal w.r.t. a Hermitian structure $h$ on $E$, the $\mathbb{U}(p,q)$-Higgs bundle is called \emph{harmonic} w.r.t. $h$ if the following Hitchin equations are satisfied:
\[R(E_1)+\beta_2^\ddag\beta_2+\beta_1\beta_1^\ddag=-i\mu Id_{E_1}\omega,\quad R(E_2)+\beta_2\beta_2^\ddag+\beta_1^\ddag\beta_1=-i\mu Id_{E_2}\omega,\]
where $R(E_i)$ is the curvature of the Chern connection in $E_i$, $\beta^\ddag$ the adjoint of $\beta$ w.r.t. the Hermitian structures on $E_1$ and $E_2$, $\omega$ a volume form on $M$ and $\mu$ a real constant\footnote{If $M$ is compact and $\omega$ normalized such that $\int_M \omega=2\pi$, $\mu$ is the slope of $E$, i.e., $\textup{deg }E/\textup{rk} E$. Here $\textup{deg }E$ is the integration of the first Chern class of $E$ over $M$.}. If for the harmonic $\mathbb{U}(p,q)$-Higgs bundle $(E_1,E_2, \beta_1, \beta_2,h)$, $\beta_1$ vanishes, we say it is a \emph{holomorphic} $\mathbb{U}(p,q)$-Higgs bundle w.r.t. $h$, which is relevant to our context. For more material on $\mathbb{U}(p,q)$-Higgs bundles, we refer the reader to \cite{bradlow2003surface, corlette1988flat, wells1980differential}.

For a simple symmetric operator $T$ with deficiency indices $(n_+,n_-)$, the characteristic vector bundle of the third kind can be viewed as a combination of $E(T)$ and $F(T)$.
\begin{definition}The characteristic vector bundle $G(T)$ of the third kind is the holomorphic vector bundle $E(T)|_{\mathbb{C}_+}\oplus F(T)|_{\mathbb{C}_+}$ over $\mathbb{C}_+$. Compared with \S\S~\ref{ssec1} and \S\S~\ref{ssec2}, the Hermitian metric on $E_\lambda(T)$ shall now be scaled by the factor $2\Im \lambda$ while that in $F_\lambda(T)$ shall be scaled by the factor $1/(2\Im \lambda)$.
\end{definition}
It seems that we haven't treated $E(T)|_{\mathbb{C}_\pm}$ on the same footing. However, by Riesz representation theorem, $F(T)|_{\mathbb{C}_+}$ can be identified with $F^\dag(T)|_{\mathbb{C}_+}$ (essentially $E(T)|_{\mathbb{C}_-}$). Since $F^\dag(T)|_{\mathbb{C}_+}$ is anti-holomorphic, if we insist on using $F^\dag(T)|_{\mathbb{C}_+}$ instead of $F(T)|_{\mathbb{C}_+}$, we should replace the complex structure on each fiber of $F^\dag(T)|_{\mathbb{C}_+}$ with the opposite one.

Now let $E_1=E(T)|_{\mathbb{C}_+}$ and $E_2=F(T)|_{\mathbb{C}_+}$. We do have a $\mathbb{U}(n_+,n_-)$-Higgs field here, which is intrinsically defined at $\lambda\in \mathbb{C}_+$ by $\beta_1(\lambda)=0$ and $\beta_2(\lambda):=K(\lambda,\bar{\lambda})d\lambda$ where $K$ is the reproducing kernel defined in \S\S~\ref{ssec2}. On $\mathbb{C}_+$ we continue to choose the Poincar$\acute{e}$ metric $\omega=\sigma$.
\begin{proposition}With the above $\beta_2$ and $\omega$, the holomorphic Hermitian vector bundle $G(T)$ is a holomorphic $\mathbb{U}(n_+,n_-)$-Higgs bundle with $\mu=0$.
\end{proposition}
\begin{proof}We only need to check that the Hitchin equations do hold in this case, but this is almost a reinterpretation of the curvature expression we have derived in \S\S~\ref{ssec1}.

Note that if a symplectic isomorphism $\Phi$ has been chosen, then according to our calculation in \S\S~\ref{ssec2}, $\beta_2(\lambda)=-iB'(\lambda)d\lambda$. In $E_\lambda(T)$, the Hermitian structure now is given by
\[(\gamma_+(\lambda)\varphi_1, \gamma_+(\lambda)\varphi_2)_{E_\lambda}=(\mathrm{K}\varphi_1, \varphi_2)_{H_+},\quad \forall \varphi_1, \varphi_2\in H_+\]
where $\mathrm{K}=Id-B^*B$ while the Hermitian structure in $F_\lambda(T)$ is given by
\[(\psi_-(\bar{\lambda})\varphi_1, \psi_-(\bar{\lambda})\varphi_2)_{F_\lambda}=(\tilde{\mathrm{K}}^{-1}\varphi_1,\varphi_2)_{H_-},\quad \forall \varphi_1, \varphi_2\in H_-\]
where $\tilde{\mathrm{K}}=Id-BB^*$. In terms of these, by definition
\[(\tilde{\mathrm{K}}^{-1}B'(\lambda)\varphi_1, \varphi_2)_{H_-}=(\mathrm{K}\varphi_1,B'(\lambda)^\ddag\varphi_2)_{H_+},\quad \varphi_1\in H_+, \quad \varphi_2\in H_-,\]
and consequently
\[B'(\lambda)^\ddag=\mathrm{K}^{-1}B'(\lambda)^*\tilde{\mathrm{K}}^{-1}.\]
Now the curvature $\hat{R}_T$ associated to the new metric on $E(T)|_{\mathbb{C}_+}$ reads
\[\hat{R}_T=\bar{\partial}(\mathrm{K}^{-1}\partial \mathrm{K})=\mathrm{K}^{-1}B'(\lambda)^*\tilde{\mathrm{K}}^{-1}B'(\lambda)d\lambda\wedge d\bar{\lambda}.\]
Then the first Hitchin equation with $\mu=0$ follows immediately. The equation for $F(T)|_{\mathbb{C}_+}$ can be checked similarly.
\end{proof}

We should point out that the converse of this proposition also holds, but the necessary conceptual preparation to prove it would deviate us from the main subject of the paper too far. The proof and some of its consequences will be spelled out elsewhere.

Thus if $(E_1, E_2, \beta_2, h)$ is a holomorphic $\mathbb{U}(n_+,n_-)$-Higgs bundle over $\mathbb{C}_+$ as in the above proposition, we can construct a vector bundle $E_2'$ over $\mathbb{C}_-$ by attaching at $\lambda\in \mathbb{C}_-$ the fiber $E_{2\bar{\lambda}}$. Then $E_1\cup E_2'^\dag$ is the characteristic vector bundle $E(T)$ for a simple symmetric operator $T$ with deficiency indices $(n_+,n_-)$. The relation of simple symmetric operators with holomorphic $\mathbb{U}(p,q)$-Higgs bundles over $\mathbb{C}_+$ should not be viewed as a coincidence. If we identify the operator $T$ with its characteristic vector bundle $G(T)$, then the correspondence in Thm.~\ref{thm12} can be regarded as a special version of the famous non-abelian Hodge correspondence.

To conclude this subsection and also provide part motivations for the next section, we shall describe the special Nevanlinna curves occurring in the study of holomorphic $\mathbb{U}(p,q)$-Higgs bundles, though such a name "Nevanlinna curve" was never used in the theory of Higgs bundles. Let $M$ be a closed Riemann surface with genus $g>1$ and $\pi_1(M)$ the fundamental group of $M$. Thus due to the uniformization theorem of Riemann surfaces, $\mathbb{C}_+$ serves as a universal cover of $M$ and $M$ can be viewed as the quotient of $\mathbb{C}_+$ under the action of $\pi_1(M)$ (then $\pi_1(M)$ is realized as a discrete subgroup of $\mathbb{PSL}(2,\mathbb{R})$). For the group $\mathbb{U}(p,q)$, people are concerned with representations of $\pi_1(M)$ into $\mathbb{PU}(p,q)$, i.e., group homomorphisms $\rho:\pi_1(M)\rightarrow \mathbb{PU}(p,q)$. With some conditions related to the Hitchin equations, associated to $\rho$ is a holomorphic map $\sigma$ from $\mathbb{C}_+$ to $\mathbb{PU}(p,q)/\mathbb{P}(\mathbb{U}(p)\times \mathbb{U}(q))$. This $\sigma$ is $\pi_1(M)$-equivariant in the sense that $\sigma(g\cdot \lambda)=\rho(g)\cdot \sigma(\lambda)$ for any $g\in \pi_1(M)$. Due to our previous argument on functional model, $\sigma$ should have its operator theoretic content concerning (projective) unitary representations of $\pi_1(M)$ in the reproducing kernel Hilbert space $\mathfrak{H}$. We shall pursue this untouched topic elsewhere.

\section{Spaces of congruence classes or weak congruence classes}\label{sec7}
It is now clear that the problem of unitary classification of simple symmetric operators with deficiency indices $(n_+,n_-)$ is equivalent to the problem of congruent classification of Nevanlinna curves of genus $(n_+,n_-)$. Perhaps the former is less hopeful partly because the unitary group of an infinite-dimensional Hilbert space is rather large. However, the latter is much more easy, at least when both $n_\pm$ are finite. A basic question raised by the correspondence in Thm.~\ref{thm12} is how the properties of a symmetric operator $T$ correspond to that of its Weyl class $[W_T(\lambda)]$.
\begin{definition}Let $\mathcal{N}(n_+,n_-)$ denote the set of congruence classes of Nevanlinna curves with genus $(n_+,n_-)$. Due to Thm.~\ref{thm12}, $\mathcal{N}(n_+,n_-)$ is called the moduli space of simple symmetric operators with deficiency indices $(n_+,n_-)$.
\end{definition}
 In applications when a symmetric operator appears, it often depends on continuous parameters. It is then natural to ask: how may the operator and its equivalence class vary if the parameters change, and is a certain property of the equivalence class of $T$ stable if the parameters are perturbed slightly? This means we can hardly view $\mathcal{N}(n_+,n_-)$ just as a set. Topological or even geometric structures must be put on $\mathcal{N}(n_+,n_-)$ to reflect the differences and connections among symmetric operators.

There are also other reasons to study $\mathcal{N}(n_+,n_-)$. We can just view a Nevanlinna curve as a holomorphic map from the hyperbolic plane to the hyperbolic space $W_+(H)$. Depending on how to realize the hyperbolic plane and $W_+(H)$, a Nevanlinna curve can be represented by a (matrix-valued) Nevanlinna function (also called a Herglotz function by some authors), a Schur function, or a Caratheodory function, etc. As objects in analysis and applied mathematics, they have actually been noticed and studied from different motivations for a long time, in spite of leaving the underlying Hermitian symmetric spaces totally neglected. As we have pointed out in \S~6, special Nevanlinna curves also appear in the guise of a special kind of Higgs bundles. In the area of several complex variables, people are also concerned with holomorphic maps between Hermitian symmetric spaces, which provide the basic examples of maps between complex hyperbolic spaces. Additionally in a series of papers starting with \cite{mok2012extension}, motivated by questions in arithmetic geometry, N. Mok investigated such maps that are isometries.

All these suggest a unified theory of Nevanlinna curves. For simplicity, in this section we only consider the case $n_\pm\in \mathbb{N}$. A thorough investigation is beyond the scope of the paper, and we only list some elementary observations.

We shall replace $\mathbb{C}_+$ by the unit disc $D$ with $0$ as its center. This is only for convenience and can be achieved by using Cayley transform. Let $\mathfrak{N}(n_+,n_-)$ be the set of Nevanlinna curves of genus $(n_+,n_-)$. We equip $\mathfrak{N}(n_+,n_-)$ with its compact-open topology, which can be metrizable. According to \cite[Thm.~3.2, Chap.~V]{kobayashi2005hyperbolic} $\mathfrak{N}(n_+,n_-)$ is locally compact. Fix a point $w_0\in W_+(H)$. Since $\mathbb{PU}(n_+,n_-)$ acts transitively on $W_+(H)$, in each congruence class there is a Nevanlinna curve $N(\lambda)$ such that $N(0)=w_0$. Consider
\[\mathfrak{N}_0(n_+,n_-):=\{N\in \mathfrak{N}(n_+,n_-)|N(0)=w_0\}.\]
$\mathfrak{N}_0(n_+,n_-)$ is certainly closed in $\mathfrak{N}(n_+,n_-)$ and still carries the residual continuous action of the isotropy subgroup $G_0$ of $w_0$, which is conjugate to $\mathbb{P}(\mathbb{U}(n_+)\times \mathbb{U}(n_-))$. Obviously, $\mathcal{N}(n_+,n_-)$ as a set can be identified with the orbit space $\mathfrak{N}_0(n_+,n_-)/G_0$. It's easy to see that the quotient topology on $\mathfrak{N}_0(n_+,n_-)/G_0$ is essentially independent of the choice of $w_0$.
\begin{theorem}\label{thm13}$\mathfrak{N}_0(n_+,n_-)/G_0$ with its quotient topology is Hausdorff, compact and contractible.
\end{theorem}
\begin{proof}Choose $N_+\in W_+(H)$ whose symplectic complement is denoted by $N_-$. In terms of these, $W_+(H)$ is realized as the open unit ball in $\mathbb{B}(N_+, N_-)$. W.l.g, we take $w_0=0\in \mathbb{B}(N_+, N_-)$. Then due to the uniform boundedness, $\mathfrak{N}_0(n_+,n_-)$ is a normal family according to Montel's theorem. Hence $\mathfrak{N}_0(n_+,n_-)$ is compact.
That $\mathfrak{N}_0(n_+,n_-)$ is Hausdorff is obvious. Since now the isotropy group $\mathbb{P}(\mathbb{U}(n_+)\times \mathbb{U}(n_-))$ is compact, the orbit space as a quotient is Hausdorff and compact.

In the above realization, $\mathfrak{N}_0(n_+,n_-)$ is contractible (because it's convex) and the map $(t, B(\cdot))\mapsto tB(\cdot)$ for $t\in [0,1]$ and $B(\cdot)\in \mathfrak{N}_0(n_+,n_-)$ gives a (null) homotopy between the identity and the constant map. By Prop.~\ref{p8}, now $(g_1,g_2)\in \mathbb{U}(n_+)\times \mathbb{U}(n_-)$ acts on $B(\lambda)$ as $g_2B(\lambda)g_1^{-1}$. The previous homotopy is equivariant w.r.t. this action and thus descends to a null homotopy in $\mathfrak{N}_0(n_+,n_-)/G_0$, i.e., $\mathfrak{N}_0(n_+,n_-)/G_0$ is contractible.
\end{proof}
\emph{Remark}. The compactness of $\mathfrak{N}_0(n_+,n_-)/G_0$ is a consequence of the hyperbolicity of $W_+(H)$. The restriction of hyperbolicity on holomorphic maps is a basic theme in hyperbolic geometry. For more information on this topic, see for example \cite{kobayashi2005hyperbolic}.
\begin{definition}$N(\cdot)\in \mathfrak{N}(n_+,n_-)$ is called generic, if the closure of $N(\mathbb{C}_+)$ in $Gr(n_+,H)$ is disjoint with the topological boundary of $W_+(H)$ in $Gr(n_+,H)$. Obviously, if $N(\cdot)$ is generic, so is any Nevalinna curve congruent to it; in this case, we say $[N(\lambda)]$ is generic. For a simple symmetric operator $T$, if $[W_T(\lambda)]$ is generic, we also say $T$ is generic.
\end{definition}
With the choice of $N_\pm$ as above, let $\mathfrak{B}$ be the space of bounded $\mathbb{B}(N_+,N_-)$-valued analytic functions on $D$. Then $\mathfrak{B}$ is a Banach space with its norm $|f|_{\mathfrak{B}}=\sup_{\lambda\in D}\|f(\lambda)\|$ for $f\in \mathfrak{B}$. Obviously, $\mathfrak{N}(n_+,n_-)\subset \mathfrak{B}$. $\mathfrak{N}(n_+,n_-)$ contains the open unit ball $\mathfrak{B}_1$ in $\mathfrak{B}$ and is contained in the closure $\overline{\mathfrak{B}_1}$, i.e., $\mathfrak{B}_1\subset \mathfrak{N}(n_+,n_-)\subset \overline{\mathfrak{B}_1}$. Then $\mathfrak{B}_1$ can be viewed as the subset of generic Nevanlinna curves. However, to some extent, generic Nevanlinna curves are not that interesting. To see this, assume that $n_+=n_-$. Then by investigation in the forthcoming \S~\ref{asing}, if $T$ is a simple symmetric operator whose Weyl curve is generic, then any self-adjoint extension of $T$ has the whole real line as its spectrum. These arguments lead to the following definition.
\begin{definition}Let $\mathcal{N}^*(n_+,n_-)\subset \mathcal{N}(n_+,n_-)$ be the space of congruence classes of Nevanlinna curves that are not generic. If $[N(\cdot)]\in \mathcal{N}^*(n_+,n_-)$, we say $N(\cdot)$ is singular.
\end{definition}
However, congruence classes in $\mathcal{N}^*(n_+,n_-)$ in a sense can produce all congruence classes in $\mathcal{N}(n_+,n_-)$: by scaling, any nonzero element in $\mathfrak{B}_1$ can be turned into a singular one, and vice versa. The complexity is that scaling seems to have no intrinsic meaning here.

Let \[\mathfrak{N}_0^*(n_+,n_-):=\{N(\cdot)\in \mathfrak{N}_0(n_+,n_-)|N(\cdot)\,\textup{is\, not\, generic}\}.\]
Then clearly, $\mathcal{N}^*(n_+,n_-)$ can be identified with $\mathfrak{N}_0^*(n_+,n_-)/\mathbb{P}(\mathbb{U}(n_+)\times \mathbb{U}(n_-))$. One can see fairly easily that $\mathfrak{N}_0^*(n_+,n_-)$ is Hausdorff and connected and so is $\mathcal{N}^*(n_+,n_-)$, but we don't know if $\mathcal{N}^*(n_+,n_-)$ is compact.

 In the spirit of Morse theory, one may consider functions on $\mathcal{N}(n_+,n_-)$ or $\mathcal{N}^*(n_+,n_-)$ to learn about these spaces. One candidate is the Chern functions $c_j$ defined in the previous section. Recall that these are functions defined on $\mathbb{C}_+$ and depend only on the congruence class by Thm.~\ref{thm12}. Then their values $C_j$ at $\lambda=i$ (if the unit disc is used, $\lambda=0$) can be viewed as functions defined on $\mathcal{N}(n_+,n_-)$ or $\mathcal{N}^*(n_+,n_-)$.
\begin{proposition}Each $C_j$ as a function on $\mathcal{N}(n_+,n_-)$ or $\mathcal{N}^*(n_+,n_-)$ is continuous.
\end{proposition}
\begin{proof}As before, we choose $N_+\in W_+(H)$ to realize each $N(\lambda)$ as a contractive operator-valued function $B(\lambda)$. We can first view $C_j$ as a function $C_j(B(\cdot))$ defined on $\mathfrak{N}_0(n_+,n_-)$. Note that $c_j$ is a polynomial in entries of $r_T$, involving derivatives of $B$ at most up to the first order. If a sequence $\{B_k(\cdot)\}\subset \mathfrak{N}_0(n_+,n_-)$ has $B(\cdot)\in \mathfrak{N}_0(n_+,n_-)$ as its limit w.r.t. the compact-open topology, then $B_k(\lambda)\rightarrow B(\lambda)$ and $B_k'(\lambda)\rightarrow B'(\lambda)$ uniformly on each compact subset in $\mathbb{C}_+$. Thus $c_j(B_k(\lambda))\rightarrow c_j(B(\lambda))$ in the same sense; in particular, we have $C_j(B_k(\cdot))\rightarrow C_j(B(\cdot))$ for evaluation at $\lambda=i$ is continuous. This shows $C_j$ is continuous on $\mathfrak{N}_0(n_+,n_-)$. Since by definition $C_j$ is constant on each $\mathbb{P}(\mathbb{U}(n_+)\times\mathbb{U}(n_-))$-orbit in $\mathfrak{N}_0(n_+,n_-)$, then $C_j$ descends to a continuous function on the quotient space $$\mathfrak{N}_0(n_+,n_-)/\mathbb{P}(\mathbb{U}(n_+)\times\mathbb{U}(n_-))=\mathcal{N}(n_+,n_-).$$
The continuity of $C_j$ on $\mathcal{N}^*(n_+,n_-)$ follows immediately.
\end{proof}
Due to Prop.~\ref{p12}, $C_1$ as a function on $\mathcal{N}(n_+,n_-)$ takes its maximum $n_+$ at a congruence class if and only if the class is that of a constant Nevanlinna curve. This fact may imply that $C_1$ can be used to detect the structure of $\mathcal{N}(n_+,n_-)$.
\begin{example}In Example \ref{ex4}, we have obtained
\[c_1(T_a)=n(1-(\frac{\Im \lambda}{\Im \lambda +a})^2).\]
Thus $C_1=n(1-(\frac{1}{1+a})^2)$. Note that each $a\geq 0$ parameterizes a congruence class of Nevanlinna curves and this expression can be viewed as $C_1$ restricted on this family of congruence classes. The example shows $C_1$ is not a trivial function.
\end{example}
Since the behaviors of different Nevanlinna curves can be rather different, it can be expected that these spaces $\mathcal{N}(n_+,n_-)$ and $\mathcal{N}^*(n_+,n_-)$ are highly singular. Our belief is that the structure of singularities of these spaces shall reflect the general pattern among simple symmetric operators.

Since $\mathbb{C}_+$ has $\mathbb{PSL}(2, \mathbb{R})$ as its automorphism group, if $N(\lambda)$ is a Nevanlinna curve, and $g\in\mathbb{PSL}(2, \mathbb{R})$, then $N(g\cdot \lambda)$ is again a Nevanlinna curve, which may be not in the same congruence class as $N(\lambda)$. This simply means the space $\mathfrak{N}(n_+,n_-)$ actually carries the natural action of the direct product group $\mathbb{PSL}(2, \mathbb{R})\times \mathbb{PU}(n_+,n_-)$.
\begin{definition}If $N_1(\cdot), N_2(\cdot)\in \mathfrak{N}(n_+,n_-)$ are in the same $\mathbb{PSL}(2, \mathbb{R})\times \mathbb{PU}(n_+,n_-)$-orbit, we say they are weakly congruent. The weak congruence class of $N(\cdot)$ will be denoted by $[N(\lambda)]_w$. The set of weak congruence classes in $\mathfrak{N}(n_+,n_-)$ is denoted by $\mathcal{N}_w(n_+,n_-)$.
\end{definition}

Then what's the meaning of this reparameterization using $\mathbb{PSL}(2,\mathbb{R})$ for simple symmetric operators?

Recall that a symmetric operator $T$ determines an isotropic subspace $A_T$ (its graph) in the standard strong symplectic Hilbert space $(H\oplus_\bot H, [\cdot, \cdot]_{new})$. For each $g=\left(
                                                                                       \begin{array}{cc}
                                                                                         a & b \\
                                                                                         c & d \\
                                                                                       \end{array}
                                                                                     \right)\in \mathbb{SL}(2,\mathbb{R})$, we can define its action on $(H\oplus_\bot H, [\cdot, \cdot]_{new})$ by \[\hat{g}(x,y)=(cy+dx, ay+bx)\] for $(x,y)\in H\oplus_\bot H$. $\hat{g}$ is obviously a symplectic isomorphism and thus transforms $A_T$ into another isotropic subspace $\hat{g}(A_T)$. Notice that if $g'=-g$, then $\widehat{g'}(A_T)=\hat{g}(A_T)$.
\begin{proposition} If $T$ is simple and $g\in \mathbb{SL}(2,\mathbb{R})$, then $\hat{g}(A_T)$ is also the graph of a simple symmetric operator.
\end{proposition}
\begin{proof}We first prove $\hat{g}(A_T)$ is really the graph of an operator. Otherwise, we can find $0\neq x\in D(T)$ such that
\[cTx+dx=0,\quad aTx+bx\neq 0. \]
  $c$ cannot be zero and hence $x$ has to be an eigenvector of $T$. This is impossible for $T$ is simple. Denote this operator by $g(T)$. $g(T)$ is certainly closed and that $\hat{g}(A_T)$ is isotropic means precisely $g(T)$ is symmetric. It is easy to see that if $x\in \textup{ker}(T^*-\lambda)$ for $\lambda\in \mathbb{C}_+\cup\mathbb{C}_-$, then $x\in \textup{ker}(g(T)^*-\frac{a\lambda+b}{c\lambda+d})$, and vice versa. This shows if $T$ is simple, so is $g(T)$.
\end{proof}
Obviously, $g(T)$ essentially makes sense for $g\in \mathbb{PSL}(2,\mathbb{R})$ and thus we call the symmetric operator $g(T)$ a M$\ddot{o}$bius transform of $T$. It has the same deficiency indices as $T$ does. The above proposition means simplicity is weakly congruently invariant.
\begin{proposition}\label{p13}If $B(\lambda)$ is a contractive Weyl function of $T$ w.r.t. a symplectic isomorphism $\Phi$, then $B(g^{-1}\cdot\lambda)$ is a contractive Weyl function of $g(T)$.
\end{proposition}
\begin{proof}This is clear from the above fact that $\textup{ker}(T^*-\lambda)=\textup{ker}(g(T)^*-\frac{a\lambda+b}{c\lambda+d})$.
\end{proof}
Thus if simple symmetric operators $T$ and $T'$ are weakly congruent, their extension theories are essentially connected with each other simply by the transform $\hat{g}$ for some $g\in \mathbb{PSL}(2,\mathbb{R})$. In this sense, $T$ and $T'$ should be reasonably identified though they may not be unitarily equivalent.

Since $\mathbb{PSL}(2,\mathbb{R})$ acts transitively on the unit disc $D$, each weak congruence class has a representative $N(\lambda)$ such that $N(0)=w_0$. Notice that the isotropy group of $0\in D$ is $\mathbb{U}(1)$. If again $N_+\in W_+(H)$ is fixed, $\mathcal{N}_w(n_+,n_-)$ can be identified with the orbit space $\mathfrak{N}_0(n_+,n_-)/(\mathbb{U}(1)\times \mathbb{P}(\mathbb{U}(n_+)\times \mathbb{U}(n_-)))$. Similar to Thm.~\ref{thm13}, we have
\begin{theorem}The space $\mathcal{N}_w(n_+,n_-)$ is Hausdorff, compact and contractible.
\end{theorem}
\begin{proof}It's similar to the proof of Thm.~\ref{thm13} and so omitted.
\end{proof}
In this case, a \emph{generic} weak congruence class still makes sense and the space $\mathcal{N}_w^*(n_+,n_-)$ can as well be defined.
\section{Symmetric operators with symmetries}\label{sec8}
One way the moduli space $\mathcal{N}(n_+,n_-)$ or $\mathcal{N}_w(n_+,n_-)$ of simple symmetric operators acquires a singularity at the congruence class $[W_T(\lambda)]$ or weak congruence class $[W_T(\lambda)]_w$ is that $T$ may have nontrivial symmetries. This is the lesson we have learned from the theory of moduli spaces of compact Riemann surfaces, where the nontrivial but finite automorphism group of a Riemman surface causes a mild singularity at the corresponding point in the moduli space. Meanwhile, symmetric operators with symmetries also occur frequently in applications, e.g., the Euclidean Laplacian $\triangle$ of the unit ball in $\mathbb{R}^n$ is rotation invariant. The goal of this section is to consider two ways in which symmetries may arise in the investigation of simple symmetric operators and the focus is on how the symmetries affect the Weyl curves. Note that unless otherwise stated we won't assume $n_\pm<+\infty$ in this section.

Let $T$ be a closed operator in a Hilbert space $H$.
\begin{definition}$U\in \mathbb{U}(H)$ is called a symmetry of a closed operator $T$, if $D(T)$ is $U$-invariant and for any $x\in D(T)$, $TUx=UTx$.
\end{definition}
All symmetries of $T$ form a subgroup $\mathbb{U}_T\subset\mathbb{U}(H)$. A symmetry $U$ of $T$ is called trivial if it is $c\times Id$ where $c\in \mathbb{U}(1)$ is a constant. In this sense, $\mathbb{U}(1)\subset \mathbb{U}_T$. If the simple symmetric operator $T$ has a symmetry $U$, then $U$ is also a symmetry of $T^*$. Consequently, this induces a representation $\rho$ of $\mathbb{U}_T$ on $\mathcal{B}_T$. Besides, for each $\lambda\in \mathbb{C}_+\cup \mathbb{C}_-$, $\textup{ker}(T^*-\lambda)$ is $U$-invariant. In particular, $\rho(U)$ fixes the Weyl curve \emph{pointwise}.
\begin{proposition}For $U\in \mathbb{U}_T$, $\rho(U)$ preserves the Hilbert space structure and the strong symplectic structure on $\mathcal{B}_T$. If $\rho(U)=c\times Id$ for a constant $c$, then $U$ is trivial.
\end{proposition}
\begin{proof}Note that $U$ also preserves the graph inner product on $D(T^*)$, and $D(T)$ is $U$-invariant and closed in $D(T^*)$ w.r.t. this inner product. Hence the Hilbert space structure on $\mathcal{B}_T$ is preserved by $\rho(U)$.

For $x,y\in D(T^*)$,
\begin{eqnarray*}[\rho(U)[x],\rho(U)[y]]_T&=&[[Ux],[Uy]]_T=(T^*Ux,Uy)_H-(Ux,T^*Uy)_H\\
&=&(UT^*x,Uy)_H-(Ux,UT^*y)_H=(T^*x,y)_H-(x,T^*y)_H\\
&=&[[x],[y]]_T.
\end{eqnarray*}

The above implies that $\rho(U)$ preserves $\textup{ker}(\mathcal{J}\mp i)\subset \mathcal{B}_T$ where $\mathcal{J}$ is the operator constructed below Def.~\ref{de1} and thus $\rho(U)$ lies in $\mathbb{U}(n_+)\times \mathbb{U}(n_-)$. Thus if $\rho(U)=c\times Id$, then $c\in \mathbb{U}(1)$. Consequently for $\lambda\in \mathbb{C}_+\cup \mathbb{C}_-$ and $x\in \textup{ker}(T^*-\lambda)$,
\[[Ux]=\rho(U)[x]=c[x]=[cx].\]
We should have $Ux-cx\in \textup{ker}(T^*-\lambda)\cap D(T)=0$, i.e., $Ux=cx$. Due to simplicity of $T$, we must have $U=c\times Id$.
\end{proof}
\emph{Remark}. The second claim has the following consequence: for a simple symmetric operator $T$ with finite deficiency indices, $\mathbb{U}_T$ is always a compact Lie group. Thus even $\mathbb{U}(H)$ is very large, simplicity of $T$ requires that $T$ cannot have too many symmetries.

 If $G$ is a subgroup of $\mathbb{U}_T$, a symplectic isomorphism $\Phi$ between $\mathcal{B}_T$ and the standard strong symplectic Hilbert space $\mathbb{H}=H_+\oplus_\bot H_-$ is called $G$-compatible, if there is a group homomorphism
$\tau: G\rightarrow \mathbb{U}(n_+,n_-)$ such that
\[\tau(g)\cdot \Phi([x])=\Phi(\rho(g)\cdot [x]),\quad \forall x\in D(T^*).\]
Note that $\tau(G)$ should lie in a subgroup conjugate to $\mathbb{U}(n_+)\times \mathbb{U}(n_-)$.
\begin{corollary}Let $G\subset \mathbb{U}_T$ be a subgroup of symmetries of $T$. If a $G$-compatible symplectic isomorphism between $\mathcal{B}_T$ and $\mathbb{H}=H_+\oplus_\bot H_-$ as above has been chosen and $B(\lambda)$ is the corresponding contractive Weyl function, then
\[\tau (g)\cdot B(\lambda)=B(\lambda),\quad \forall g\in G, \quad \forall \lambda\in \mathbb{C}_+.\]
\end{corollary}
\begin{proof}This can be checked easily and so the proof is omitted.
\end{proof}
This of course means that symmetries restrict the possible form of the contractive Weyl function; in particular, if the image of $\tau(G)$ in $\mathbb{PU}(n_+,n_-)$ is nontrivial, then $W_T(\cdot)\in \mathfrak{N}(n_+,n_-)$ has a nontrivial isotropy subgroup. Conversely we have
\begin{proposition}If $W_T(\cdot)$ has a nontrivial isotropy subgroup in $\mathfrak{N}(n_+,n_-)$ under the action of $\mathbb{PU}(n_+,n_-)$, then $T$ has a nontrivial symmetry.
\end{proposition}
\begin{proof}W.l.g, we can assume $T$ to be the multiplication operator in our functional model. Let $g=(g_1, g_2)\in \mathbb{U}(n_+)\times \mathbb{U}(n_-)$ fix the contractive Weyl function $B(\lambda)$ and the image of $g$ in $\mathbb{P}(\mathbb{U}(n_+)\times \mathbb{U}(n_-))$ be nontrivial. By Prop.~\ref{p8}, we see $g_2B(\lambda)=B(\lambda)g_1$ for all $\lambda\in \mathbb{C}_+$. This identity can be used to produce a nontrivial gauge transformation of the characteristic vector bundle $F(T)$ and consequently a unitary operator on the reproducing kernel Hilbert space $\mathfrak{H}$. This unitary operator is a symmetry of $T$ simply because the gauge transformation is linear fiberwise.
\end{proof}

We shall now turn to another kind of symmetries the symmetric operator $T$ may have. For simplicity, we should introduce the following definition.
\begin{definition}If a simple symmetric operator $T$ has only trivial symmetries, we say $T$ is irreducible.
\end{definition}
Irreducible simple symmetric operators shall be viewed as the building blocks of general simple symmetric operators.
\begin{proposition}A simple symmetric operator $T$ in $H$ is irreducible, if and only if it can't be written as the orthogonal direct sum of two simple symmetric operators.
\end{proposition}
\begin{proof}If $H=H_1\oplus_\bot H_2$ and $T=T_1\oplus_\bot T_2$ such that $T_1$ (resp. $T_2$) is a symmetric operator in $H_1$ (resp. $H_2$). For any $x=(x_1,x_2)\in H_1\oplus_\bot H_2$, set $Ux=(x_1,-x_2)$. Then clearly $U$ is a symmetry of $T$, which is not a multiple of $Id$. Thus $\mathbb{U}_T\neq \mathbb{U}(1)$.

Conversely, if $U$ is a nontrivial symmetry of $T$, $U$ also induces a unitary action on $H\oplus_\bot H$ via $\hat{U}(y_1,y_2)=(Uy_1,Uy_2)$. Clearly, $\hat{U}$ preserves both the Hilbert space structure and strong symplectic structure $[\cdot, \cdot]_{new}$ on $H\oplus_\bot H$. By definition, this action also preserves the graph $A_T$ of $T$ and $A_T^{\bot_s}$.

Since $U$ is not a multiple of $Id$, the action of $\hat{U}$ on $A_T^{\bot_s}$ is also not by multiplying by a constant $c$. Otherwise, $U$ acts on each $\textup{ker}(T^*-\lambda)$ for $\lambda\in \mathbb{C}_+\cup \mathbb{C}_-$ by scaling by the same constant $c$ and simplicity of $T$ then implies that $U=c\times Id$ on $H$. Using the spectral decomposition of $\hat{U}$ on $A_T^{\bot_s}$, we can find two nontrivial $\hat{U}$-invariant closed subspaces $S_1$ and $S_2$ in $A_T^{\bot_s}$ such that $A_T^{\bot_s}=S_1\oplus_\bot S_2$. Due to the special form of $\hat{U}$, each $S_i$ is a closed subspace of $H_i\oplus_\bot H_i$ for a certain closed subspace $H_i\subset H$ such that $H=H_1\oplus_\bot H_2$. Let $A_i=A_T\cap S_i$. Then it can be proved that $A_T=A_1\oplus_\bot A_2$. Then $A_i$ is the graph of a symmetric operator $T_i$ in $H_i$. Both $T_1$ and $T_2$ should be simple just because $T$ is. It should be pointed out that one of $A_1$ and $A_2$ can be 0. For example, if $A_2=0$, then $T_2$ is the zero symmetric operator in $H_2$.
\end{proof}
\emph{Remark}. By this proposition, if $\mathbb{U}_T\neq \mathbb{U}(1)$, $T$ is the orthogonal direct sum of two simple symmetric operators $T_1$ and $T_2$; If furthermore neither deficiency index of $T_i$ ($i=1,2$) is zero, then by choosing a suitable symplectic isomorphism, the contractive Weyl function can be written as a block-diagonal operator-valued function.

To summarize, the following three statements are equivalent for a simple symmetric operator: (1) $T$ has only trivial symmetries; (2) $W_T(\cdot)\in \mathfrak{N}(n_+,n_-)$ has trivial isotropy subgroup under the action of $\mathbb{PU}(n_+,n_-)$; (3) $T$ cannot be written as the orthogonal direct sum of two simple symmetric operators.

Motivated by the investigation on homogeneous Cowen-Douglas operators in \cite{koranyi2011classification}, we have the following definition.
\begin{definition}\label{de2}A simple symmetric operator $T$ is called homogeneous if for any $g\in \mathbb{PSL}(2,\mathbb{R})$,  $g(T)$ is unitarily equivalent to $T$, i.e., there is a unitary operator $U$ such that $UD(T)=D(g(T))$ and for each $x\in D(T)$, $g(T)Ux=UTx$.
\end{definition}
If $U$ in Def.\ref{de2} exists, it surely is not unique, but if $T$ is irreducible, $U$ is unique up to a constant factor $u\in \mathbb{U}(1)$.
\begin{proposition}Let $T$ be a simple symmetric operator, $B(\lambda)$ its contractive Weyl function w.r.t. a certain symplectic isomorphism $\Phi$ and $g\in \mathbb{PSL}(2,\mathbb{R})$. Then $g(T)$ is unitarily equivalent to $T$ if and only if there exists $\tilde{g}\in \mathbb{PU}(n_+,n_-)$ such that $B(g^{-1}\cdot \lambda)=\tilde{g}\cdot B(\lambda)$ for all $\lambda\in \mathbb{C}_+$.
\end{proposition}
\begin{proof}This is clear from Thm.~\ref{thm12} and Prop.~\ref{p13}.
\end{proof}
\begin{example}By the proposition, the symmetric operators in Example \ref{ex5} are homogenous; If the deficiency indices are $(1,1)$, they are irreducible. It's interesting to know whether there are other nontrivial (irreducible) examples.
\end{example}
\begin{proposition}Let $T$ be an irreducible simple symmetric operator and $G$ a subgroup of $\mathbb{PSL}(2,\mathbb{R})$. If $g(T)$ is unitarily equivalent to $T$ for any $g\in G$, then associated with $T$ is a conjugacy class of homomorphisms $\tau: G\rightarrow \mathbb{PU}(n_+,n_-)$.
\end{proposition}
\begin{proof}Fix a boundary triplet for $T$ and let $B(\lambda)$ be the corresponding contractive Weyl function. Then due to irreducibility of $T$, for each $g\in G$ there is a unique $\tau(g)\in \mathbb{PU}(n_+,n_-)$ such that $B(g\cdot \lambda)=\tau(g)\cdot B(\lambda)$ for all $\lambda\in \mathbb{C}_+$. It is easy to check that $\tau$ is a group homomorphism from $G$ to $\mathbb{PU}(n_+,n_-)$. A different boundary triplet will produce a new homomorphism conjugate to $\tau$.
\end{proof}

Recall that a subgroup $G\subset \mathbb{PSL}(2,\mathbb{R})$ is called discrete if the subspace topology of $G$ in $\mathbb{PSL}(2,\mathbb{R})$ is discrete.
\begin{definition}Let $G$ be a discrete subgroup of $\mathbb{PSL}(2,\mathbb{R})$. An irreducible simple symmetric operator $T$ is called $G$-automorphic if for each $g\in G$, $g(T)$ is unitarily equivalent to $T$.
\end{definition}
\begin{example}As we have noticed in \S\S~\ref{sec6}, associated with a holomorphic $\mathbb{U}(p,q)$-Higgs bundle over a compact Rieamann surface $M$ with genus $g>1$ is a $\pi_1(M)$-equivariant Nevanlinna curve. In the Fuchsian model of $M$, $\pi_1(M)$ is realized as a discrete subgroup of $\mathbb{PSL}(2,\mathbb{R})$. Then the multiplication operator in the associated functional model is necessarily $\pi_1(M)$-automorphic.

We should also point out that symmetric operators with symmetries associated to specific subgroups of $\mathbb{PSL}(2,\mathbb{R})$ have been noticed in the literature before, see \cite{bekker2013selfadjoint} for instance and references therein.
\end{example}
\section{Contractive Weyl function and spectral analysis}\label{sec9}
From now on, we only consider symmetric operators $T$ with deficiency indices $(n,n)$. This case is of particular interest in spectral theory. For later convenience, we shall just use $Gr(n,2n)$ to denote the Grassmannian of closed subspaces of dimension $n$ in a strong symplectic Hilbert space of dimension $2n$. This can be the intrinsic $Gr(n, \mathcal{B}_T)$ or $Gr(n, G\oplus_\bot G)$ if a specific boundary triplet is chosen. In particular, if $n=+\infty$, $Gr(n,2n)$ consists of closed subspaces of infinite dimension and co-dimension.

For simplicity right now, let's first assume $n\in \mathbb{N}$. Then all closed extensions of $T$ are parameterized by linear subspaces of $\mathbb{C}^{2n}$. There are several topological components of this parameter space indexed by the dimension of the subspaces parameterized by them. These are $Gr(r, 2n)$ ($0\leq r\leq 2n$), i.e., the Grassmannian of $r$-dimensional subspaces in $\mathbb{C}^{2n}$. The complex dimension of $Gr(r, 2n)$ is $r(2n-r)$, taking its maximum $n^2$ when $r=n$. In this sense, we say abstract boundary conditions parameterized by points in $Gr(n, 2n)$ are \emph{generic}.

 Another fundamental fact is that the two-branched Weyl curve now lies in $Gr(n, 2n)$. This is in sharp contrast to the case $n_+\neq n_-$, because in the latter case the two branches lie in $Gr(n_+, n_++n_-)$ and $Gr(n_-,n_++n_-)$ separately. Though these two manifolds are of the same dimension $n_+\times n_-$, they are different topological components and generally not of the maximal dimension among all topological components. Consequently, a special phenomenon may occur when the indices are $(n,n)$: the two-branched Weyl curve may admit an analytic continuation across part of the real line and the two branches then touch each other there. Due to Prop.~\ref{p11}, this is possible. In the best case, the Weyl curve may even be analytic along the whole real line.

 In the literature, the Weyl function $M(\lambda)$ is often used in spectral analysis of extensions of $T$. However, a disadvantage of $M(\lambda)$ is that, the spectrum of $T_0$ generally appears as singularities of $M(\lambda)$. If $T_0$ happens to have a discrete spectrum, then this means that the Weyl curve admits an analytic continuation across the real line except for at most countably infinite points. Even better, the singularities on $\mathbb{R}$ may not be geometric, i.e., they are there only because a "bad" coordinate chart is in use. To deal with problems concerning spectra, the contractive Weyl function $B(\lambda)$ may sometimes provide a better choice. Thus the basic goal of this section is spectral analysis in terms of the contractive Weyl function $B(\lambda)$ rather than the Weyl function $M(\lambda)$. Of course, inevitably the presentation in this section has some overlaps with the existing literature, but our emphasis is on the geometric formulation and the usefulness of $B(\lambda)$.

Recall that $\lambda\in \mathbb{C}$ is of regular type for a symmetric operator $T$ if there is a constant $c_\lambda>0$ such that
\begin{equation}\|(T-\lambda)x\|\geq c_\lambda \|x\|\label{eq2}\end{equation}
for any $x\in D(T)$. Denote the set of $\lambda$ of regular type for $T$ by $\Theta(T)$. It is open on $\mathbb{C}$ and if further $\Theta(T)=\mathbb{C}$, $T$ is called a \emph{regular} operator. Regular operators are necessarily simple.

From the above inequality (\ref{eq2}), $\textup{Ran}(T-\lambda)$ has to be closed for $\lambda\in\Theta(T)$. It is well-known that $\mathbb{C}_+\cup \mathbb{C}_-\subset \Theta(T)$. For $\lambda\in \Theta(T)$, we denote the image of $\textup{ker}(T^*-\lambda)$ in $\mathcal{B}_T$ by $\mathrm{M}(\lambda)$.
\begin{lemma}If $\lambda_0\in \Theta(T)\cap\mathbb{R} $, then $\mathrm{M}(\lambda_0)$ is Lagrangian in $\mathcal{B}_T$.
\end{lemma}
\begin{proof} Since for any $x,y\in \textup{ker}(T^*-\lambda_0)$,
\[(T^*x,y)-(x,T^*y)=\lambda_0(x,y)-\lambda_0(x,y)=0.\]
This shows $\mathrm{M}(\lambda_0)$ is isotropic, i.e., $\mathrm{M}(\lambda_0)\subset \mathrm{M}(\lambda_0)^{\bot_s}$. If $\check{y}\in \mathrm{M}(\lambda_0)^{\bot_s}$, let $y$ be a pre-image of $\check{y}$ in $D(T^*)$, then by definition for all $x \in \textup{ker}(T^*-\lambda_0)$
\[0=(T^*x,y)-(x,T^*y)=\lambda_0(x,y)-(x,T^*y)=(x,(\lambda_0-T^*)y).\]
Thus, $(\lambda_0-T^*)y\in (\textup{ker}(T^*-\lambda_0))^\bot=\textup{Ran}(T-\lambda_0)$. Consequently, there is a $y_0\in D(T)$ such that
\[T^*y-\lambda_0y=Ty_0-\lambda_0y_0.\]
This shows that $y-y_0\in \textup{ker}(T^*-\lambda_0)$. Note that the image of $y-y_0$ in $\mathcal{B}_T$ is again $\check{y}$. This implies that $\check{y}\in \mathrm{M}(\lambda_0)$ and hence $\mathrm{M}(\lambda_0)^{\bot_s}\subset \mathrm{M}(\lambda_0)$. The proof is completed.
\end{proof}
\begin{lemma}\label{lem1}Let $\tilde{T}$ be a closed extension of $T$, and $L$ the image of $D(\tilde{T})$ in $\mathcal{B}_T$. Then for $\lambda_0\in \Theta(T)$, $\lambda_0\in \rho(\tilde{T})$ if and only if $L$ and $\mathrm{M}(\lambda_0)$ are transversal, i.e.,
\[L\cap \mathrm{M}(\lambda_0)=0,\quad \textup{and}\quad L\oplus \mathrm{M}(\lambda_0)=\mathcal{B}_T.\]
\end{lemma}
\begin{proof}"$\Rightarrow$". If $L\cap \mathrm{M}(\lambda_0)\neq 0$, then there is a nonzero $x=\textup{ker}(T^*-\lambda_0)\cap D(\tilde{T})$, i.e., $\lambda_0$ is an eigenvalue of $\tilde{T}$. This is a contradiction. For any $x\in D(T^*)$, note that
\[x=(\tilde{T}-\lambda_0)^{-1}(T^*-\lambda_0)x+[x-(\tilde{T}-\lambda_0)^{-1}(T^*-\lambda_0)x]\]
and
\[(T^*-\lambda_0)(x-(\tilde{T}-\lambda_0)^{-1}(T^*-\lambda_0)x)=(T^*-\lambda_0)x-(T^*-\lambda_0)x=0.\]
So $D(T^*)=D(\tilde{T})\oplus \textup{ker}(T^*-\lambda_0)$, implying that $\mathcal{B}_T=L\oplus \mathrm{M}(\lambda_0)$.

"$\Leftarrow$". $L\cap \mathrm{M}(\lambda_0)=0$ implies that $\tilde{T}-\lambda_0$ is injective and $\lambda_0$ is not an eigenvalue of $\tilde{T}$. It suffices to prove $\textup{Ran}(\tilde{T}-\lambda_0)=H$, because the closed graph theorem shall imply that $\lambda_0\in \rho(\tilde{T})$. Since $\textup{Ran}(T-\lambda_0)$ is closed, so is $\textup{Ran}(T-\overline{\lambda_0})$ no matter $\lambda_0\in \mathbb{R}$ or not. Thus $\textup{Ran}(T^*-\lambda_0)$ is closed due to the closed range theorem. We know that $\textup{Ran}(T^*-\lambda_0)=H$ because $\textup{ker}(T-\overline{\lambda_0})=0$ (for $T$ is simple).

For any $y\in H$, choose $x$ such that $(T^*-\lambda_0)x=y$. Since $ L\oplus \mathrm{M}(\lambda_0)=\mathcal{B}_T$ , we can choose $x_1\in D(\tilde{T})$ and $x_2\in \textup{ker}(T^*-\lambda_0)$ such that $x=x_1+x_2$. Therefore,
\[(T^*-\lambda_0)x=(T^*-\lambda_0)x_1=(\tilde{T}-\lambda_0)x_1=y.\]
This shows $\textup{Ran}(\tilde{T}-\lambda_0)=H$ and the proof is finished.
\end{proof}
\begin{theorem}\label{t1}
If $\lambda_0\in \Theta(T)\cap\mathbb{R}$, then the two-branched Weyl curve $W_T(\lambda)$ continues analytically around $\lambda_0$.
\end{theorem}
The theorem will be proved in terms of the contractive Weyl function $B(\lambda)$. This works because different symplectic isomorphisms only introduce holomorphic isomorphisms among the ambient Grassmanianns. Let a symplectic isomorphism $\Phi$ be fixed. We begin with an extension of Prop.~\ref{res}. Recall that $T_\pm$ are extensions of $T$ determined by the boundary conditions $\Gamma_\pm x=0$ respectively.
\begin{proposition}If $\lambda_0\in \Theta(T)\cap\mathbb{R}$, then $\lambda_0\in \rho(T_\pm)$.
\end{proposition}
\begin{proof}We only prove the result for $T_-$. Note that $\mathrm{M}(\lambda_0)$ is Lagrangian in $\mathcal{B}_T$ and the image of $D(T_-)$ in $\mathcal{B}_T$ is a maximal completely positive-definite subspace. Due to Prop.~\ref{p3} and Lemma \ref{lem1}, the conclusion follows.
\end{proof}
\begin{proof}\emph{Proof of Thm.~\ref{t1}}. Since $\Theta(T)$ is open, if $\lambda_0\in \mathbb{R}$ is of regular type, so are points in a certain open interval $I\subset \mathbb{R}$ containing $\lambda_0$. Since $\rho(T_+)$ is open as well, it can be easily seen that Prop.~\ref{p1} continues to hold for points in $I$ and points on $\mathbb{C}_-$ very near $I$. Then the map $\gamma_+(\lambda)$ can be extended to an enlarged set $K\supset \mathbb{C}_+$ such that $K$ contains an open disc centered at $\lambda_0$. Therefore Prop.~\ref{p2} and Corol.~\ref{c1} still hold for $\lambda\in K$. Similar results hold for $\gamma_-$. This shows that the Weyl curve $W_T(\lambda)$ admits an analytic continuation around $\lambda_0$ and the proof of Thm.~\ref{t1} is thus completed.
\end{proof}
 As a consequence of Thm.~\ref{t1}, the characteristic vector bundles $E(T)$ and $F(T)$ are now both well-defined over $\Theta(T)$. If furthermore $T$ is regular, these are holomorphic vector bundles over $\mathbb{C}$. Since $\mathbb{C}$ is a Stein manifold, due to the famous Oka-Grauert principle \cite[Thm.~9.5]{leiterer1990holomorphic}, $E(T)$ and $F(T)$ are both holomorphically trivial and in particular the reproducing kernel Hilbert space $\mathfrak{H}$ is essentially a space of (vector-valued) entire functions. Indeed it has been known for several decades that some regular simple symmetric operators with deficiency indices $(1,1)$ have the multiplication operators in certain de Branges spaces of entire functions as their models. We refer the reader to \cite{silva2013branges} for a modern treatment of this connection.

Now for $\lambda_0\in \Theta(T)\cap\mathbb{R}$, we can find a disc $D$ with its center $\lambda_0$, sufficiently small such that both $\gamma_\pm$ are well-defined over $D$. Then for $\lambda\in D$, the image $\mathrm{M}(\lambda)$ of $\textup{ker}(T^*-\lambda)$ in the standard strong symplectic Hilbert space has two kinds of descriptions: \[\{(\phi,B(\lambda)\phi)\in G\oplus_\bot G|\varphi\in G\}\quad \textup{and}\quad\{(B(\bar{\lambda})^*\phi,\phi)\in G\oplus_\bot G|\phi\in  G\}.\]
On the disc $D$, it is necessary that
\[B(\lambda)=(B(\bar{\lambda})^*)^{-1}.\]
 Note that for $\lambda\in \Theta(T)\cap\mathbb{R}$, $B(\lambda)$ is a unitary operator, just as expected.

Recall that a pair $(V_1,V_2)$ of closed subspaces in a Hilbert space $H$ is called a Fredholm pair if $V_1\cap V_1$ is of finite dimension and $V_1+V_2$ is closed and of finite codimension in $H$. For a Fredholm pair $(V_1,V_2)$ in $H$, $$\textup{ind}(V_1,V_2):=\textup{dim}(V_1\cap V_2)-\textup{codim}(V_1+V_2)$$ is called the index of the pair.
\begin{theorem}\label{thm4}If $\tilde{T}$ is a closed extension of $T$ and $L$ is the image of $D(\tilde{T})$ in $\mathcal{B}_T$, then for $\lambda\in \Theta(T)$,\\
i) $\lambda\in \sigma_p(\tilde{T})$ if and only if $L\cap \mathrm{M}(\lambda)\neq 0$. If it is the case, \[\dim\textup{ker}(\tilde{T}-\lambda)=\dim (L\cap \mathrm{M}(\lambda)).\]
ii) $\lambda\in \sigma_c(\tilde{T})$ if and only if $L\cap \mathrm{M}(\lambda)=0$ and $L+\mathrm{M}(\lambda)$ is a proper dense subspace of $\mathcal{B}_T$.\\
iii) $\lambda\in \sigma_r(\tilde{T})$ if and only if $L\cap \mathrm{M}(\lambda)=0$ and $\overline{L+\mathrm{M}(\lambda)}\neq \mathcal{B}_T$.\\
iv) $\tilde{T}-\lambda$ is a Fredholm operator if and only if $(L,\mathrm{M}(\lambda))$ is a Fredholm pair in $\mathcal{B}_T$. If it is the case, then $\textup{ind}\tilde{T}=\textup{ind}(L,\mathrm{M}(\lambda))$.
\end{theorem}
\begin{proof}i) The first statement is clear. For the second, note that \[\textup{ker}(\tilde{T}-\lambda)=D(\tilde{T})\cap \textup{ker}(T^*-\lambda),\quad \textup{ker}(T^*-\lambda)\cap D(T)=0.\]
These show that the quotient map from $\textup{ker}(\tilde{T}-\lambda)\subset D(T^*)$ to $\mathcal{B}_T$ is injective and the image of $\textup{ker}(\tilde{T}-\lambda)$ is precisely $L\cap \mathrm{M}(\lambda)$. Then the claim follows.

ii) Due to the trichotomy of spectra, this will follow if we can prove iii).

iii) $\lambda\in \sigma_r(\tilde{T})$ if and only if $\lambda\not\in \sigma_p(\tilde{T})$ and $\bar{\lambda}\in \sigma_p(\tilde{T}^*)$. Note that $\tilde{T}^*$ is itself a closed extension of $T$ and the image of $D(\tilde{T}^*)$ in $\mathcal{B}_T$ is precisely $L^{\bot_s}$. Due to i), $\lambda\in \sigma_p(\tilde{T}^*)$ if and only if $L^{\bot_s}\cap \mathrm{M}(\bar{\lambda})\neq 0$ while the latter is equivalent to $\overline{L+\mathrm{M}(\lambda)}\neq \mathcal{B}_T$.

iv) The proof here is an abstraction of that of Thm.~1.3.4 in \cite{frey2005non}. We claim that $\textup{Ran}(\tilde{T}-\lambda)$ is closed if and only if $L+\mathrm{M}(\lambda)$ is closed in $\mathcal{B}_T$. Indeed we have a natural map
\[\tau:\mathcal{B}_T/\mathrm{M}(\lambda)\rightarrow \textup{Ran}(T^*-\lambda)/\textup{Ran}(T-\lambda)=H/\textup{Ran}(T-\lambda)\] given by $$\tau([x]+\mathrm{M}(\lambda))=(T^*x-\lambda x)+\textup{Ran}(T-\lambda).$$ It's easy to see $\tau$ is a topological linear isomorphism, sending $(L+\mathrm{M}(\lambda)/\mathrm{M}(\lambda)$ to $\textup{Ran}(\tilde{T}-\lambda)/\textup{Ran}(T-\lambda)$. Note that $L+\mathrm{M}(\lambda)$ is closed in $\mathcal{B}_T$ if and only if $(L+\mathrm{M}(\lambda))/\mathrm{M}(\lambda)$ is closed in $\mathcal{B}_T/\mathrm{M}(\lambda)$ and similarly that $\textup{Ran}(\tilde{T}-\lambda)/\textup{Ran}(T-\lambda)$ is closed in $H/\textup{Ran}(T-\lambda)$ if and only if $\textup{Ran}(\tilde{T}-\lambda)$ is closed in $H$. The claim then follows.

If $\textup{Ran}(\tilde{T}-\lambda)$ is closed, then $\dim\textup{ker}(\tilde{T}-\lambda)=\dim(L\cap \mathrm{M}(\lambda))$ by i). Additionally,
\[\dim\textup{coker}(\tilde{T}-\lambda)=\dim(H/\textup{Ran}(\tilde{T}-\lambda)).\]
Note that
\[H/\textup{Ran}(\tilde{T}-\lambda)=[H/\textup{Ran}(T-\lambda)]/[\textup{Ran}(\tilde{T}-\lambda)/\textup{Ran}(T-\lambda)].\]
We see that
\begin{eqnarray*}\dim(H/\textup{Ran}(\tilde{T}-\lambda))&=&\dim([\mathcal{B}_T/\mathrm{M}(\lambda)]/[(L+\mathrm{M}(\lambda))/\mathrm{M}(\lambda)])\\
&=&\dim(\mathcal{B}_T/(L+\mathrm{M}(\lambda))).\end{eqnarray*}
This shows that $\tilde{T}-\lambda$ is Fredholm if and only if $(L,\mathrm{M}(\lambda))$ is a Fredholm pair and in particular, $\textup{ind}\tilde{T}=\textup{ind}(L,\mathrm{M}(\lambda))$.
\end{proof}
\emph{Remark}. It's clear from the theorem that whether $\lambda\in \sigma(\tilde{T})$ or not is actually an intersection property of $L$ and $\mathrm{M}(\lambda)$ in $Gr(n,2n)$. This geometric flavour will be enhanced in \S~\ref{entire}.

If a trivialization $\Phi$ has been chosen, we are interested in self-adjoint extensions $T_U$ determined by the boundary condition $\Gamma_-x=U\Gamma_+x$ where $U\in \mathbb{U}(G)$ and $x\in D(T^*)$. Spectral analysis of $T_U$ can be carried out in terms of the boundary data $U$ and $B(\lambda)$. In the literature on boundary value problems of elliptic differential operators, this procedure is called "reduction to the boundary".
\begin{corollary}\label{coro1}If $T_U$ is the extension of $T$ defined as above and $\lambda\in \Theta(T)\cap \mathbb{R}$, then\\
i) $\lambda\in \rho(T_U)$ if and only if $U-B(\lambda)$ is invertible.\\
ii) $\lambda\in \sigma_p(T_U)$ if and only if $\textup{ker}(U-B(\lambda))\neq 0$. In particular, $\dim \textup{ker}(T_U-\lambda)=\dim \textup{ker}(U-B(\lambda))$.\\
iii) $\lambda\in \sigma_c(T_U)$ if and only if $\textup{ker}(U-B(\lambda))=0$ and $U-B(\lambda)$ is not invertible.\\
iv) $T_U-\lambda$ is Fredholm if and only if $U-B(\lambda)$ is Fredholm.
\end{corollary}
\begin{proof}For brevity, $\Phi(\mathrm{M}(\lambda))\subset G\oplus_\bot G$ shall still be denoted by $\mathrm{M}(\lambda)$.

i) Due to Lemma \ref{lem1}, $\lambda\in \rho(T_U)$ if and only if $L_U$ and $\mathrm{M}(\lambda)$ are transversal, where $L_U$ is the Lagrangian subspace
$\{(\varphi, U\varphi)\in G\oplus_\bot G|\varphi\in G\}$. Since $\mathrm{M}(\lambda)$ is of the same form with $U$ replaced by $B(\lambda)$. Thus $\lambda\in \rho(T_U)$ is equivalent to the invertibility of the operator matrix $\left(
                   \begin{array}{cc}
                     Id & U \\
                     Id & B(\lambda) \\
                   \end{array}
                 \right)
$, which is further equivalent to the invertibility of $U-B(\lambda)$.

ii) $\varphi\in \textup{ker}(U-B(\lambda))$ if and only if $(\varphi, U\varphi)\in L_U\cap \mathrm{M}(\lambda)$. The result then follows due to Thm.~\ref{thm4}, i).

iii) This follows from i) and ii).

iv) Let $P_U$ and $P_{B(\lambda)}$ be the orthogonal projections onto $L_U$ and $\mathrm{M}(\lambda)$ in $G\oplus_\bot G$ respectively. W.r.t. the decomposition $G\oplus_\bot G$,
\[P_U=\frac{1}{2}\left(
        \begin{array}{cc}
          Id & U^* \\
          U & Id \\
        \end{array}
      \right),\quad P_{B(\lambda)}=\frac{1}{2}\left(
        \begin{array}{cc}
          Id & B(\lambda)^* \\
          B(\lambda) & Id \\
        \end{array}
      \right).
\] From the general theory of Fredholm pairs (see for example the Appendix of \cite{frey2005non}), $(L_U, \mathrm{M}(\lambda))$ is a Fredholm pair if and only if $Id-P_U:\mathrm{M}(\lambda)\rightarrow \textup{ker}P_U$ is Fredholm. The conclusion then follows.
\end{proof}
\begin{corollary}$\gamma_+(\lambda)$ (resp. $\gamma_-(\lambda)$) has an analytic continuation on $\rho(T_+)$ (resp. $\rho(T_-)$). In particular, $B(\lambda)$ is analytically defined on $\rho(T_+)$.
\end{corollary}
\begin{proof}We only prove the claim for $\gamma_+(\lambda)$. By our previous arguments, we only need to prove $\gamma_+(\lambda)$ is well-defined and analytic on $\rho(T_+)\cap \mathbb{C}_-$. Note that for $\lambda\in \rho(T_+)\cap \mathbb{C}_-$ (of course $\lambda\in \Theta(T)$ as well), by Lemma \ref{lem1}, the operator  matrix $\left(
                   \begin{array}{cc}
                     0 & Id \\
                     B(\bar{\lambda})^* & Id \\
                   \end{array}
                 \right)
$ is invertible. Thus $B(\bar{\lambda})^*$ is invertible. This implies that Prop.~\ref{p2} and Corol.~\ref{c1} continue to hold on $\rho(T_+)\cap \mathbb{C}_-$. The result then follows.
\end{proof}

We now establish a resolvent formula of Krein-Naimark type in terms of $B(\lambda)$ for later use. Generally, such a formula is written in terms of Weyl function $M(\lambda)$ and linear relations. The advantage of our version is its convenience for spectral theory of self-adjoint extensions of $T$. Recall that the map $\gamma_+(\lambda)$ is from $G$ to $\textup{ker}(T^*-\lambda)$. However, we will view it as a map from $G$ to $H$ and the conjugate $\gamma_+(\lambda)^*$ is used in this sense.
\begin{lemma}\label{lem5}If $\lambda\in \Theta(T)\cap \mathbb{R}$, then
\[\gamma_+(\lambda)^*=-iB(\lambda)^{-1}\Gamma_-(T_+-\lambda)^{-1}.\]
\end{lemma}
\begin{proof}For any $x\in H$, let $y=(T_+-\lambda)^{-1}x\in D(T_+)$. Then for any $\varphi\in G$, we have
\begin{eqnarray*}(\gamma_+(\lambda)^*(T_+-\lambda)y,\varphi)_G&=&((T_+-\lambda)y, \gamma_+(\lambda)\varphi)_H\\&=&(T_+y,\gamma_+(\lambda)\varphi)_H-(y,\lambda\gamma_+(\lambda)\varphi)_H\\
&=&(T^*y,\gamma_+(\lambda)\varphi)_H-(y, T^*\gamma_+(\lambda)\varphi)_H\\
&=&i(\Gamma_+y, \Gamma_+\gamma_+(\lambda)\varphi)_G-i(\Gamma_-y, \Gamma_-\gamma_+(\lambda)\varphi)_G\\
&=&-i(\Gamma_-y, B(\lambda)\varphi)_G=-i(B(\lambda)^{-1}\Gamma_-y, \varphi)_G.
\end{eqnarray*}
The formula then follows, where we have used the facts $\Gamma_+y=0$ and $B(\lambda)^*=B(\lambda)^{-1}$ for $\lambda\in \Theta(T)\cap \mathbb{R}$.
\end{proof}
\emph{Remark}. Note that on the orthogonal complement of $\textup{ker}(T^*-\lambda)$, $\gamma_+(\lambda)^*$ is zero. Along the same line, one can prove that for any $\lambda\in \rho(T_+)$,
\begin{equation}\label{e1}\gamma_-(\bar{\lambda})^*=-i\Gamma_-(T_+-\lambda)^{-1}.
\end{equation}

If $I\subset \Theta(T)\cap \mathbb{R}$ is an open interval, then $B(\lambda)$ over $I$ is a real analytic curve in $\mathbb{U}(G)$. Consequently, $B(\lambda)^{-1}\frac{dB(\lambda)}{d\lambda}$ takes its values in the Lie algebra of $\mathbb{U}(G)$, i.e., the Lie algebra formed by skew-Hermitian operators.
\begin{proposition}\label{p4}Let $I$ be as above. Then for $\lambda\in I$
\[-iB(\lambda)^{-1}\frac{dB(\lambda)}{d\lambda}=\gamma_+(\lambda)^*\gamma_+(\lambda).\]
For $\lambda\in \rho(T_+)$, we have
\[\frac{dB(\lambda)}{d\lambda}=i\gamma_-(\bar{\lambda})^*\gamma_+(\lambda).\]
\end{proposition}
\begin{proof}Due to Corol.~\ref{c1}, for $\lambda\in I$,
\[\frac{dB(\lambda)}{d\lambda}=\Gamma_-(T_+-\lambda)^{-1}\gamma_+(\lambda).\]
The first formula then follows from Lemma \ref{lem5} and the second from Eq.~(\ref{e1}).
\end{proof}
\emph{Remark}. For real $\lambda\in \Theta(T)$, since $\gamma_+(\lambda)$ is an isomorphism between $G$ and $\textup{ker}(T^*-\lambda)$, $\gamma_+(\lambda)^*\gamma_+(\lambda)$ is positive-definite and invertible on $G$. If the deficiency indices of $T$ are $(1,1)$, then $B(\lambda)$ is of the form $e^{i\theta(\lambda)}$ for a real analytic "phase" function $\theta(\lambda)$ on $I$. The above formula in this case simply implies $\theta'>0$. Therefore, as $\lambda$ increases in $I$, so does the angle $\theta$.

Note that by the above proposition for $\lambda_0\in \Theta(T)\cap \mathbb{R}$, $B'(\lambda_0)$ is invertible. This fact implies the following result revealing the particularity of $\lambda_0$ in terms of the curvature of $T$.
\begin{proposition}If $\lambda_0\in \Theta(T)\cap\mathbb{R}$, then $\lim_{\widehat{\lambda\rightarrow \lambda_0}}r_T(\lambda)=0$. Here $\widehat{\lambda\rightarrow \lambda_0}$ means $\lambda$ approaches $\lambda_0$ in $\mathbb{C}_+$.
\end{proposition}
\begin{proof}We assume a boundary triplet has been fixed and the contractive Weyl function is analytic and invertible around $\lambda_0\in \mathbb{C}$. If
$B(\lambda)=\sum_{i=0}^\infty B_i(\lambda-\lambda_0)^i$
is the Taylor expansion of $B(\lambda)$ around $\lambda_0$, then $B_0=B(\lambda_0)$ and $B_1=B'(\lambda_0)$. For $\lambda\in \mathbb{C}_+$, the formula $B(\lambda)=(B(\bar{\lambda})^*)^{-1}$ implies
$B_0B_1^*+B_1B_0^*=0$,
and
\[Id-B^*B=2i\Im \lambda (B_1^*B_0+\varphi_1(\lambda)),\quad Id-BB^*=-2i\Im \lambda(B_1B_0^*+\varphi_2(\lambda))\]
where both $\varphi_1(\lambda)$ and $\varphi_2(\lambda)$ approach zero as $\lambda$ goes to $\lambda_0$ in $\mathbb{C}_+$.

Recall from Thm.~\ref{thm5}, in terms of the fixed boundary triplet,
\[r_T=Id-4(\Im{\lambda})^2\mathrm{K}^{-1}B'^*\tilde{\mathrm{K}}^{-1}B'\]
where $\mathrm{K}=Id-B^*B$ and $\tilde{\mathrm{K}}=Id-BB^*$. Thus
\[\lim_{\widehat{\lambda\rightarrow \lambda_0}}4(\Im{\lambda})^2\mathrm{K}^{-1}B'^*\tilde{\mathrm{K}}^{-1}B'=(iB_1^*B_0)^{-1}B_1^*(-iB_1B_0^*)^{-1}B_1=Id.
\]
The conclusion then follows.
\end{proof}

Now we come to the spectral analysis of self-adjoint extensions of $T$. For the fixed symplectic isomorphism $\Phi$, these extensions are parameterized by $U\in \mathbb{U}(G)$, which correspond to boundary conditions $\Gamma_-x=U\Gamma_+x$. We denote the extension corresponding to $U$ by $T_U$. Recall that $T_{Id}=T_0$ and and $T_{-Id}=T_1$.
\begin{lemma}\label{lem2} For $\lambda\in \Theta(T)\cap \mathbb{R}$, the domain $D(T_U)$ of $T_U$ can be characterized in the following way:
\begin{eqnarray*}D(T_U)&=&\{x=(T_+-\lambda)^{-1}(y+w)+\gamma_+(\lambda)\varphi|\varphi\in G, w\in\textup{ker}(T^*-\lambda), \\&y&\in \textup{Ran}(T-\lambda), iB(\lambda)\gamma_+(\lambda)^*w=(U-B(\lambda))\varphi\}\end{eqnarray*}
\end{lemma}
\begin{proof}Denote the right hand side by $S$. If $x\in D(T_U)$, then $x=x_\lambda+z_\lambda$ where $x_\lambda\in D(T_+)$ and $z_\lambda\in \textup{ker}(T^*-\lambda)$. This is due to the decomposition $D(T^*)=D(T_+)\oplus\textup{ker}(T^*-\lambda)$ (see the proof of Lemma \ref{lem1}). Since $H=\textup{Ran}(T-\lambda)\oplus_\bot \textup{ker}(T^*-\lambda)$, there are $y\in \textup{Ran}(T-\lambda)$ and $w\in \textup{ker}(T^*-\lambda)$ such that $(T_+-\lambda)x_\lambda=y+w$. Set $y=(T-\lambda)v$ for $v\in D(T)$. Then we have $(T_+-\lambda)(x_\lambda-v)=w$, implying
$x_\lambda=v+(T_+-\lambda)^{-1}w$. Thus
\[\Gamma_-x_\lambda=\Gamma_-(T_+-\lambda)^{-1}w=iB(\lambda)\gamma_+(\lambda)^*w.\]
Besides, we also have $\Gamma_-z_\lambda=B(\lambda)\Gamma_+z_\lambda$. Since $\Gamma_-x=U\Gamma_+x$, we obtain
$$\Gamma_-x_\lambda+\Gamma_-z_\lambda=U\Gamma_+z_\lambda$$
 where the fact $\Gamma_+x_\lambda=0$ was used. Combining these facts together and setting $\varphi=\Gamma_+z_\lambda$, we finally have $D(T_U)\subset S$. The inclusion $S\subset D(T_U)$ can be checked directly.
\end{proof}
Note that if $x=(T_+-\lambda)^{-1}(y+w)+\gamma_+(\lambda)\varphi\in D(T_U)$, then
\[(T_U-\lambda)x=(T^*-\lambda)x=y+w.\]
\begin{theorem}\label{reso}(Resolvent formula) If $\lambda\in \rho(T_U)\cap \mathbb{R}$ for $U\in \mathbb{U}(G)$, then $U-B(\lambda)$ is invertible and
\[(T_U-\lambda)^{-1}-(T_+-\lambda)^{-1}=i\gamma_+(\lambda)(B(\lambda)^{-1}U-Id)^{-1}\gamma_+(\lambda)^*.\]
\end{theorem}
\begin{proof}Certainly $\lambda$ lies in $\Theta(T)$. If $z\in H$, set $x=(T_U-\lambda)^{-1}z\in D(T_U)$. Then from the characterization of $D(T_U)$ in Lemma \ref{lem2},
\begin{eqnarray*}(T_U-\lambda)^{-1}z&=&x=(T_+-\lambda)^{-1}(T_U-\lambda)x+\gamma_+(\lambda)\Gamma_+x\\
&=&(T_+-\lambda)^{-1}z+\gamma_+(\lambda)\Gamma_+x.\end{eqnarray*}
Again due to Lemma \ref{lem2},
\[\Gamma_+x=i(U-B(\lambda))^{-1}B(\lambda)\gamma_+(\lambda)^*w=i(U-B(\lambda))^{-1}B(\lambda)\gamma_+(\lambda)^*(y+w)\]
for $\gamma_+(\lambda)^*y=0$. Since $y+w=(T_U-\lambda)x=z$, the resolvent formula then follows.
\end{proof}
\emph{Remark}. In the traditional Krein formula (see for example \cite[Thm.~2.6.1]{behrndt2020boundary}), $T_0$, rather than our $T_+$, plays the role of a background operator. Furthermore, in that formula, $\lambda$ should lie in $\rho(T_U)\cap\rho(T_0)$ and the right hand side of the formula often involves the inverse of a linear relation. This inconvenience can be avoided in our formula.
\begin{corollary}If $\lambda\in \rho(T_{U_1})\cap \rho(T_{U_2})\cap\mathbb{R}$ where $U_1, U_2\in \mathbb{U}(G)$, then
\[(T_{U_1}-\lambda)^{-1}-(T_{U_2}-\lambda)^{-1}=i\gamma_+(\lambda)[(B(\lambda)^{-1}U_1-Id)^{-1}-(B(\lambda)^{-1}U_2-Id)^{-1}]\gamma_+(\lambda)^*.\]
\end{corollary}
\begin{proof}This is clear from the resolvent formula.
\end{proof}

The above Krein-type resolvent formula can be extended slightly: if $\lambda\in \rho(T_+)$ and a boundary triplet has been chosen such that a generic boundary condition is parameterized by $Y\in \mathbb{B}(G)$, i.e., \[D(T_Y)=\{x\in D(T^*)|\Gamma_-x=Y\Gamma_+x\},\] then in terms of $T_+$,
\begin{eqnarray*}D(T_Y)&=&\{x=(T_+-\lambda)^{-1}(y+w)+\gamma_+(\lambda)\varphi|\varphi\in G, w\in\textup{ker}(T^*-\bar{\lambda}), \\&y&\in \textup{Ran}(T-\lambda), i\gamma_-(\bar{\lambda})^*w=(Y-B(\lambda))\varphi\}\nonumber.\end{eqnarray*}
This can be proved along the same line of Lemma \ref{lem2}. Using this statement and Eq.~(\ref{e1}), if $\lambda\in \rho(T_+)\cap \rho(T_Y)$, then we obtain
\begin{equation}(T_Y-\lambda)^{-1}-(T_+-\lambda)^{-1}=i\gamma_+(\lambda)(Y-B(\lambda))^{-1}\gamma_-(\bar{\lambda})^*.\label{e2}\end{equation}
If $\lambda\in \rho(T_-)\cap \rho(T_Y)$, by replacing $T_+$ with $T_-$, a similar formula holds.

\section{Analytic singularity of symmetric operators}\label{asing}
If $T$ is a simple symmetric operator with deficiency indices $(n,n)$, we denote the set of all self-adjoint extensions of $T$ by $\mathcal{S}$. There are conceptually three closed subsets of $\mathbb{R}$ associated with $T$. The first is $\cap_{\tilde{T}\in \mathcal{S}}\sigma(\tilde{T})$. The other two are as follows.
\begin{definition}(\cite{akhiezer2013theory}) The set $\sigma_k(T):=\mathbb{C}\backslash\Theta(T)$ is called the spectral kernel of $T$.
\end{definition}
\begin{definition}$\lambda_0\in \mathbb{R}$ is called an analytic singular point of $T$ if the two-branched Weyl map $W_T(\lambda)$ doesn't admit an analytic continuation around $\lambda_0$ in $\mathbb{C}$. The set $\sigma_a(T)\subset \mathbb{R}$ of all analytic singular points of $T$ is called the analytic spectrum of $T$. In particular, if $\sigma_a(T)=\emptyset$, we call $T$ an entire operator.
\end{definition}
It should be emphasized that though with a boundary triplet the Weyl map may be represented by the Weyl function $M(\lambda)$, $M(\lambda)$ may not admit an analytic continuation to where $W_T(\lambda)$ admits such a continuation. The singularities of $M(\lambda)$ is a mixture of contributions from both $W_T(\lambda)$ and the boundary condition $U=Id$. After all, $M(\lambda)$ is a basic tool to study the spectrum of $T_0$.
\begin{proposition}\label{p14}For a simple symmetric operator $T$ with deficiency indices $(n,n)$, $\sigma_k(T)=\cap_{\tilde{T}\in \mathcal{S}}\sigma(\tilde{T})=\sigma_a(T)$.
\end{proposition}
\begin{proof}$\sigma_k(T)=\cap_{\tilde{T}\in \mathcal{S}}\sigma(\tilde{T})$ seems to be a classical result, but we cannot find a suitable reference. For the convenience of the reader, we give an argument in our formalism. Obviously, $\cap_{\tilde{T}\in \mathcal{S}}\sigma(\tilde{T})\supset\sigma_k(T)$. Let $B(\lambda)$ be the contractive Weyl function of $T$ w.r.t. a chosen boundary triplet $(G, \Gamma_\pm)$. Now, if $\lambda_0\in \Theta(T)\cap \mathbb{R}$, then by our previous discussion $B(\lambda)$ is analytic around $\lambda_0$ in $\mathbb{C}$. We can choose $U\in \mathbb{U}(G)$ such that $U-B(\lambda_0)$ is invertible. Due to Coro.~\ref{coro1}, $\lambda_0\in \rho(T_U)$. Thus $\lambda_0\not\in \cap_{\tilde{T}\in \mathcal{S}}\sigma(\tilde{T})$.

By our results in the previous section, the two-branched Weyl map $W_T(\lambda)$ admits an analytic continuation to $\Theta(T)$. Thus, $\sigma_a(T)\subset \sigma_k(T)$.

 Conversely, if the real number $\lambda_0\in \mathbb{R}\backslash \sigma_a(T)$, then with a certain boundary triplet $(G, \Gamma_\pm)$ the contractive Weyl function $B(\lambda)$ is holomorphic around $\lambda_0$ in $\mathbb{C}$ and due to Prop.~\ref{p11} $B(\lambda)$ is unitary if $\lambda\in \mathbb{R}$. We can choose a self-adjoint extension parameterized by a unitary operator $U\in \mathbb{U}(G)$ w.r.t. the chosen boundary triplet such that $U-B(\lambda_0)$ is invertible. Since the pseudounitary group $\mathbb{U}(n, n)$ acts transitively on the space $\mathcal{L}(G\oplus_\bot G)$ of Lagrangians in $G\oplus_\bot G$, we can change the boundary triplet such that $U$ is transformed into $Id$. Certainly $B(\lambda)$ also changes accordingly, but  $Id-B(\lambda_0)$ is still invertible because invertibility of $U-B(\lambda_0)$ simply means $L_U$ and $\mathrm{M}(\lambda_0)$ in Thm.~\ref{thm4} are transversal, and transversality is unchanged under the transformation. By continuity, $Id-B(\lambda)$ is invertible for $\lambda$ near $\lambda_0$. Note that, with the new boundary triplet, for $\lambda\in \mathbb{C}_+$
\[M(\lambda)=i(Id+B(\lambda))(Id-B(\lambda))^{-1}.\]
This expression shows that $M(\lambda)$ admits an analytic continuation around $\lambda_0$ in $\mathbb{C}$. It is a standard fact that this means that $\lambda_0\in \rho(T_0)$ (see, for example,  \cite[Thm.~3.6.1]{behrndt2020boundary}) and thus $\lambda_0\in \Theta(T)$. Therefore $\sigma_k(T)\subset \sigma_a(T)$.
\end{proof}
\emph{Remark}. The result $\sigma_k(T)=\sigma_a(T)$ is not trivial at all. To our best knowledge, this has only been established before for the case $n=1$ \cite[pp.420, Thm.~2]{akhiezer2013theory}. In this special case, it's natural to view $M(\mathbb{C}_+\cup \mathbb{C}_-)$ as a subset of $\mathbb{C}$ and consequently we know very well where to continue (if possible) the Weyl function $M(\lambda)$ for $\lambda\in \mathbb{R}$. Actually, $M(\lambda)$ is a meromorphic function on $\Theta(T)$. In the general case, it is the fact $W_T(\mathbb{C}_+\cup \mathbb{C}_-)\subset Gr(n,2n)$ that tells us where to continue the Weyl curve naturally.
\begin{proposition}Let $T$ have deficiency indices $(n,n)$ with $n<+\infty$. If $M(\lambda)$ is the Weyl function of $T$ w.r.t. a chosen boundary triplet, then $M(\lambda)$ is meromorphic on $\Theta(T)$, and the corresponding contractive Weyl function $B(\lambda)$ can be continued to be meromorphic on $\Theta(T)$ (of course analytic at points in $\Theta(T)\cap \mathbb{R}$) as well.
\end{proposition}
\begin{proof}
The first claim is clear. Note that for $\lambda\in \mathbb{C}_+$, $B(\lambda)=(M(\lambda)-i)(M(\lambda)+i)^{-1}$. Since $M(\lambda)+i$ is invertible on $\mathbb{C}_+$, the holomorphic Fredholm theorem implies that $(M(\lambda)+i)^{-1}$ is analytic on $\Theta(T)\backslash(\mathcal{P}(M)\cup S)$ where $\mathcal{P}(M)$ is the set of poles of $M(\lambda)$ on $\Theta(T)\cap \mathbb{R}$ and $S$ is a discrete subset of $\Theta(T)\backslash\mathcal{P}(M)$. Around a point $\lambda_0\in S$, $B(\lambda)=Q(\lambda)/\det(M(\lambda)+i)$ for an analytic matrix-valued function $Q(\lambda)$. Therefore $B(\lambda)$ is meromorphic at $\lambda_0$. Due to results in our previous section, $B(\lambda)$ is analytic at points in $\mathcal{P}(M)$. The conclusion then follows.
\end{proof}
\emph{Remark}. It is not hard to find that the poles of $B(\lambda)$ (zeros of $\det(M(\lambda)+i$) on $\mathbb{C}_-$ are precisely eigenvalues of $T_+$ on $\mathbb{C}_-$ while zeros of $\det B(\lambda)$ on $\mathbb{C}_+$ are precisely eigenvalues of $T_-$ on $\mathbb{C}_+$.

\begin{corollary}The set of entire operators with deficiency indices $(n,n)$ in $H$ coincides with the set of regular operators with deficiency indices $(n,n)$ in $H$.
\end{corollary}
\begin{proof}This is obvious from the above Prop.~\ref{p14}.
\end{proof}
Consequently by definition the Weyl curve of an entire operator is actually an entire curve in $Gr(n, 2n)$.

\emph{Remark}. There is another notion of entire operators introduced by M. G. Krein \cite{gorbachuk2012mg}. Such an operator is regular, and additionally a specific choice of a "gauge" or "generalized gauge" is also involved in this notion. Thus our definition is more general. Besides, though regular operators are precisely those that are entire, we still suggest the name "entire operator" for the obvious reason.

If $n\in \mathbb{N}$, $Gr(n,2n)$ is a projective manifold. In a projective manifold $M$, an entire curve $f:\mathbb{C}\rightarrow M$ is called algebraically non-degenerate if $f(\mathbb{C})$ lies in no proper subvariety of $M$, i.e., the Zariski closure of $f(\mathbb{C})$ is the whole $M$. This motivates the following definition.
\begin{definition}
An entire operator $T$ with deficiency indices $(n,n)$ where $n\in \mathbb{N}$ is called algebraically non-degenerate if its Weyl curve is algebraically non-degenerate in $Gr(n,2n)$.
\end{definition}
This definition and its variant will prove important in the following section.
\section{Entire operators with finite deficiency index}\label{entire}
If $T$ is an entire operator with deficiency indices $(n,n)$, for brevity we shall just say $T$ has deficiency index $n$. Throughout this section, we assume $T$ is an entire operator with \emph{finite} deficiency index $n$.
\subsection{Entire operators and value distribution theory}
 Let $T$ be an entire operator with deficiency index $n$. Then its Weyl curve $W_T(\lambda)$ is an entire curve in $Gr(n,2n)$. Entire curves in algebraic varieties have been studied in value distribution theory  for several decades, although it's hard to say the theory has reached its final stage. A Second Main Theorem on entire curves in projective algebraic manifolds has only been obtained in 2009 by Ru in \cite{ru2009holomorphic}, which was a major breakthrough in this field. Our geometric viewpoint towards Weyl functions then naturally leads us to introduce value distribution theory into our picture of simple symmetric operators.

Let us recall briefly the basics of value distribution theory. For a detailed account, we refer the reader to the recently published book \cite{ru2021nevanlinna}. In its original form developed by R. Nevanlinna in 1920s, this theory was concerned with the value distribution of meromorphic functions on $\mathbb{C}$. Since then, the theory has been extended in several different directions, of which the most relevant for us is the following geometric generalization.

Note that a meromorphic function on $\mathbb{C}$ can be viewed as a holomorphic map from $\mathbb{C}$ to the projective line $\mathbb{CP}^1$ or the Riemann sphere in complex analysis. We can replace $\mathbb{CP}^1$ with a general projective manifold $M$ and consider an entire curve $f: \mathbb{C}\rightarrow M$ as a generalization of meromorphic functions. Points in $\mathbb{CP}^1$ can be interpreted as divisors, and their counterparts in the generalized setting are (Cartier) divisors in $M$.

Recall that a Cartier divisor $D$ on $M$ can be described as a collection $\{U_i, f_i\}$, where $\{U_i\}$ is an open cover of $M$ and $f_i$ a meromorphic function on $U_i$ such that on each overlap $U_i\cap U_j$, $f_i/f_j$ is holomorphic and non-vanishing pointwise. If each $f_i$ can be chosen to be holomorphic, the divisor is called \emph{effective}. Associated with a Cartier divisor is a holomorphic line bundle $\mathcal{O}(D)$ whose transition functions are given by those $f_i/f_j$. $D$ is called \emph{ample}, if $\mathcal{O}(D)$ has a positive metric $h$, i.e., the first Chern form $c_1(\mathcal{O}(D),h)$ is a K$\ddot{a}$hler metric on $M$. \footnote{For line bundles, this notion of positivity coincides with Griffiths' notion of positivity.} Roughly speaking, an effective divisor $D$ in $M$ is a collection of hypersurfaces (possibly with singularities and multiplicities) in $M$, locally defined by the equations $f_i=0$. If $\mathcal{O}(D)$ is topologically trivial, $D$ is called principal ($D$ actually comes from a meromorphic function on $M$). If two divisors $D_1=\{U_i, f_i\}$, $D_2=\{U_i, g_i\}$ are such that $\{U_i, f_i/g_i\}$ is a principal divisor, then $D_1$ and $D_2$ are called linearly equivalent and this relation is denoted by $D_1\sim D_2$.

 Given an entire curve $f:\mathbb{C}\rightarrow M$ in a projective manifold $M$. Assume that $L$ is a holomorphic line bundle over $M$ equipped with a Hermitian metric $h$ and $s$ a holomorphic section of $L$. We denote the zero locus of $s$ by $(s)$, which is an effective Cartier divisor. The following three functions defined in terms of these data are basic objects of value distribution theory.

 \begin{definition} The characteristic function of $f$ w.r.t. $(L, h)$ is defined as
\[T_{f,L}(r)=\int_0^r\frac{dt}{t}\int_{B_t}f^*(c_1(L,h)),\]
where $c_1(L,h)$ is the first Chern form of $(L, h)$, $B_t$ the open disc centered at $0\in \mathbb{C}$ and with radius $t>0$, and $f^*$ the pull-back of differential forms.
\end{definition}
It can be proved that if $L$ is ample and $T_{f,L}(r)$ is bounded, then $f$ must be constant. If $T_{f,L}(r)=O(\ln r)$ as $r$ goes to $+\infty$, then $f$ is rational in the sense that in each affine coordinate chart, $f$ is rational in $\lambda\in \mathbb{C}$, see \cite[Thm.~2.5.28]{noguchi2013nevanlinna}.

\begin{definition}If $s\circ f\not\equiv 0$, the proximity function of $f$ w.r.t. $D=(s)$ is defined to be
\[m_f(r, D)=-\int_0^{2\pi}\ln\|s(f(re^{i\theta}))\|\frac{d\theta}{2\pi}.\]
\end{definition}
$m_f(r,D)$ measures how close $f$ is, on average, to $D$ on the circle $\partial B_r$. It is important to note that both $T_{f, L}(r)$ and $m_f(r, D)$ are essentially independent of the choice of the Hermitian metric $h$, because a different $h$ only modifies these functions with a bounded term in $r$. Besides, $m_f(r,D)$ is bounded from below. It should be pointed out that $m_f(r,D)$ is not completely determined by $h$ and $D$, because for a nonzero constant $c\in \mathbb{C}$ the section $c\times s$ gives the same divisor. However, this ambiguity only leaves a constant to be added to $m_f(r, D)$.
\begin{definition}
If $s\circ f\not\equiv 0$, the counting function of $f$ w.r.t. $D=(s)$ is defined to be
\[N_f(r, D)=\int_0^r [n_f(t, D)-n_f(0,D)]\frac{dt}{t}+n_f(0,D)\ln r,\]
where $n_f(t, D)$ is the number of roots of the equation $s\circ f=0$ in the disc $B_t$ and $n_f(0,D)=\lim_{r\rightarrow 0+}n_f(r,D)$. Note that roots are counted according to their analytic multiplicity.
\end{definition}
The above three functions are interrelated in the following manner.
\begin{theorem}(First Main Theorem) Let $f:\mathbb{C}\rightarrow M$ be an entire curve in the projective manifold $M$ and $L$ a Hermitian line bundle over $M$. If $s$ is a holomorphic section of $L$ with $D=(s)$ and $s\circ f\not\equiv 0$, then
\[T_{f,L}(r)=m_f(r, D)+N_f(r, D)+O(1),\]
where $O(1)$ is a bounded term in $r$.
\end{theorem}
The First Main Theorem is a version of the famous Poincar$\acute{e}$-Lelong formula
\[dd^c[\ln \|s\|^2]=-c_1(L,h)+[D],\]
which is an identity of currents. A careful examination shall reveal that the bounded term $O(1)$ is actually independent of $r$ but does depend on $s$.

We are now in a position to explain the relation of abstract boundary value problems with value distribution theory.

An entire operator $T$ of course provides its Weyl curve serving as the entire curve in the above argument and $M=Gr(n,2n)$. $Gr(n,2n)$ has a tautological holomorphic vector bundle $E$ of rank $n$, whose fiber at $x\in Gr(n,2n)$ is precisely the subspace $V_x\subset \mathcal{B}_T$ of dimension $n$ parameterized by $x$. Let $l$ be the determinant line bundle of $E$, i.e., $l=\det E=\wedge^nE$, and $L:=l^*$ is the dual of $l$. It can be easily checked that $L$ is ample. A generic boundary condition $y$ is a point in $Gr(n,2n)$, determining a holomorphic section $s_y$ of $L$ up to a constant factor. Indeed, if $\{e_1,\cdots, e_n\}$ is a basis of $V_y$, then we have a nonzero element $\sigma:=e_1\wedge\cdots\wedge
e_n\in\textup{det}V_y$. Let $\nu$ be a fixed nonzero element of $\textup{det}\mathcal{B}_T=\wedge^{2n}\mathcal{B}_T$. Then for any $\tau\in l_x$, we define $s_y(\tau)\in \mathbb{C}$ by $\sigma\wedge \tau=s_y(\tau)\nu$. The (well-defined) zero locus of $s_y$ then provides the Cartier divisor $D$--sometimes it is also called the Schubert hyperplane determined by $y$. We use $T_y$ to denote the closed extension of $T$ parameterized by $y\in Gr(n,2n)$. If $y_1,y_2\in Gr(n,2n)$ are such that $V_{y_1}$ and $V_{y_2}$ are transversal in $\mathcal{B}_T$, we also say $y_1$ and $y_2$ are transversal boundary conditions.

The relative position of $W_T$ w.r.t. $(s_y)$ for $y\in Gr(n,2n)$ is clearly of great importance. We say $y\in Gr(n,2n)$ is degenerate w.r.t. $T$ if $s_y\circ W_T\equiv 0$. The following lemma seems elementary, but we don't know where to find a suitable reference.
\begin{lemma}Let $D$ be a disc in $\mathbb{C}$. If $W_T(D)$ lies in an algebraic subvariety $V\subset Gr(n,2n)$ of codimension 1, then $W_T(\mathbb{C})\subset V$.
\end{lemma}
\begin{proof}Let $U\supset D$ be the maximal connected subset of $\mathbb{C}$ such that $W_T(U)\subset V$. Obviously, $U$ is a closed subset. If $\lambda_0$ is a boundary point of $U$, then there is a sequence $\{\lambda_j\}\subset U$ such that $\lim_{j\rightarrow \infty}\lambda_j=\lambda_0$. If $f=0$ is a defining equation of $V$ around $W_T(\lambda_0)$, then by definition $f(W_T(\lambda_j))=0$ for sufficiently large $j$. By the Identity Theorem, $f(W_T(\lambda))\equiv 0$ around $\lambda_0$. A contradiction! This shows that $U$ is also open. Consequently $U=\mathbb{C}$.
\end{proof}

\begin{proposition}\label{p9}For $y\in Gr(n,2n)$, $T_y$ only has eigenvalues (if $\sigma(T_y)\neq \varnothing$). $\lambda\in \mathbb{C}$ is an eigenvalue of $T_y$ if and only if $s_y(W_T(\lambda))=0$. In particular, if $y$ is non-degenerate w.r.t. $T$, then any eigenvalue of $T_y$ is isolated.
\end{proposition}
\begin{proof}The first two statements are clear from the definition of $s_y$ and Thm.~\ref{thm4}. The third follows from the above lemma.
\end{proof}
If $y$ is a degenerate boundary condition, then $s_y(W_T(\lambda))\equiv 0$ and hence $\sigma(T_y)=\mathbb{C}$. This can really happen and we will give simple examples in \S~\ref{sturm}. In the following sense degenerate boundary conditions are scarce.
\begin{proposition}The locus $\mathcal{D}$ of degenerate boundary conditions in $Gr(n,2n)$ is an algebraic subvariety of complex codimension  $\leq 2$.
\end{proposition}
\begin{proof}Note that
\[\mathcal{D}=\{y\in Gr(n,2n)|s_y\wedge s_{W_T(\lambda)}=0,\quad \forall \lambda\in \mathbb{C}\}.\]
$\mathcal{D}$ is clearly a Zariski closed subset of $Gr(n,2n)$. By Chow's famous theorem in algebraic geometry, $\mathcal{D}$ is an algebraic subvariety of $Gr(n,2n)$. Since $W_T(\lambda)$ is never constant, there are $\lambda_1\neq \lambda_2$ such that $W_T(\lambda_1)\neq W_T(\lambda_2)$. Therefore, $\mathcal{D}$ has to be a proper subvariety of the Schubert hyperplane determined by $W_T(\lambda_1)$ and the claim follows.
\end{proof}

Thus for $n=1$ no degenerate boundary condition exists.

Motivated by value distribution theory, we can introduce several useful concepts for our entire operator $T$.
\begin{definition}An entire operator $T$ is said to be weakly algebraically non-degenerate if any $y\in Gr(n,2n)$ is non-degenerate w.r.t. $T$. \end{definition}
If $W_T(\lambda)$ is algebraically degenerate in $Gr(n,2n)$, we call the Zariski closure $\mathcal{Z}_T$ of $W_T(\mathbb{C})$ in $Gr(n,2n)$ the Weyl variety of $T$. Obviously, $\mathcal{Z}_T$ is the natural ambient space suitable for investigating the value distribution theory of $W_T(\lambda)$. It's interesting to know how the geometry of $\mathcal{Z}_T$ reflects properties of $T$.
\begin{definition}
If $\lambda_0\in \mathbb{C}$ is an isolated eigenvalue of $T_y$ for $y\in Gr(n, 2n)$, the analytic multiplicity of $\lambda_0$ is defined to be the order of $\lambda_0$ as the zero of $s_y(W_T(\lambda))=0$.
\end{definition}
As far as we know, this notion was only defined before for some ordinary differential operators.
\begin{definition}\label{d1}
The height $h_T(r)$ of an entire operator $T$ is the characteristic function of its Weyl curve $W_T(\lambda)$ w.r.t. the line bundle $L=l^*$. If
\[\liminf_{r\rightarrow +\infty}\frac{h_T(r)}{\ln r}=+\infty,\]
we say $T$ is transcendental. The Weyl order $\rho_T$ of $T$ is defined by
\[\rho_T=\limsup_{r\rightarrow+\infty}\frac{\ln h_T(r)}{\ln r}.\]
If $\rho_T$ is finite, we say $T$ is of finite Weyl order, and the number
\[\tau_T=\limsup_{r\rightarrow+\infty}\frac{h_T(r)}{r^{\rho_T}}\]
is called the Weyl type of $T$.
\end{definition}
It should be emphasized that $h_T(r)$ in essence only depends on the unitary equivalence class of $T$, and so do the Weyl order and type. If $T$ is not transcendental, we say $T$ is rational, for the Weyl curve $W_T(\lambda)$ has to be a rational curve in $Gr(n,2n)$ in this case.

\begin{definition}The proximity function $m_T(r, y)$ of $T$ w.r.t. the non-degenerate boundary condition $y\in Gr(n,2n)$ is the proximity function $m_{W_T}(r, D)$ where $D=(s_y)$.
\end{definition}
Since $\mathcal{B}_T$ has a natural inner product, there is a canonical way to choose a metric $h$ on $L$. However, when a boundary triplet $(G, \Gamma_\pm)$ has been fixed and $\mathcal{B}_T$ is  identified with $G\oplus_\bot G$, we prefer to equip $L$ with the natural metric induced from the inner product on $G\oplus_\bot G$. Then $c_1(L, h)$ is a K$\ddot{a}$hler metric $\omega$ on $Gr(n, 2n)$, which is clearly $\mathbb{U}(2n)$-invariant. When a boundary triplet is fixed, we always make such a choice.
\begin{definition}The counting function $N_T(r, y)$ of $T$ w.r.t. the non-degenerate boundary condition $y\in Gr(n,2n)$ is $N_{W_T}(r, D)$ where $D=(s_y)$.
\end{definition}
Due to Prop.~\ref{p9}, $n_T(r,y):=n_{W_T}(r, D)$ for $D=(s_y)$ is the number of eigenvalues of $T_y$ in $B_r$, counting analytic multiplicities. Thus $N_T(r, y)$ measures the distribution of eigenvalues of $T_y$ in $B_r$. In this setting, the First Main Theorem demonstrates the general pattern for abstract boundary value problems of $T$: while $y$ can vary in $Gr(n, 2n)$ and consequently both $m_T(r, y)$ and $N_T(r, y)$ vary accordingly, the summation of the two essentially doesn't depend on the chosen boundary condition $y$ and is only controlled by $T$ itself. Without value distribution theory, this observation is far from obvious.

\begin{definition}Let $y\in Gr(n,2n)$ be non-degenerate w.r.t. the transcendental entire operator $T$ and $\{\lambda_i\}$ the nonzero eigenvalues of $T_y$, counting multiplicity and ordered such that $|\lambda_i|\leq |\lambda_{i+1}|$. Then the infimum $\alpha_T(y)$ of positive numbers $\alpha$ such that
\[\sum_i\frac{1}{|\lambda_i|^\alpha}< +\infty\]
is called the Weyl exponent of $y$ w.r.t. $T$.
\end{definition}
\begin{proposition}For a non-degenerate boundary condition $y\in Gr(n,2n)$ w.r.t. the transcendental entire operator $T$, $\alpha_T(y)\leq \rho_T$.
\end{proposition}
\begin{proof}The proof follows from the discussion in \cite[\S~1, Chap.~2]{gol'dberg2008value}. Note that though the book only considers value distribution of meromorphic functions, the argument there only depends on a general study on growth category and the First Main Theorem.
\end{proof}
\begin{definition}Let $T$ be a transcendental entire operator with deficiency index $n$. The Nevanlinna defect $\delta_T(y)$ of a non-degenerate boundary condition $y\in Gr(n,2n)$ w.r.t. $T$ is defined by
\[\delta_T(y)=\liminf_{r\rightarrow +\infty}\frac{m_T(r, y)}{h_T(r)}.\]
The Valiron defect $\Delta_T(y)$ of $y$ w.r.t. $T$ is defined by
\[\Delta_T(y)=\limsup_{r\rightarrow +\infty}\frac{m_T(r, y)}{h_T(r)}.\]
\end{definition}
Note that $0\leq \delta_T(y)\leq \Delta_T(y)\leq1$. If $\delta_T(y)\neq 0$ (resp. $\Delta_T(y)\neq 0$), we say $y$ is a Nevanlinna (resp. Valiron) exceptional boundary condition. Obviously by definition a Nevanlinna exceptional boundary condition is always a Valiron exceptional boundary condition. In particular, if $T_y$ has only finite eigenvalues, then $\delta_T(y)=1$ and $y$ is called a Picard exceptional boundary condition. The importance of determining whether a boundary condition is Valiron exceptional or not is due to:
\begin{proposition}(Abstract Weyl Law) If $y\in Gr(n,2n)$ is non-degenerate and not Valiron exceptional for the transcendental entire operator $T$, then
\[\lim_{r\rightarrow +\infty}\frac{N_T(r,y)}{h_T(r)}=1.\]
\end{proposition}
\begin{proof}This is clear from the First Main Theorem.
\end{proof}
\emph{Remark}. If there is a Picard exceptional boundary condition $y$ without any eigenvalue at all, we can choose another boundary condition $y'\in Gr(n,2n)$ such that $V_y\oplus V_{y'}=\mathcal{B}_T$. This gives rise to an affine coordinate chart in $Gr(n,2n)$ where $y$ and $y'$ are represented by matrices $(Id|0)$ and $(0|Id)$ respectively, while the Weyl curve is represented by $(W(\lambda)|Id)$ for a matrix-valued entire function $W(\lambda)$. In terms of $y$ and $y'$, a generic boundary condition $z\in Gr(n,2n)$ can be represented by two (not unique) constant $n\times n$ matrices $A, C$  such that $\textup{rk}(A|C)=n$. Then eigenvalues of $T_z$ are precisely zeros of the entire function
\[\det\left(
        \begin{array}{cc}
          W(\lambda) & Id \\
          A & C \\
        \end{array}
      \right).
\]

The following simple example shows (Picard) exceptional boundary conditions can really occur.
\begin{example}\label{ex1}Consider $T=-i\frac{d}{dx}$ with domain $H_0^1(I)\subset L^2(I)$ where $I=[0,1]$ and $H_0^1(I)$ is the usual Sobolev space
 \[\{f\in L^2(I)|f'\in L^2(I), f(0)=f(1)=0\}.\] It's a classical result that $T$ is a symmetric operator in $L^2(I)$ and $$D(T^*)=H^1(I):=\{f\in L^2(I)|f'\in L^2(I)\}.$$ $T$ is an entire operator with deficiency index 1. A natural boundary triplet can be chosen by setting $\Gamma_+\varphi=\varphi(1)/\sqrt{2}$ and $\Gamma_-\varphi=\varphi(0)/\sqrt{2}$ for $\varphi\in D(T^*)$. Then the corresponding contractive Weyl function is $B(\lambda)=e^{i\lambda}$. Thus the Weyl curve misses the two divisors $0$, $\infty\in \mathbb{CP}^1$. These two correspond to the boundary conditions $\varphi(1)=0$ and $\varphi(0)=0$ respectively. The corresponding extensions have no eigenvalues at all. Both of the two boundary conditions thus have (Nevanlinna or Valiron) defect 1 and all other generic boundary conditions should have infinite eigenvalues--this certainly can be derived by solving the corresponding boundary value problems explicitly, but our emphasis is that this is almost the general patten even if we cannot solve all the boundary value problems explicitly.  \end{example}

In the following we shall give explicit formulae for $h_T(r)$ and $m_T(r,y)$ in terms of the contractive Weyl function $B(\lambda)$ associated with a chosen boundary triplet. For $G\oplus_\bot G$, we choose an orthonormal basis $\{e_i\}_{i=1}^n$ for the first copy of $G$ but denote the same basis by $\{f_j\}_{j=1}^n$ for the second copy. Following the traditional notation, for an $n$ by $n$ matrix $F$ we identify $(Id|F)$ with the $n$-dimensional subspace of $G\oplus_\bot G$ spanned by $\{e_i+\sum_{j=1}^nF_{ij}f_j\}_{i=1}^n$. Then $\{(Id|F)|F\in M_{n\times n}\}$ forms a coordinate chart $U_+$ for $Gr(n,2n)$, where by $M_{n\times n}$ we mean the space of $n$ by $n$ complex matrices.  Similarly, $\{(F|Id)|F\in M_{n\times n}\}$ forms another coordinate chart $U_-$ for $Gr(n,2n)$, where $(F|Id)$ is identified with the subspace spanned by $\{\sum_{j=1}^nF_{ij}e_j+f_i\}_{i=1}^n$. Both $U_+$ and $U_-$ are affine coordinate charts in $Gr(n,2n)$ and the Weyl curve $W_T$ lies in $U_+\cup U_-$.
\begin{proposition}Each entry $B_{ij}(\lambda)$ of $B(\lambda)$ is a meromorphic function on $\mathbb{C}$ with $T_{B_{ij}}(r)\leq h_T(r)+O(1)$. Here $T_{B_{ij}}(r)$ is the characteristic function of $B_{ij}(\lambda)$ viewed as an entire curve in $\mathbb{CP}^1$.
\end{proposition}
\begin{proof}In terms of $\{e_i\}$ and $\{f_i\}$ as above, we can embed $Gr(n,2n)$ into $\mathbb{CP}^{C_{2n}^n-1}$ (the so-called Pl$\ddot{u}$cker embedding). Then the entries of $B(\lambda)$ are essentially part of the Pl$\ddot{u}$cker coordinates of $W_T(\lambda)$. Note that $h_T(r)$ is precisely the characteristic function of $W_T(\lambda)$ embedded as a curve in $\mathbb{CP}^{C_{2n}^n-1}$. The conclusion then follows from \cite[Thm.~A5.1.2]{ru2021nevanlinna}.
\end{proof}

Note that $\{s_i:=e_i+\sum_{j=1}^nF_{ij}f_j\}_{i=1}^n$ is a local frame over $U_+$ for the tautological bundle $E$ over $Gr(n, 2n)$. Then $\sigma:=s_1\wedge\cdots \wedge s_n$ is a local frame of $l$. The inner product on $G\oplus_\bot G$ introduces a natural Hermitian metric $h$ on $l$; in particular,
\[h(\sigma, \sigma)=\det ((s_i,s_j))=\det(Id+FF^*).\]
Then the first Chern form of $L$ can be expressed over $U_+$ by
\[c_1(L,h)=dd^c\ln \det(Id+FF^*),\]
where $d^c=\frac{i}{4\pi}(\bar{\partial}-\partial)$. A similar formula over $U_-$ holds. In this way, over a neighbourhood of $\overline{\mathbb{C}_+}$ (the closure of $\mathbb{C}_+$),
\[W_T^*(c_1(L, h))=dd^c\ln \det(Id+B(\lambda)B(\lambda)^*).\]
Similarly, over a neighbourhood of $\overline{\mathbb{C}_-}$ in $\mathbb{C}$,
\[W_T^*(c_1(L,h))=dd^c\ln \det(Id+B(\overline{\lambda})^*B(\overline{\lambda}))=dd^c\ln \det(Id+B(\overline{\lambda})B(\overline{\lambda})^*),\]
where the second equality follows from the fact that $BB^*$ and $B^*B$ have the same nonzero eigenvalues. For $\lambda\in \overline{\mathbb{C}_+}$, we write $\lambda=u+iv$ where $u,v\in \mathbb{R}$.
\begin{theorem}In terms of the contractive Weyl function $B(\lambda)$,
\begin{equation}h_T(r)=\frac{1}{2\pi i}\int_0^r\frac{dt}{t}\int_{-t}^t(\ln\det B(u))'du+O(1),\label{e3}\end{equation}
where $O(1)$ is a bounded term in $r$.
\end{theorem}
\begin{proof}We use $B_t^+$ to denote the upper half of the disc $B_t$. It is not hard to find that
\[\int_{B_t}W_T^*(c_1(L,h))=2\int_{B_t^+}W_T^*(c_1(L,h))=2\int_{B_t^+}dd^c\ln \det(Id+B(\lambda)B(\lambda)^*).\]
By Stokes' theorem, we have
\begin{eqnarray*}&\quad&\int_{B_t^+}dd^c\ln \det(Id+B(\lambda)B(\lambda)^*)=\int_{\partial B_t^+}d^c\ln \det(Id+B(\lambda)B(\lambda)^*)\\
&=& \int_{C_t^+}d^c\ln \det(Id+B(\lambda)B(\lambda)^*)+\int_{I_t}d^c\ln \det(Id+B(\lambda)B(\lambda)^*),
\end{eqnarray*}
where $C_t^+$ is the semicircle part of $\partial B_t^+$ while $I_t$ is the interval $[-t, t]$ on the real line. We note that in polar coordinates for a smooth function $f$ on $\mathbb{C}$
\[d^cf=\frac{1}{4\pi}(\rho\frac{\partial f}{\partial \rho}d\theta-\frac{1}{\rho}\frac{\partial f}{\partial \theta}d\rho),\]
while in rectangular coordinates,
\[d^c f=\frac{1}{4\pi}(\frac{\partial f}{\partial u}dv-\frac{\partial f}{\partial v}du).\]
Therefore,
\[4\pi\int_{C_t^+}d^c\ln \det(Id+B(\lambda)B(\lambda)^*)=\int_0^\pi t\frac{\partial}{\partial t}\ln \det(Id+B(te^{i\theta})B(te^{i\theta})^*)d\theta\]
and
\begin{eqnarray*}&\quad&4\pi\int_0^r\frac{dt}{t}\int_{C_t^+}d^c\ln \det(Id+B(\lambda)B(\lambda)^*)\\&=&\int_0^rdt\frac{d}{dt}\int_0^\pi \ln \det(Id+B(te^{i\theta})B(te^{i\theta})^*)d\theta\\
&=&\int_0^\pi \ln \det(Id+B(re^{i\theta})B(re^{i\theta})^*)d\theta\\
&-&\pi\ln \det(Id+B(0)B(0)^*)=O(1). \end{eqnarray*}
The last equality is because $\|B(\lambda)\|<1$ for $\lambda\in \mathbb{C}_+$.

On the other side,
\[4\pi \int_{I_t}d^c\ln \det(Id+B(\lambda)B(\lambda)^*)=-\int_{-t}^t[\frac{\partial}{\partial v}\ln \det(Id+B(\lambda)B(\lambda)^*)]|_{v=0}du.\]
Note that
\[\frac{\partial}{\partial v}\ln \det(Id+B(\lambda)B(\lambda)^*)=\textup{Tr}[(Id+B(\lambda)B(\lambda)^*)^{-1}\frac{\partial}{\partial v}(B(\lambda)B(\lambda)^*)]\]
and
\[\frac{\partial}{\partial v}(B(\lambda)B(\lambda)^*)=iB'(\lambda)B(\lambda)^*-iB(\lambda)B'(\lambda)^*.\]
Using the fact that on the real line $B(u)$ is unitary, we find
\[[\frac{\partial}{\partial v}\ln \det(Id+B(\lambda)B(\lambda)^*)]|_{v=0}=i\textup{Tr}(B'(u)B(u)^*).\]
Note that $\textup{Tr}(B'(u)B(u)^*)=(\ln\det B(u))'$. The claimed formula then follows.
\end{proof}
\emph{Remark}. From the proof we can see the bound of the bounded term $O(1)$ can be chosen to depend only on the deficiency index $n$.

The formula (\ref{e3}) has a clear spectral theoretic meaning. Following the remark after Prop.~\ref{p4}, we assume $n=1$, then $B(u)=e^{i\theta(u)}$ for a strictly increasing phase function $\theta(u)$. Note that the zeros of $e^{i\theta(u)}-1=0$ are precisely the eigenvalues of $T_0$. Then
\begin{equation}h_T(r)=\frac{1}{2\pi}\int_0^r\frac{\theta(t)-\theta(-t)}{t}dt+O(1).\label{eq8}\end{equation}
To some extent, $h_T(r)$ measures how much the phase $\theta(u)$ changes when $u$ goes through the interval $[-r, r]$. The more $\theta(u)$ changes, the more eigenvalues $T_0$ will have in the interval $[-r,r]$. In the general case, $\det B(u)$ takes values in $\mathbb{U}(1)$ and we can lift it to a real-valued total phase function $\phi(u)$ such that $\det B(u)=e^{i\phi(u)}$ and hence $(\ln \det B(u))'=i\phi'(u)$. In particular, an equation similar to Eq.~(\ref{eq8}) holds.
\begin{proposition}\label{p18}If $T$ is a transcendental entire operator with deficiency index 1, then the spectra of self-adjoint extensions of $T$ are pairwise interlaced.
\end{proposition}
\begin{proof}
If $\tilde{T}$ and $\tilde{T}'$ are any two self-adjoint extensions of $T$, then a boundary triplet can be chosen such that $\tilde{T}=T_0$ and $\tilde{T}'=T_1$ (\cite[Thm.~2.5.9]{behrndt2020boundary}). In terms of the phase function $\theta(u)$ in the above argument, we see that the eigenvalues of $\tilde{T}$ and $\tilde{T}'$ are precisely the zeros of $e^{i\theta(u)}-1=0$ and $e^{i\theta(u)}+1=0$ respectively. The claim then follows since $\theta(u)$ is analytic and strictly increasing.
\end{proof}
This result is not new, but we think our proof is conceptually very clear. The proposition will be generalized in \S\S~\ref{spe}.
\begin{theorem}Let $y\in Gr(n,2n)$ be represented by $(Id|\Omega)$ in $U_+$. Then in terms of the contractive Weyl function $B(\lambda)=B(re^{i\theta})$,
\begin{equation}
m_T(r,y)=-\int_0^\pi\ln [|\det(B-\Omega)\det(Id-B^*\Omega)|]\frac{d\theta}{2\pi}+O(1).
\end{equation}
In particular, if $\Omega=Id$, then
\[m_T(r,y)=-\int_0^\pi\ln |\det(B-Id)|\frac{d\theta}{\pi}+O(1).\]
\end{theorem}
\begin{proof}Now $V_y=\textup{span}\{\sigma_i:=e_i+\sum_{j=1}^n\Omega_{ij}f_j,\, j=1,\cdots, n.\}$. Then we can set $s_y=\sigma_1\wedge\cdots \sigma_n$ and by definition on $U_+$
\[s_y(\sigma)=\det\left(
                         \begin{array}{cc}
                           Id & \Omega \\
                           Id & F \\
                         \end{array}
                       \right)=\det(F-\Omega),
\]
implying
\[\|s_y\|^2=\frac{|\det(F-\Omega)|^2}{\det(Id+FF^*)}.\]
Similarly, on $U_-$
\[\|s_y\|^2=\frac{|\det(Id-F\Omega)|^2}{\det(Id+FF^*)}.\]
Then by definition $-2m_T(r,y)$ equals to
\begin{eqnarray*}&\quad&\int_0^\pi\ln \frac{|\det(B(re^{i\theta})-\Omega)|^2}{\det(Id+B(re^{i\theta})B(re^{i\theta})^*)}\frac{d\theta}{2\pi}+\int_\pi^{2\pi}\ln \frac{|\det(Id-B(re^{-i\theta})^*\Omega)|^2}{\det(Id+B(re^{-i\theta})^*B(re^{-i\theta}))}\frac{d\theta}{2\pi}\\
&=&\int_0^\pi\ln \frac{|\det(B(re^{i\theta})-\Omega)|^2}{\det(Id+B(re^{i\theta})B(re^{i\theta})^*)}\frac{d\theta}{2\pi}+\int_0^{\pi}\ln \frac{|\det(Id-B(re^{i\theta})^*\Omega)|^2}{\det(Id+B(re^{i\theta})^*B(re^{i\theta}))}\frac{d\theta}{2\pi}\\
&=&\int_0^\pi \ln[|\det(B(re^{i\theta})-\Omega)\det(Id-B(re^{i\theta})^*\Omega)|^2]\frac{d\theta}{2\pi}+O(1).\end{eqnarray*}
The last equality is again because $\|B(\lambda)\|<1$ for any $\lambda\in \mathbb{C}_+$.
\end{proof}
\begin{example}In Example \ref{ex1}, we have $B(\lambda)=e^{i\lambda}$, and thus $h_T(r)=\frac{r}{\pi}+O(1)$. In particular, $T$ is of Weyl order 1 and has Weyl type $1/\pi$ in the sense of Def.~\ref{d1}.
\end{example}

For solving boundary value problems, it's important to investigate the structure of the set of exceptional boundary conditions. We denote the sets of Picard, Nevanlinna and Valiron exceptional boundary conditions by $\mathrm{B}_P(T)$, $\mathrm{B}_N(T)$ and $\mathrm{B}_V(T)$ respectively. By definition, $\mathrm{B}_P(T)\subseteq \mathrm{B}_N(T)\subseteq \mathrm{B}_V(T)$. We claim that these "exceptional" boundary conditions are really exceptional in the following sense.
\begin{theorem}\label{thm9}If the entire operator $T$ is transcendental and weakly algebraically non-degenerate, then $\mathrm{B}_V(T)$ is of measure zero w.r.t. to any smooth measure on $Gr(n,2n)$.
\end{theorem}

The proof of the theorem depends on a Crofton formula for Schubert hyperplanes in $Gr(n,2n)$. A similar formula for entire curves in $\mathbb{CP}^n$ is well-known and a basic result in value distribution theory. We equip $Gr(n,2n)$ with the K$\ddot{a}$hler structure $\omega$ by fixing a boundary triplet. We normalize $\mu_y=\omega^{n^2}$ such that $\int_{Gr(n,2n)}\mu_y=1$.

\begin{lemma}(\cite{griffiths1978complex})If $f$ is a non-constant entire curve in $Gr(n,2n)$, then
\begin{equation}T_f(r, L)=\int_{Gr(n,2n)}N_f(r, (s_y))\mu_y.
\end{equation}
\end{lemma}
\begin{proof} In \cite{griffiths1978complex} Griffiths has actually proved such a formula for general Schubert cycles in general Grassmannians. For the reader's convenience, we outline Griffiths' argument in this specific setting. Since the locus of degenerate boundary conditions in $Gr(n, 2n)$ is of codimension at least 2, w.l.g, we can assume $f$ is weakly algebraically non-degenerate.

We have the Poincar$\acute{e}$-Lelong formula
\[dd^c[\ln \|s_y\|^2]=-c_1(L,h)+[(s_y)].\]
By Stokes' Theorem, we have
\[\int_{\partial B_t}d^c[\ln \|s_y(f(\lambda))\|^2]=-\int_{B_t}f^*(c_1(L,h))+n_f(r,(s_y)).\]
Integrating the above formula over $Gr(n,2n)$, and interchanging the order of integration, we get
\[\int_{\partial B_t}f^*(\eta)=-\int_{B_t}f^*(c_1(L,h))+\int_{Gr(n,2n)}n_f(r,(s_y))\mu_y,\]
where $\eta=\int_{Gr(n,2n)}d^c[\ln \|s_y\|^2]\mu_y$ is a current on $Gr(n,2n)$ of degree 1, i.e., for any smooth $(2n^2-1)$-form $\alpha$,
\[\eta(\alpha)=\int_{Gr(n,2n)\times Gr(n,2n)}d^c[\ln \|s_y\|^2]\wedge \alpha \mu_y.\]
We claim that $\eta\equiv0$ (and the lemma then follows immediately). To see this we only need to show $\eta$ is $\mathbb{PU}(2n)$-invariant. This is because $Gr(n,2n)$ is a compact Hermitian symmetric space (in particular a Riemannian symmetric space) whose invariant currents are precisely the harmonic forms. Since harmonic forms of odd degree on $Gr(n,2n)$ are trivial, $\eta$ has to be zero for it's of degree 1.

For any $g\in \mathbb{U}(2n)$, if $\sigma_x$ is a normalized basis of $l_x$ for $x\in Gr(n, 2n)$, then by definition,
\[s_{g\cdot y}(g\cdot \sigma_x)=\det g \times s_y(\sigma_x).\]
Since $g$ preserves the inner product, $g\cdot \sigma_x$ is a normalized basis of $l_{g\cdot x}$ and hence
\[\|s_y\|(x)=\|s_{g\cdot y}\|(g\cdot x),\quad \textup{or} \quad (g^{-1})^*(\|s_y\|)=\|s_{g\cdot y}\|.\]
Consequently,
\begin{eqnarray*}(g^*\eta)(\alpha)&=&\int_{Gr(n,2n)}\mu_y\int_{Gr(n,2n)}d^c[\ln \|s_y\|^2]\wedge (g^{-1})^*\alpha\\
&=&\int_{Gr(n,2n)}\mu_y\int_{Gr(n,2n)}d^c[\ln \|s_{g\cdot y}\|^2]\wedge (g^{-1})^*\alpha\\
&=&\int_{Gr(n,2n)}\mu_y\int_{Gr(n,2n)}(g^{-1})^*[d^c[\ln \|s_y\|^2]\wedge \alpha]\\
&=&\int_{Gr(n,2n)}\mu_y\int_{Gr(n,2n)}d^c[\ln \|s_y\|^2]\wedge \alpha=\eta(\alpha),
\end{eqnarray*}
where the invariance of $\mu_y$ and the connectedness of $\mathbb{PU}(2n)$ are used. Thus $\eta$ is really $\mathbb{PU}(2n)$-invariant.
\end{proof}
We state an inequality which comes from an easy refinement of the First Main Theorem. Its proof can be obtained by a slight modification of that of \cite[Thm.~5.1]{Ochiai1974some}. The notation in \cite{Ochiai1974some} is a bit different. Ochiai's $r$ is actually our $\ln r$, his $T_f(r)$ is our $T_f(e^r)-T_f(1)$. Similarly his $N_f(\phi, r)$ is our $N_f(e^r, (\phi))-N_f(1,(\phi))$ where $(\phi)$ is the divisor provided by a section $\phi$ of $L$. Thus for $r>1$ we define
\[\tilde{h}_T(r)=h_T(r)-h_T(1),\quad \tilde{N}_T(r,y)=N_T(r,y)-N_T(1,y).\]
Then we still have the Crofton formula
\[\tilde{h}_T(r)=\int_{Gr(n,2n)}\tilde{N}_T(r,y)\mu_y.\]
As Thm.~5.1 in \cite{Ochiai1974some}, we have
\begin{lemma}\label{Oh}If $T$ is a transcendental and weakly algebraically non-degenerate entire operator with deficiency index $n$, then there is a positive constant $K$ such that for $r>1$ and any $y\in Gr(n,2n)$
\[\tilde{N}_T(r,y)<\tilde{h}_T(r)+K.\]
Note that $K$ is independent of $r$ and $y$.
\end{lemma}
\begin{proof}In terms of $\tilde{h}_T(r)$ and $\tilde{N}_T(r,y)$, the First Main Theorem takes the following form (integrating the Poincar$\acute{e}$-Lelong formula over the annulus $1<|\lambda|<r$)
\[\tilde{h}_T(r)=\tilde{N}_T(r,y)+m_T(r,y)-m_T(1,y),\quad r>1.\]
Note that $s_y$ is only determined by $y\in Gr(n, 2n)$ up to a nonzero factor. W.l.g, we can require that $\max_{x\in Gr(n,2n)}\|s_y(x)\|=1$. Then $m_T(r,y)\geq 0$ and
\[\tilde{N}_T(r,y)\leq \tilde{h}_T(r)+m_T(1,y).\]
To be precise, we write $m_T(1,y)$ as $m_T(1,s_y)$, which can be viewed as a function on the set
\[\mathcal{R}:=\{s_y|y\in Gr(n,2n),\, \max_{x\in Gr(n,2n)}\|s_y(x)\|=1\}.\]
This set is actually a circle bundle over $Gr(n,2n)$ and thus compact. Then the conclusion follows if one can prove $m_T(1,s_y)$ is continuous on $\mathcal{R}$. This continuity can be proved along the same line of the proof of Lemma 5.1 in \cite{Ochiai1974some} and we refer the reader to it for the details.
\end{proof}

After these preliminary works, we can now turn to the proof of Thm.~\ref{thm9}. We follow closely the idea in the proof of \cite[Thm.~2.1, Chap.~4]{gol'dberg2008value}.

\begin{proof}Proof of Thm.~\ref{thm9}. We use the smooth measure $\mu_y$ on $Gr(n,2n)$.

If $D\subset Gr(n,2n)$ is measurable, by integration we obtain
\[\int_{D}\tilde{N}_T(r,y)\mu_y\leq \tilde{h}_T(r)|D|+K,\]
where $|D|$ is the $\mu_y$-measure of $D$. Let $\eta\in (0,1/2)$ and for $r>1$
\[B(r):=\{y\in Gr(n,2n)|\tilde{N}_T(r,y)\leq \tilde{h}_T(r)-\tilde{h}_T^{\frac{1}{2}+\eta}(r)\}.\]
Since $\tilde{N}_T(r,y)$ is measurable as a function of $y$, $B(r)$ is measurable. We claim that
\[|B(r)|<\frac{2K}{\tilde{h}_T^{\frac{1}{2}+\eta}(r)}.\]
Assume the contrary, i.e., $|B(r)|\geq \frac{2K}{\tilde{h}_T^{\frac{1}{2}+\eta}(r)}$. Let $B_1(r)\subset B(r)$ such that
\[|B_1(r)|=\frac{2K}{\tilde{h}_T^{\frac{1}{2}+\eta}(r)}\]
and $C(r)=Gr(n,2n)\backslash B_1(r)$. Then by the Crofton formula and Lemma \ref{Oh},
\begin{eqnarray*}\tilde{h}_T(r)&=&\int_{B_1(r)}\tilde{N}_T(r,y)\mu_y+\int_{C(r)}\tilde{N}_T(r,y)\mu_y\\
&\leq& [\tilde{h}_T(r)-\tilde{h}_T^{\frac{1}{2}+\eta}(r)]|B_1(r)|+\tilde{h}_T(r)|C(r)|+K\\
&=&\tilde{h}_T(r)-\tilde{h}_T^{\frac{1}{2}+\eta}(r)|B_1(r)|+K\\
&=&\tilde{h}_T(r)-K,
\end{eqnarray*}
a contradiction!

Note that for $x>1$, the function $x-x^{\frac{1}{2}+\eta}$ is strictly increasing. We can choose $r_1>1$ such that $\tilde{h}_T(r_1)>1$, and define a sequence $\{r_n\}$ by the following recurrence: $r_n$ is uniquely determined by the solution of
\[\tilde{h}_T(r)-\tilde{h}_T^{\frac{1}{2}+\eta}(r)=\tilde{h}_T(r_{n-1}),\quad n=2,3,\cdots.\]
Since $\tilde{h}_T(r_n)-\tilde{h}_T(r_{n-1})=\tilde{h}_T^{\frac{1}{2}+\eta}(r_n)>1$, we know that as $n$ approaches $+\infty$, $\tilde{h}_T(r_n)\rightarrow +\infty$ and $r_n$ goes to $+\infty$ increasingly.

If $y$ is not in $B(r_n)$, then for $r_n\leq r\leq r_{n+1}$, we have
\begin{eqnarray*}
\tilde{N}_T(r,y)&\geq& \tilde{N}_T(r_n,y)> \tilde{h}_T(r_n)-\tilde{h}_T^{\frac{1}{2}+\eta}(r_n)\\
&\geq& \tilde{h}_T(r_n)-\tilde{h}_T^{\frac{1}{2}+\eta}(r_n)+[\tilde{h}_T(r)-\tilde{h}_T^{\frac{1}{2}+\eta}(r)]-[\tilde{h}_T(r_{n+1})-\tilde{h}_T^{\frac{1}{2}+\eta}(r_{n+1})]\\
&=&\tilde{h}_T(r)-\tilde{h}_T^{\frac{1}{2}+\eta}(r)-\tilde{h}_T^{\frac{1}{2}+\eta}(r_n)\geq \tilde{h}_T(r)-2\tilde{h}_T^{\frac{1}{2}+\eta}(r).
\end{eqnarray*}
Hence if $y\in Gr(n,2n)\backslash \cup_{k=n}^\infty B(r_k)$, then for any $r>r_n$ we shall have
\[\tilde{N}_T(r,y)\geq\tilde{h}_T(r)-2\tilde{h}_T^{\frac{1}{2}+\eta}(r),\]
and consequently
\[\liminf_{r\rightarrow +\infty}\frac{N_T(r,y)}{h_T(r)}=\liminf_{r\rightarrow +\infty}\frac{\tilde{N}_T(r,y)}{\tilde{h}_T(r)}\geq 1,\]
which surely implies $\liminf_{r\rightarrow +\infty}\frac{N_T(r,y)}{h_T(r)}=1$ and $y$ is not Valiron exceptional. So we have $\mathrm{B}_V(T)\subset \cup_{k=n}^\infty B(r_k)$ for any $n\in \mathbb{N}$.

On the other hand, we can obtain
\begin{eqnarray*}
|\cup_{k=n}^\infty B(r_k)|&\leq& \sum_{k=n}^\infty|B(r_k)|\leq 2K \sum_{k=n}^\infty \tilde{h}_T^{-\frac{1}{2}-\eta}(r_k)
=2K\sum_{k=n}^\infty\frac{\tilde{h}_T(r_k)-\tilde{h}_T(r_{k-1})}{\tilde{h}_T^{1+2\eta}(r_k)}\\
&=&2K\sum_{k=n}^\infty\int_{r_{k-1}}^{r_k}\frac{d\tilde{h}_T(r)}{\tilde{h}_T^{1+2\eta}(r_k)}\leq 2K\sum_{k=n}^\infty\int_{r_{k-1}}^{r_k}\frac{d\tilde{h}_T(r)}{\tilde{h}_T^{1+2\eta}(r)}\\
&=&2K\int_{r_{n-1}}^\infty \frac{d\tilde{h}_T(r)}{\tilde{h}_T^{1+2\eta}(r)}=\frac{K}{\eta}\frac{1}{\tilde{h}_T^{2\eta}(r_{n-1})}.
\end{eqnarray*}
This shows $|\mathrm{B}_V(T)|=0$ as required.
\end{proof}
\emph{Remark}. The proof only uses the fact that the Weyl curve is transcendental and weakly algebraically non-degenerate. Thus the theorem holds for any transcendental and weakly algebraically non-degenerate entire curve in $Gr(n,2n)$.

The distribution of exceptional boundary conditions in $Gr(n,2n)$ is, to some extent, controlled by the Second Main Theorems, which are viewed as the much deeper results in value distribution theory.

For a finite set $I$, by $\sharp I$ we mean the number of elements in $I$.

\begin{definition}Let $X$ be a complex projective manifold of dimension $m\geq 1$. Let $D_1,\cdots, D_q$ be effective divisors on $X$. These divisors are said to be in general position in $X$ if for any subset $I\subset \{1,2,\cdots, q\}$ with $\sharp I\leq m+1$
\[\textup{dim}(\cap_{j\in I}\textup{supp}D_j)\leq m-\sharp I,\]
where $\textup{supp}D$ means the support of the divisor $D$. Here by convention, $\textup{dim}\varnothing=-1$.
\end{definition}
A basic result of Ru in \cite{ru2009holomorphic} is the following Second Main Theorem on entire curves in projective algebraic varieties.
\begin{theorem}
Let $X$ be a complex projective manifold of dimension $m\geq 1$. Let $D_1,\cdots, D_q$ be effective divisors on $X$, located in general position. Suppose that there exists an ample divisor $A$ on $X$ and positive integers $d_j$ such that $D_j\sim d_j A$ for $j=1,2,\cdots, q$. Let $f: \mathbb{C}\rightarrow X$ be an algebraically non-degenerate entire curve. Then for every $\varepsilon>0$,
\[\sum_{j=1}^qd_j^{-1}m_f(r, D_j)\leq (m+1+\varepsilon)T_{f, \mathcal{O}(A)}(r)||_E.\]
\end{theorem}
Here $||_E$ means the above inequality holds except for $r>0$ in a set of finite Lebesgue measure.

\begin{definition}$y_i\in Gr(n, 2n)$, $i=1,\cdots, q$, are said to be boundary conditions in general position if the corresponding divisors $D_i=(s_{y_i})$ are in general position in $Gr(n,2n)$.
\end{definition}

Note that all $(s_y)$ are linearly equivalent. Applying Ru's result to our context, we have
\begin{theorem}\label{thm7}If $T$ is a transcendental entire operator with deficiency index $n$, which is algebraically non-degenerate,  and $y_i\in Gr(n, 2n)$, $i=1,\cdots, q$ ($q>n^2$) are boundary conditions in general position, then for any given $\varepsilon>0$
\begin{equation}
\sum_{j=1}^qm_T(r, y_j)\leq (n^2+1+\varepsilon)h_T(r)||_E.
\end{equation}
\end{theorem}
As a consequence, we have the following defect inequality.
\begin{corollary}Let the assumption of Thm.~\ref{thm7} hold. Then
\begin{equation}\label{e4}\sum_{j=1}^q\delta_T(y_j)\leq n^2+1.\end{equation}
\end{corollary}
This inequality puts a constraint on the number of Nevanlinna exceptional boundary conditions in general position. For instance, if $n=1$, as soon as these $y_i$'s are different from each other, they are in general position. Then there are at most countably infinite Nevanlinna exceptional boundary conditions. In particular, the number of Picard exceptional boundary conditions can at most be 2. This latter statement, of course, is only a disguised form of Picard's famous theorem.

If the Weyl curve $W_T(\lambda)$ is algebraically degenerate, then it has more chances to omit divisors; see an example in \S~\ref{sturm}. In this case, for $n>1$, applying Thm.~A7.3.4 in \cite{ru2021nevanlinna} (see also the remark after this theorem), the right hand side of the above defect inequality (\ref{e4}) should be replaced by $n^2(n^2+1)/2$.

One should note that there are some gaps between Ru's Second Main Theorem and our setting: what interest us are not divisors associated to general holomorphic sections of $L$ but Schubert hyperplanes. Thus we'd better use the notion of weak algebraic non-degeneracy rather than that of algebraic non-degeneracy. This is closer to the context of H. Cartan's original work on entire curves in projective space $\mathbb{CP}^n$ (see Thm. A5.1.7 in \cite{ru2021nevanlinna}), where only hyperplanes are considered and the notion of linear non-degeneracy is enough. Further work is needed in this respect.

\subsection{Multiplicities of eigenvalues}
via the discussion before, it's possible to translate the problem of counting eigenvalues into a problem of value distribution theory. However, there is a subtlety that must be dealt with first: for a self-adjoint operator, an eigenvalue $\lambda_i$ is always counted according to its \emph{geometric multiplicity}, i.e., the dimension of the corresponding eigensubspace, while in value distribution theory a zero is always counted according to the analytic multiplicity, i.e., the multiplicity of $\lambda_i$ as the zero of an analytic function. Do the two notions of multiplicity coincide?
\begin{theorem}\label{thm8}Assume that $T$ is an entire operator with deficiency index $n$ and that $T_0$ is a self-adjoint extension of $T$. Then for any eigenvalue of $T_0$, the geometric multiplicity and the analytic multiplicity coincide.
\end{theorem}
\begin{proof}Though an alternative proof can be given if we can establish the following Thm.~\ref{thm17}, we prefer to provide an independent proof here.

We can choose a boundary triplet such that $T_0$ is parameterized by the identity matrix $Id$. Let $B(\lambda)$ be the corresponding contractive Weyl function. Then $\lambda_0\in \mathbb{R}$ is an eigenvalue of $T_0$ if and only if $Id-B(\lambda_0)$ is not invertible, i.e., $\lambda_0$ is a zero of the analytic function $\det(Id-B(\lambda))$. Let $n_a(\lambda_0)$ (resp. $n_g(\lambda_0)$) denote the analytic (resp. geometric) multiplicity of $\lambda_0$. Then by Coro.~\ref{coro1},
\[n_g(\lambda_0)=n-\textup{rk}(Id-B(\lambda_0)).\]
We want to calculate $n_a(\lambda_0)$. Note that $n_a(\lambda_0)$ is the residue of $\frac{d}{d\lambda}\ln\det(Id-B(\lambda))$ at $\lambda_0$ and that for $\lambda\in \mathbb{R}$
\begin{eqnarray*}\frac{d}{d\lambda}\ln\det(Id-B(\lambda))&=&-\textup{Tr}[(Id-B(\lambda))^{-1}\frac{dB}{d\lambda}]\\
&=&-i\textup{Tr}[(Id-B(\lambda))^{-1}B(\lambda)\gamma_+(\lambda)^*\gamma_+(\lambda)]\\
&=&-i\textup{Tr}[(B(\lambda)^{-1}-Id))^{-1}\gamma_+(\lambda)^*\gamma_+(\lambda)]\end{eqnarray*}
due to Prop.~\ref{p4}. From our resolvent formula in Thm.~\ref{reso}, one can get
\begin{eqnarray*}
\Gamma_-[(T_0-\lambda)^{-1}-(T_+-\lambda)^{-1}]\gamma_+(\lambda)&=&i\Gamma_-\gamma_+(\lambda)(B(\lambda)^{-1}-Id)^{-1}\gamma_+(\lambda)^*\gamma_+(\lambda)\\
&=&iB(\lambda)(B(\lambda)^{-1}-Id)^{-1}\gamma_+(\lambda)^*\gamma_+(\lambda).
\end{eqnarray*}
Then we get
\[i(B(\lambda)^{-1}-Id)^{-1}\gamma_+(\lambda)^*\gamma_+(\lambda)=B(\lambda)^{-1}\Gamma_-[(T_0-\lambda)^{-1}-(T_+-\lambda)^{-1}]\gamma_+(\lambda).\]
So we have
\[\frac{d}{d\lambda}\ln\det(Id-B(\lambda))=-\textup{Tr}(B(\lambda)^{-1}\Gamma_-[(T_0-\lambda)^{-1}-(T_+-\lambda)^{-1}]\gamma_+(\lambda)).\]
Note that
\[n_a(\lambda_0)=\lim_{\lambda\rightarrow \lambda_0}(\lambda-\lambda_0)\frac{d}{d\lambda}\ln\det(Id-B(\lambda)).\]
Let $P_{\lambda_0}$ be the orthogonal projection onto the eigensubspace of $T_0$ corresponding to $\lambda_0$. Then
\[\lim_{\lambda\rightarrow \lambda_0}\textup{Tr}((\lambda-\lambda_0)B(\lambda)^{-1}\Gamma_-[(T_0-\lambda)^{-1}(Id-P_{\lambda_0})-(T_+-\lambda)^{-1}]\gamma_+(\lambda))=0\]
because the term
\[B(\lambda)^{-1}\Gamma_-[(T_0-\lambda)^{-1}(Id-P_{\lambda_0})-(T_+-\lambda)^{-1}]\gamma_+(\lambda)\]
 has no singularity at $\lambda_0$ and the finite-dimensional trace is continuous. Thus
\begin{eqnarray*}\lim_{\lambda\rightarrow \lambda_0}(\lambda-\lambda_0)\frac{d}{d\lambda}\ln\det(Id-B(\lambda))&=&\lim_{\lambda\rightarrow \lambda_0}\textup{Tr}(B(\lambda)^{-1}\Gamma_-P_{\lambda_0}\gamma_+(\lambda))\\
&=&\textup{Tr}(B(\lambda_0)^{-1}\Gamma_-P_{\lambda_0}\gamma_+(\lambda_0))\\
&=&\textup{Tr}(B(\lambda_0)^{-1}\Gamma_-\gamma_+(\lambda_0)\Gamma_+P_{\lambda_0}\gamma_+(\lambda_0)).\end{eqnarray*}
The last equality is because on $\textup{ker}(T^*-\lambda)$, $\gamma_+(\lambda_0)\Gamma_+$ is the identity. Furthermore,
\begin{eqnarray*}&\quad&\textup{Tr}(B(\lambda_0)^{-1}\Gamma_-\gamma_+(\lambda_0)\Gamma_+P_{\lambda_0}\gamma_+(\lambda_0))=\textup{Tr}(B(\lambda_0)^{-1}B(\lambda_0)\Gamma_+P_{\lambda_0}\gamma_+(\lambda_0))\\
&=&\textup{Tr}(\Gamma_+P_{\lambda_0}\gamma_+(\lambda_0))=\textup{Tr}(\gamma_+(\lambda_0)\Gamma_+P_{\lambda_0})\\
&=&\textup{Tr}(P_{\lambda_0})=n_g(\lambda_0).
\end{eqnarray*}
Note that "$\textup{Tr}$" in the last two lines may be the trace of operators in an infinite-dimensional Hilbert space. However, the third equality holds because the operator after $\textup{Tr}$ is of finite rank.
\end{proof}
To the best of our knowledge, in the literature, the analytic multiplicity is never defined in this generality and the above theorem has only been established for some special ordinary differential operators.

\begin{proposition}An entire operator is rational if and only if $\textup{dim}H<+\infty$.
\end{proposition}
\begin{proof}If $\textup{dim}H<+\infty$, the rationality of $W_T(\lambda)$ can be easily seen by choosing a basis of $H$ and computing $\textup{ker}(T^*-\lambda)$ explicitly.

Conversely, if $W_T(\lambda)$ is rational, then the contractive Weyl function $B(\lambda)$ is rational on $\mathbb{C}$ for any chosen boundary triplet. The equation $\det(B(\lambda)-Id)=0$ has only finite real zeros, counting (analytic) multiplicities. According to the above theorem and the spectral resolution of $T_0$, $H$ is of finite dimension.
\end{proof}
\begin{corollary}An entire operator $T$ with deficiency index $n$ is transcendental if and only if one (thus any) of its self-adjoint extensions has countably infinite eigenvalues.
\end{corollary}
\begin{proof}If $T$ is transcendental, then $\textup{dim}H=+\infty$. Let $\tilde{T}$ be any self-adjoint extension of $T$. Since each eigenvalue of $\tilde{T}$ is of finite multiplicity, the spectral decomposition of $\tilde{T}$ implies that it has countably infinite eigenvalues. The converse is obvious.
\end{proof}

Recall the definition of algebraic multiplicity of an eigenvalue of a (bounded or unbounded) operator $A$ in a Hilbert space $H$; see \cite{behrndt2022generalized}.
\begin{definition}Vectors $\varphi_0,\cdots, \varphi_{k-1}\in D(A)$ form a Jordan chain of length $k$ for $A$ at $\lambda_0\in \sigma_p(A)$ if and only if $\varphi_0\neq 0$ and
\[(A-\lambda_0)\varphi_0=0,\quad (A-\lambda_0)\varphi_j=\varphi_{j-1},\quad j=1,\cdots, k-1.\]
If $\{\varphi_{0,j}\}_{1\leq j\leq N}$, $N\in \mathbb{N}\cup \{+\infty\}$, is a basis of $\textup{ker}(A-\lambda_0)$ and $\{\varphi_{0,j},\cdots, \varphi_{k_j-1,j}\}$ the corresponding Jordan chain of maximal length $k_j$. Then the algebraic multiplicity $m_{al}(\lambda_0)$ of $\lambda_0$ is $\sum_{j=1}^Nk_j$. If one of the eigenvectors has a Jordan chain of arbitrary length, we simply define $m_{al}(\lambda_0)=+\infty$.
\end{definition}
By definition, we have $m_g(\lambda_0)\leq m_{al}(\lambda_0)$. If $\lambda_0$ is an isolated point in $\sigma(A)$, we can define
\[P(\lambda_0)=\frac{i}{2\pi}\oint_{\partial D(\lambda_0, \varepsilon)}(A-\lambda)^{-1}d\lambda,\]
where $\partial D(\lambda_0, \varepsilon)$ is the positively oriented boundary of the disc $D(\lambda_0, \varepsilon)\subset \mathbb{C}$ with center $\lambda_0$ and radius $\varepsilon$ sufficiently small. $P(\lambda_0)$ is a projection and called the Riesz projection for $A$ at $\lambda_0$. If $P(\lambda_0)$ is of finite rank, then $\lambda_0$ has to be an eigenvalue and $\textup{Ran}P(\lambda_0)$ is the algebraic eigensubspace of $A$ at $\lambda_0$; in particular, $m_{al}(\lambda_0)=\textup{Tr}(P(\lambda_0))$.
\begin{lemma}For any $y\in Gr(n,2n)$, there exists a maximal completely negative-definite subspace $N_-$ of $\mathcal{B}_T$ such that $V_y$ and $N_-$ are transversal.
\end{lemma}
\begin{proof}Let $s_y$ be the holomorphic section of $L$ determined by $y\in Gr(n,2n)$. If there were no such $N_-$ claimed in the lemma, then $s_y$ vanishes on $W_-(\mathcal{B}_T)\subset Gr(n, 2n)$ by definition. Since $W_-(\mathcal{B}_T)$ is open in $ Gr(n, 2n)$, the Identity Theorem implies that $s_y\equiv 0$. This is a contradiction.
\end{proof}
As a result of the lemma, if the boundary condition $y\in Gr(n,2n)$ is fixed, then we can always choose a boundary triplet $(G, \Gamma_\pm)$ in term of which the boundary condition $y$ can be parameterized by some $Y\in \mathbb{B}(G)$, i.e.,
\[D(T_y)=\{x\in D(T^*)|\Gamma_-x=Y\Gamma_+x\}.\]
\begin{theorem}\label{thm17}For an entire operator $T$ with deficiency index $n$ and any non-degenerate boundary condition $y\in Gr(n,2n)$, if $\lambda_0$ is an eigenvalue of $T_y$, then $n_{al}(\lambda_0)$ is finite and $n_{al}(\lambda_0)=n_a(\lambda_0)$.
\end{theorem}
\begin{proof}W.l.g, we assume $\lambda_0\in \rho(T_+)$. Due to the above lemma, a suitable boundary triplet can be chosen such that the boundary condition $y$ can be written as $\Gamma_-x=Y\Gamma_+x$ for some $Y\in \mathbb{B}(G)$.

Along the same line of proof of Thm.~\ref{thm8}, we find for $\lambda\in \rho(T_+)\cap \rho(T_y)$
\[\frac{d}{d\lambda}\ln\det(Y-B(\lambda))=-i\textup{Tr}((Y-B(\lambda))^{-1}\gamma_-(\bar{\lambda})^*\gamma_+(\lambda)),\]
where Prop.~\ref{p4} is used.

If $D(\lambda_0,\varepsilon)$ is sufficiently small, we note that by the resolvent formula (\ref{e2}),
\begin{eqnarray*}P(\lambda_0)&=&\frac{i}{2\pi}\oint_{\partial D(\lambda_0,\varepsilon)}(T_Y-\lambda)^{-1}d\lambda\\
&=&\frac{i}{2\pi}\oint_{\partial D(\lambda_0,\varepsilon)}[(T_+-\lambda)^{-1}+i\gamma_+(\lambda)(Y-B(\lambda))^{-1}\gamma_-(\bar{\lambda})^*]d\lambda\\
&=&\frac{-1}{2\pi}\oint_{\partial D(\lambda_0,\varepsilon)}\gamma_+(\lambda)(Y-B(\lambda))^{-1}\gamma_-(\bar{\lambda})^*d\lambda,
\end{eqnarray*}
where the last equality is because $(T_+-\lambda)^{-1}$ is analytic on $\rho(T_+)$. Note that the integral in the last line is a Bochner integral in the Banach algebra of trace class operators. Thus $P(\lambda_0)$ has a finite trace. Consequently, $P(\lambda_0)$ has to be of finite rank and $m_{al}(\lambda_0)<+\infty$. Furthermore,
\begin{eqnarray*}
n_{al}(\lambda_0)&=&\textup{Tr}(P(\lambda_0))=\frac{-1}{2\pi}\oint_{\partial D(\lambda_0,\varepsilon)}\textup{Tr}[\gamma_+(\lambda)(Y-B(\lambda))^{-1}\gamma_-(\bar{\lambda})^*]d\lambda\\
&=&\frac{-1}{2\pi}\oint_{\partial D(\lambda_0,\varepsilon)}\textup{Tr}[(Y-B(\lambda))^{-1}\gamma_-(\bar{\lambda})^*\gamma_+(\lambda)]d\lambda\\
&=&\frac{1}{2\pi i}\oint_{\partial D(\lambda_0,\varepsilon)}\frac{d}{d\lambda}\ln\det(Y-B(\lambda))d\lambda\\
&=& n_a(\lambda_0).
\end{eqnarray*}
Note that the first equality holds because $\textup{Tr}$ is continuous on the Banach algebra of trace class operators in $\mathbb{B}(H)$ and consequently commutes with the Bochner integral. Thus the claim is proved.
\end{proof}
Due to the local nature of the three notions of multiplicity, the results in this section can obviously be extended to more general setting without the assumption that $T$ is entire.
\subsection{On spectral theory of self-adjoint extensions}\label{spe}
Note that the subset of self-adjoint boundary conditions is a (totally) real analytic submanifold of $Gr(n,2n)$, whose $\mu_y$-measure is certainly zero. For a weakly algebraically non-degenerate transcendental entire operator $T$, one may wonder if there is any exceptional self-adjoint boundary condition and whether almost all self-adjoint boundary conditions are not exceptional in a reasonable sense. In this subsection, we show how our viewpoint towards boundary value problems leads to new results on spectral theory of self-adjoint extensions of an entire operator. One of our goals is to generalize Prop.~\ref{p18}.

Our first task is to generalize the interlacing property for a general deficiency index $n\in \mathbb{N}$. This was given in \cite{reiffenstein2020matrix} and our presentation follows closely that paper.
\begin{definition}Let $\Omega\subset \mathbb{C}$ be an open domain in $\mathbb{C}$ and $f$ a meromorphic function on $\Omega$. If $f(z)=\sum_{j=N}^\infty a_j(z-z_0)^j$ is the Laurent expansion of $f$ around $z_0\in \Omega$ where $N\in \mathbb{Z}$ and $a_N\neq 0$, define $\theta_f(z_0)=N$. This way we get a function $\theta_f:\Omega\rightarrow \mathbb{Z}$, called the divisor function of $f$ on $\Omega$.
\end{definition}
The next lemma is clear and gives an alternative characterization of interlacing property.
\begin{lemma}Let $f$ be a meromorphic function on the domain $\Omega\supset \mathbb{R}$. Then the real zeros and poles of $f$ are all simple and interlace if and only if for every finite interval $(a,b)\subset \mathbb{R}$ the inequality
\[|\sum_{x\in(a,b)}\theta_f(x)|\leq 1\]
holds.
\end{lemma}
This motivates the following definition.
\begin{definition}Let $\theta: \mathbb{R}\rightarrow \mathbb{Z}$ have discrete support. $\theta$ is called $n$-interlacing for $n\in \mathbb{N}$ if for any finite interval $(a,b)$ we have
\[|\sum_{x\in(a,b)}\theta(x)|\leq n.\]
\end{definition}
If $f$ is meromorphic on $\Omega\supset \mathbb{R}$, we say $f$ satisfies the $n$-interlacing condition if the zeros and poles of $f$ on $\Omega$ are all real and $\theta_f|_\mathbb{R}$ is $n$-interlacing. The next two lemmas are basic for our purpose. We include Reiffenstein's elegant proofs for the convenience of the reader. Note that our matrix-valued Nevanlinna functions are uniformly strict. This simplifies the proof of the second lemma a bit.
\begin{lemma}\label{lem6}(\cite{reiffenstein2020matrix}) Let $\theta:\mathbb{R}\rightarrow \mathbb{Z}$ have discrete support and $n\in \mathbb{N}$. Then $\theta$ is $n$-interlacing if and only if there are 1-interlacing functions $\theta_1,\cdots, \theta_n$ such that $\theta=\sum_{j=1}^n \theta_j$.
\end{lemma}
\begin{proof}The "if" part is clear. The "only if" part can be proved by constructing these $\theta_j$ explicitly. W.l.g, we can assume that $\theta(0)=0$ and define
\begin{equation*}\Xi(x)=\left \{
\begin{array}{ll}\sum_{t\in (0,x)}\theta(t),\quad\quad x>0\\
-\sum_{t\in (x,0)}\theta(t),\quad x<0\\
\quad 0, \quad \quad \quad \quad \quad \quad x=0.
\end{array}
 \right.
 \end{equation*}
Since $\theta$ has discrete support, $\Xi$ is well-defined. We use $\Xi(x+)$ and $\Xi(x-)$ to denote $\lim_{t\searrow x}\Xi(t)$ and $\lim_{t\nearrow x}\Xi(t)$ respectively. For each $j\in \mathbb{Z}$, we define
\begin{equation*}\theta_j(x)=\left \{
\begin{array}{ll}1,\quad \Xi(x+)>j\geq \Xi(x-),\\
-1,\quad \Xi(x-)> j\geq \Xi(x+),\\
 0, \quad \textup{otherwise}.
\end{array}
 \right.
 \end{equation*}
Each $\theta_j$ is 1-interlacing. Indeed, if $x<y$ and $\theta_j(x)=\theta_j(y)=1$, then by definition
\[\Xi(x-)\leq j<\Xi(x+), \quad \Xi(y-)\leq j<\Xi(y+),\]
and consequently $\Xi(x+)>j \geq \Xi(y-)$. Let $t_0:=\inf \{t>x|\Xi(t+)\leq j\}$. Then $t_0\in (x,y)$. Since $\Xi$ is a step function, it can be derived from this that $\Xi(t_0-)>j\geq \Xi(t_0+)$ and thus $\Xi(t_0)=-1$ by definition. Similarly, if $x<y$ and $\theta_j(x)=\theta_j(y)=-1$, we can find $t_0\in (x,y)$ such that $\theta_j(t_0)=1$. Hence $\theta_j$ is 1-interlacing.

Note that for $x<y$, we have
\[|\Xi(y)-\Xi(x)|=|\sum_{t\in(x,y)}\theta(t)|\leq n\]
for $\theta$ is $n$-interlacing. This implies that all but at most $n$ of these $\theta_j$'s vanish. It can be seen easily that $\theta=\sum_{j=1}^n \theta_j$.
\end{proof}
\begin{lemma}\label{lem7}(\cite{reiffenstein2020matrix}) Let $Q$ be an $n$ by $n$ (uniformly strict) matrix-valued Nevanlinna function meromorphic on $\mathbb{C}$. Then there exist $n$ scalar Nevanlinna functions $q_1,\cdots, q_n$ meromorphic on $\mathbb{C}$ such that $\det Q=q_1\cdot \cdots \cdot q_n$. In particular, $\det Q$ is $n$-interlacing.
\end{lemma}
\begin{proof}
We can prove the first claim by induction. The case of $n=1$ is trivial. Assume the result holds for $n=k$. Now let $Q$ be a $(k+1)\times (k+1)$ matrix-valued Nevanlinna function. Then $\det Q(z)\neq 0$ for any $z\in \mathbb{C}_+$. Denote the submatrix obtained by deleting the $j$-th row and $j$-th column by $Q_{(j)}(z)$. Obviously, $Q_{(j)}(z)$ is a $k\times k$ matrix-valued Nevanlinna function. Then
\[-Q(z)^{-1}=\left(
               \begin{array}{cccc}
                 -\frac{\det Q_{(1)}(z)}{\det Q(z)} & \ast & \ast & \ast \\
                 \ast & -\frac{\det Q_{(2)}(z)}{\det Q(z)} & \ast & \ast \\
                 \ast & \ast & \ddots & \ast \\
                 \ast & \ast & \ast & -\frac{\det Q_{(k+1)}(z)}{\det Q(z)} \\
               \end{array}
             \right).
\]
Note that $-Q(z)^{-1}$ itself is a Nevanlinna function and in particular each $-\frac{\det Q_{(j)}(z)}{\det Q(z)}$ has to be a scalar Nevanlinna function. Set $q(z)=-\frac{\det Q_{(1)}(z)}{\det Q(z)}$. Then
\[\det Q(z)=\det Q_{(1)}(z)\times (-q(z)^{-1})\]
and the result follows from the induction hypothesis.

Obviously, $\theta_{\det Q}=\sum_{j=1}^n\theta_{q_j}$ and the last claim follows from Lemma~\ref{lem6}.
\end{proof}
Now we can prove one of our basic theorems in this subsection.
\begin{theorem}\label{thm10}Let $T$ be a transcendental entire operator with deficiency index $n$ and $y_1, y_2\in Gr(n, 2n)$ two self-adjoint boundary conditions for $T$. Then for any $r>0$, we have
\[|n_T(r,y_1)-n_T(r, y_2)|\leq n,\]
where $n_T(r,y)$ was defined as in the proceeding subsection.
\end{theorem}
\begin{proof}If $y_1$ and $y_2$ are transversal, we can choose a boundary triplet such that $T_{y_1}=T_0$ and $T_{y_2}=T_1$ \cite[Thm.~2.5.9]{behrndt2020boundary}. Let $M(\lambda)$ and $B(\lambda)$ be the corresponding Weyl function and contractive Weyl function respectively. Then eigenvalues (counting multiplicity) of $T_0$ (resp. $T_1$) are precisely the roots of $\det (Id-B(\lambda))=0$ (resp. $\det (Id+B(\lambda))=0$) around the real line. Therefore
\begin{eqnarray*}n_T(r,y_2)-n_T(r,y_1)&=&\sum_{\lambda\in(-r, r)}\theta_{\det(Id+B)}(\lambda)-\sum_{\lambda\in(-r, r)}\theta_{\det(Id-B)}(\lambda)\\
&=&\sum_{\lambda\in(-r, r)}[\theta_{\det(Id+B)}(\lambda)-\theta_{\det(Id-B)}(\lambda)]\\
&=&\sum_{\lambda\in(-r, r)}\theta_{\det M}(\lambda).
\end{eqnarray*}
The last line is due to the fact that around the real line $M=i(Id+B)(Id-B)^{-1}$ and thus $\det M=i^n\det(Id+B)/\det(Id-B)$. The theorem then follows from Lemma \ref{lem7}.

If $y_1$ and $y_2$ are not transversal, we can choose a self-adjoint boundary condition $y_3$ close enough to $y_2$ such that $y_1$ and $y_3$ are transversal.\footnote{The subset of self-adjoint boundary conditions that are transversal to $y_1$ is open and dense in the Lagrangian Grassmannian.} Then
\[|n_T(r,y_1)-n_T(r, y_3)|\leq n\]
and the conclusion follows from the stability of spectrum in a bounded domain.
\end{proof}
\emph{Remark}. It is worthwhile to mention that the theorem holds in more general setting: even if $T$ is not entire, the inequality is still correct in any finite open interval in $\Theta(T)\cap \mathbb{R}$. The above proof still goes through with minor modifications. We also note that the theorem has been already established in \cite[Sec.~3, Chap.~9]{birman2012spectral} by using the Krein resolvent formula. In a more specified case, the theorem was also proved in \cite{hillairet2010finite} by min-max principle.

Thm.~\ref{thm10} demonstrates the uniform patten of the distribution of eigenvalues of all self-adjoint extensions, in sharp contrast with how wild a non-self-adjoint extension may be. The following theorem strengthens this remark.
\begin{theorem}\label{thm16}Let $T$ be a transcendental entire operator with deficiency index $n\in \mathbb{N}$ and $y\in Gr(n,2n)$ a self-adjoint boundary condition for $T$. If a boundary triplet is chosen such that
$T_0=T_y$ and $B(\lambda)$ is the contractive Weyl function, then
\[|\frac{1}{2\pi i}\int_{-r}^rd\ln \det B(s)-n_T(r,y)|\leq n.\]
\end{theorem}
\begin{proof}Consider the map $\Upsilon(s,t)=\det(B(s)-tId)\in \mathbb{C}$ for $(s,t)\in [-r,r]\times [0, \mathcal{T}]$ for sufficiently large $\mathcal{T}>1$. $\Upsilon$ is a homotopy between the two curves $\Upsilon(s,0)=\det B(s)$ and $\Upsilon(s,\mathcal{T})=\det( B(s)-\mathcal{T}Id)$ on $\mathbb{C}$. We also assume that $\det(B(\pm r)-Id)\neq 0$, i.e., $\pm r$ are not eigenvalues of $T_0$ (if not so, we can replace $r$ with a number a bit smaller and use a limiting argument).

Let $C$ be the oriented boundary of the rectangle $R:=[-r,r]\times [0, \mathcal{T}]$ on the $st$-plane. The four pieces are denoted by $C_j$, $j=1,2,3,4$ (see Fig. 1).
~\\
 \begin{center}
\begin{picture}(160, 80)
\vector(1,0){150}
\put(-80,0){\vector(1,0){10}}
\put(-80,0){\vector(0,1){90}}
\put(-40,0){\vector(0,1){70}}
\put(-120,70){\vector(0,-1){70}}
\put(-40,70){\vector(-1,0){80}}
\put(5,0){$s$}
\put(-78,85 ){$t$}
\put(-80,5){$C_1$}
\put(-40,30){$C_2$}
\put(-135,30){$C_4$}
\put(-80,60){$C_3$}
\put(-70,30){$D$}
\put(-85, -15){\textup{Fig.~1}}
\end{picture}
\end{center}

 Now consider the integral
\begin{equation}\frac{1}{2\pi i}\oint_C\Upsilon^*(\frac{dz}{z})=\frac{1}{2\pi }\oint_C\Im \Upsilon^*(\frac{dz}{z})=\frac{1}{2\pi }\sum_{j=1}^4\int_{C_j}\Im\Upsilon^*(\frac{dz}{z}).\label{e5}\end{equation}
Note that there is no singularity on $\Upsilon(C)$ for $dz/z$ and the only singularities in $\Upsilon(R)$ correspond precisely to roots of $\det(B(s)-Id)=0$ in $(-r,r)$. By the Argument Principle, we know that
\begin{equation}\frac{1}{2\pi }\oint_C\Im \Upsilon^*(\frac{dz}{z})=n_T(r,y).\label{e7}\end{equation}
Indeed,
\[\det(B(s)-t)=\prod_{j=1}^n(\varrho_j(s)-t)\]
 where $\varrho_j(s)$ are eigenvalues of $B(s)$. If $s=\lambda_0\in (-r,r)$ is an eigenvalue of $T_0$ with geometric multiplicity $k$, then due to Thm.~\ref{thm8}, as $s\rightarrow \lambda_0$ and $t\rightarrow 1$,
\[\det(B(s)-t)=q(s,t)\prod_{j=1}^k(a_j(s-\lambda_0)-t+1+o(s-\lambda_0))\]
where $q(s,t)$ is smooth and non-vanishing, $a_j$ are complex numbers with $\Im a_j\neq 0$ and $o(s-\lambda_0)$ is a term of order higher than $s-\lambda_0$. This local form means the contribution of $\lambda_0$ to the left side of Eq.~(\ref{e7}) is $k$. Therefore, contributions from all eigenvalues of $T_0$ in $(-r,r)$ is $n_T(r,y)$.

Let's estimate each term on the right hand side of Eq.~(\ref{e5}). We find that
\[\frac{1}{2\pi }\int_{C_1}\Im \Upsilon^*(\frac{dz}{z})=\frac{1}{2\pi i}\int_{-r}^r d\ln \det B(s).\]
Note that for $\mathcal{T}$ large enough, we have
\[|\det (B(s)-\mathcal{T})-(-1)^n\mathcal{T}^n|< c\mathcal{T}^{n-1}\] where $c$ is a positive constant independent of $s\in \mathbb{R}$. In particualr, $|\Im \det (B(s)-\mathcal{T})|<c\mathcal{T}^{n-1}$. Furthermore, $\mathcal{T}$ can be chosen so that either $\Re \det (B(s)-\mathcal{T})>\mathcal{T}^n/2$ (if $n$ is even, see Fig.~2.) or $<-\mathcal{T}^n/2$ (if $n$ is odd).
\begin{center}
\begin{picture}(80,100)
\put(-40, 15){\vector(1,0){120}}
\put(-30,-25){\line(0,1){80}}
\put(-10,-25){\line(0,1){80}}
\put(-30,13){\vector(3,4){40}}
\put(-30,13){\vector(2,-3){40}}
\bezier{40}(10,60)(39,12)(13,-40)
\put(-25,68){$\Re z=\mathcal{T}^n/2$}
\put(10,55){$\Upsilon(r, \mathcal{T})$}
\put(13,-45){$\Upsilon(-r, \mathcal{T})$}
\put(9,65){\circle*{3}}
\put(10,-45){\circle*{3}}
\put(27,15){$\Upsilon(\cdot, \mathcal{T})$}
\put(13,-65){\textup{Fig.~2}}
\end{picture}
\end{center}
~\\~\\

In either case, we can choose a single-valued branch of $\ln z$ to find
\[|\frac{1}{2\pi }\int_{C_3}\Im \Upsilon^*(\frac{dz}{z})]|<\frac{1}{\pi}\arctan\frac{2c}{\mathcal{T}}.\]

Let us turn to
\[\int_{C_2}\Im \Upsilon^*(\frac{dz}{z})=\int_0^\mathcal{T} \Im d\ln \det (B(r)-t)=\Sigma_{j=1}^n\int_0^\mathcal{T}\Im d\ln (\rho_j(r)-t),\]
where $\rho_j(r)\in \mathbb{U}(1)$, $j=1,\cdots, n$, are the eigenvalues of $B(r)$. Again, we can choose a single-valued branch of $\ln z$ such that (see Fig.~3)
\[|\frac{1}{2\pi }\int_{C_2}\Im d\ln (\rho_j(r)-t)]|<\frac{1}{2}.\]

 \begin{center}
\begin{picture}(80, 60)
\put(15,10){O}
\put(15,15){\circle{100}}
\put(15,15){\line(1,0){100}}
\put(20,35){\circle*{3}}
\put(7,40){$\rho_j(r)$}
\put(15,14){\vector(1,3){7}}
\put(100,5){$\mathcal{T}$}
\put(102,15){\circle*{3}}
\put(100,15){\vector(-4,1){80}}
\put(15,-20){\textup{Fig.~3}}
\end{picture}
\end{center}
~\\

Consequently,
\[|\frac{1}{2\pi }\int_{C_2}\Im \Upsilon^*(\frac{dz}{z})|<\frac{n}{2}.\]
Similarly,
\[|\frac{1}{2\pi}\int_{C_4}\Upsilon^*(\frac{dz}{z})|<\frac{n}{2}.\]

Combining all the estimates together and letting $\mathcal{T}$ go to $+\infty$ then completes the proof.
\end{proof}
\emph{Remark}. If $\tilde{y}$ is the self-adjoint boundary condition parameterized by $-Id$, then along the same way we have
\[|\frac{1}{2\pi i}\int_{-r}^rd\ln \det B(s)-n_T(r,\tilde{y})|\leq n.\]
Consequently,
\[|n_T(r,y)-n_T(r, \tilde{y})|\leq 2n.\]
This is a weaker version of Thm.~\ref{thm10}.
\begin{corollary}Let assumptions be as in Thm.~\ref{thm16}. If $T_0$ has no eigenvalues in the interval $(a,b)$ (one of the end could be infinity), then
\[0<\frac{1}{2\pi i}\int_a^bd\ln \det B(s)\leq n.\]
\end{corollary}
\begin{proof}If $(a,b)$ is a finite interval, the proof of Thm.~\ref{thm16} carries over here. If $(a,b)$ is infinite, a limiting argument will suffice. Note that the integration is positive, because the integrand is so (see Prop.~\ref{p4}).
\end{proof}

\begin{corollary}Let $T$ be a transcendental entire operator with deficiency index $n\in \mathbb{N}$ and $0<\rho_T<+\infty$. Then for any self-adjoint boundary condition $y\in Gr(n,2n)$, we have $\alpha_T(y)=\rho_T$.
\end{corollary}
\begin{proof}Generally for a function $\alpha(r)$ positive and non-decreasing for sufficiently large $r$, the order of $\alpha(r)$ is defined to be $\limsup_{r\rightarrow +\infty}\frac{\ln\alpha(r)}{\ln r}$ \cite[\S~1, Chap.~2]{gol'dberg2008value}. By \cite[Thm.~1.8, Chap.~2]{gol'dberg2008value}, $\alpha_T(y)$ is actually the order of $n_T(r,y)$. W.l.g, we can assume a boundary triplet is chosen such that $T_0=T_y$. Then the above theorem and \cite[Thm.~1.2, Chap.~2]{gol'dberg2008value} imply $n_T(r,y)$ has the same order as $\frac{1}{2\pi i}\int_{-r}^rd\ln \det B(s)$. The conclusion follows from \cite[Thm.~1.1, Chap.~2]{gol'dberg2008value}.
\end{proof}

The following easy result may explain why spectral theory of self-adjoint extensions are well-developed, compared with that of non-self-adjoint ones: the distribution of eigenvalues for any self-adjoint extension is well-reflected by $T$ itself.
\begin{corollary}If $T$ is a transcendental entire operator with deficiency index $n\in \mathbb{N}$, then there is no Varilon exceptional self-adjoint boundary condition for $T$.
\end{corollary}
\begin{proof}By definition, this is clear from the above theorem and the First Main Theorem.
\end{proof}
\subsection{Unicity theorems and inverse spectral theory }
In value distribution theory for meromorphic functions, two meromorphic functions $f,g$ are said to share the value $a\in \mathbb{CP}^1$ if $f^{-1}(\{a\})=g^{-1}(\{a\})$. Using his now famous Second Main Theorem, Nevanlinna discovered the following Five-Value-Theorem \cite[Coro.~A5.2.5]{ru2021nevanlinna}.
\begin{proposition}If two non-constant meromorphic functions on $\mathbb{C}$ share five distinct values, then $f\equiv g$.
\end{proposition}
This result can be improved by taking account of analytic multiplicity of each $\lambda$ such that $f(\lambda)=g(\lambda)=a$. Such theorems are called Unicity Theorems in value distribution theory and much progress has been made for the case of meromorphic functions. However, for general entire curves in projective algebraic manifolds, much less is known up to now. We would like to point out the meaning of Unicity Theorems in inverse spectral theory. We view all the Weyl curves of entire operators with deficiency index $n$ to lie in the same Grassmannian $Gr(n,2n)$.
\begin{definition}Two entire operators $T$, $T'$ with deficiency index $n$ are said to share the same boundary condition $y\in Gr(n,2n)$ if $W_{T}^{-1}(Z_y)=W_{T'}^{-1}(Z_y)$ where $Z_y=\{x\in Gr(n,2n)|s_y(x)=0\}$.
\end{definition}
Note that $W_T^{-1}(Z_y)$ is precisely the spectrum (without counting multiplicities) of the extension $T_y$. Thus entire operators $T$ and $T'$ share the boundary condition $y$ if and only if the two operators $T_y$ and $T'_y$ have the same spectrum. If $n=1$, that $T$ and $T'$ share the boundary condition $y$ means precisely the Weyl curves share the value $y\in \mathbb{CP}^1$. In inverse spectral theory, one is basically concerned with determining a symmetric operator from the spectra of its different (self-adjoint) extensions. A classical example was suggested by G. Borg in \cite{Borg1946umkehrung}, see also our \S~\ref{real}. Another example that can be put into this framework is M. Kac's famous problem "can one hear the shape of a Drum?" \cite{kac1966can}.

The Five-Value-Theorem has a beautiful counterpart in spectral theory.
\begin{proposition}\label{p15}For two entire operators $T$ and $T'$ with deficiency index 1, if their Weyl curves share two self-adjoint boundary conditions, then $T$ and $T'$ are unitarily equivalent.
\end{proposition}
\begin{proof}We can choose boundary triplets such that the two self-adjoint boundary conditions are parameterized by $\pm 1$. Let $M(\lambda)$ and $\tilde{M}(\lambda)$ be the corresponding Weyl functions. Note that $M(\lambda)$ and $\tilde{M}(\lambda)$ are meromorphic function on $\mathbb{C}$, sharing the common simple real zeros and poles. $M(\lambda)$ (or $\tilde{M}(\lambda)$) is determined by its zeros and poles up to a positive constant factor \cite[Coro.~3.1.5]{reiffenstein2020matrix}. Then $\tilde{M}(\lambda)=cM(\lambda)$ for a positive constant $c>0$ and consequently $T$ and $T'$ have the same curvature. Due to Thm.~\ref{thm14}, $T$ and $T'$ are unitarily equivalent. (It is also easy to see $\tilde{M}(\lambda)$ is a M$\ddot{o}$bius transform of $M(\lambda)$.)
\end{proof}
\emph{Remark}. More information is determined by the two spectra. By replacing the boundary triplet $(\mathbb{C}, \Gamma_0, \Gamma_1)$ of $T$ with $(\mathbb{C}, c^{-1/2}\Gamma_0, c^{1/2}\Gamma_1)$ (this won't change $T_0$ and $T_1$), we can assume that $\tilde{M}(\lambda)=M(\lambda)$. Therefore, $T_0$ (resp. $T_1$) is unitarily equivalent to $T'_0$ (resp. $T'_1$). We will show in \S~\ref{sturm} that Borg's result can be viewed as an enhanced version of Prop.~\ref{p15} in the setting of Sturm-Liouville operators.

We would like to know for $n>1$ whether we can prove $T$ and $T'$ are unitarily equivalent if they share sufficiently many self-adjoint boundary conditions.
\section{Real symmetric operators}\label{real}
According to \cite{dieudonne1983history}, "after 1913, almost all papers on spectral theory in Hilbert space dealt exclusively with complex Hilbert space". However, the extension theory of \emph{real} symmetric operators in \emph{real} Hilbert spaces does have some different flavours and the underlying real structure really affects the spectral theoretic content.

If $H$ is a real Hilbert space, a real strong symplectic structure $\omega$ on $H$ can be defined in a manner similar to the complex case, simply by replacing "sesquilinear" with "bilinear" and the condition ii) with $$\omega(x,y)=-\omega(y, x)$$ for any $x,y\in H$. In this case, $\textup{dim}H$ should be even and we assume it to be $2n$. Then $n$ is a complete symplectic invariant of $\omega$. That means all real strong symplectic Hilbert spaces of dimension $2n$ are isomorphic (If $n$ is finite, this is a standard result in symplectic geometry. We leave the case of $n=+\infty$ to the interested reader.).
\begin{example}If $G$ is a real Hilbert space and $\mathbb{G}=G\oplus_\bot G$, then on $\mathbb{G}$ there is a standard strong symplectic structure: if $x=(x_1,x_2), y=(y_1,y_2)\in \mathbb{G}$, then
\[\omega(x,y)=(x_2,y_1)_G-(x_1,y_2)_G.\]
\end{example}

We can complexify $H$ to get a complex Hilbert space $H_\mathbb{C}:=H\otimes_\mathbb{R}\mathbb{C}$. For $x=x_1+ix_2\in H_\mathbb{C}$ where $x_1,x_2\in H$, we denote $x_1-ix_2$ by $\bar{x}$ or $\mathfrak{j}x$. We call $x_1$ (resp. $x_2$) the real (resp. imaginary) part of $x$. $\mathfrak{j}$ is usually called a real structure on $H_\mathbb{C}$ and $H$ is recovered by taking the $\mathfrak{j}$-invariant part of $H_\mathbb{C}$. Note that for $x_1+ix_2, y_1+iy_2\in H_\mathbb{C}$ ($x_i, y_i\in H$, $i=1,2$), the inner product on $H_\mathbb{C}$ is
\[(x_1+ix_2, y_1+iy_2)_{H_\mathbb{C}}=(x_1,y_1)+(x_2,y_2)+i(x_2,y_1)-i(x_1,y_2).\]
In particular, $\|x_1+ix_2\|_{H_\mathbb{C}}=\|x_1-ix_2\|_{H_\mathbb{C}}$. If $M\subset H_\mathbb{C}$ is a subspace, denote the subspace $\{x\in H_\mathbb{C}| \bar{x}\in M\}$ by $\bar{M}$.

$\omega$ can also be complexified in a bilinear (not sesquilinear) way to be defined on $H_\mathbb{C}$. Denote this extension by $\omega_\mathbb{C}$.
\begin{definition}A maximal isotropic complex subspace of $H_\mathbb{C}$ w.r.t. $\omega_\mathbb{C}$ is called Lagrangian (its dimension has to be $n$). A complex Lagrangian subspace $L$ is transversal if we have $H_\mathbb{C}=L\oplus \bar{L}$.
 \end{definition}Denote the space of all complex Lagrangian subspaces by $\mathfrak{L}(H_\mathbb{C})$. It is called the Lagrangian Grassmanian of $(H_\mathbb{C}, \omega_\mathbb{C})$, which is a smooth subvariety of $Gr(n, 2n)$.
\begin{definition}
$L\in \mathfrak{L}(H_\mathbb{C})$ is called completely positive-definite if there is a constant $c_L$ such that
\[-i\omega_\mathbb{C}(x, \bar{x})\geq c_L \|x\|^2_{H_\mathbb{C}}\]
for all $x\in L$.
\end{definition}
We denote the space of all completely positive-definite complex Lagrangian subspaces by $\mathcal{L}_+(H_\mathbb{C})$. Similarly, we can define completely negative-definite complex Lagrangian subspaces and the corresponding space $\mathcal{L}_-(H_\mathbb{C})$. Obviously $\mathcal{L}_\pm(H_\mathbb{C})\subset \mathfrak{L}(H_\mathbb{C})$. These subspaces do exist: as before, by Riesz representation theorem, $\omega(x,y)=(Ax,y)_H$ for a linear isomorphism $A$ such that $A^*=-A$. Then $\mathcal{J}:=|A|^{-1}A$ satisfies $\mathcal{J}^2=-Id$. $\mathcal{J}$ can be complexified to be a complex linear map on $H_\mathbb{C}$. Let $M:=\textup{ker}(\mathcal{J}-i)$. Then $M\in \mathcal{L}_+(H_\mathbb{C})$ and $\bar{M}=\textup{ker}(\mathcal{J}+i)\in \mathcal{L}_-(H_\mathbb{C})$. Note that $\omega_\mathbb{C}$ can be used to define a strong symplectic structure on $H_\mathbb{C}$ in the following manner:
\[[x,y]=\omega_\mathbb{C}(x,\bar{y})\]
for any $x, y\in H_\mathbb{C}$.
\begin{lemma}If $L\in \mathcal{L}_+(H_\mathbb{C})$, then $\bar{L}\in \mathcal{L}_-(H_\mathbb{C})$, and $L$ is transversal.
\end{lemma}
\begin{proof}If $x\in \bar{L}$, then
\[i\omega_\mathbb{C}(x,\bar{x})=-i\omega_\mathbb{C}(\bar{x}, x)\geq c_L \|\bar{x}\|^2_{H_\mathbb{C}},\]
implying $L\in \mathcal{L}_-(H_\mathbb{C})$.

If $0\neq x\in L\cap \bar{L}$, then also $\bar{x}\in L\cap \bar{L}$ and consequently both the real and imaginary parts of $x$ are in $L$. W.l.g, assume $x\in H$. Then $-i\omega_{\mathbb{C}}(x,\bar{x})=-i\omega_{\mathbb{C}}(x,x)=0$, contradicting the positive-definiteness of $L$. So $L\cap \bar{L}=0$  (If $n\in \mathbb{N}$, the proof is finished).

Note that w.r.t. the strong symplectic structure $[\cdot, \cdot]$, $L$ is a maximal completely positive-definite subspace, and $L^{\bot_s}=\bar{L}$. Thus $H_\mathbb{C}=L\oplus \bar{L}$.
\end{proof}
The theorem shows that $\mathcal{L}_\pm(H_\mathbb{C})\subset W_\pm(H_\mathbb{C})$. Thus if $M\in \mathcal{L}_+(H_\mathbb{C})$ is fixed, any $L\in \mathcal{L}_+(H_\mathbb{C})$ shall be parameterized by an operator $B\in \mathbb{B}(M,\bar{M})$ with $\|B\|< 1$. However, $B$ should satisfy an additional condition.
\begin{theorem}\label{thm15}If $N\in \mathcal{L}_+(H_\mathbb{C})$ is fixed, then $L\in W_+(H_\mathbb{C})$ is in $\mathcal{L}_+(H_\mathbb{C})$ if and only if $L=\{(x, Bx)\in N\oplus \bar{N}| x\in N\}$, where
 $B\in \mathbb{B}(N,\bar{N})$, $\|B\|< 1$, and $B^*\mathfrak{j}=\mathfrak{j}B$.
\end{theorem}
\begin{proof}We only need to prove the statement relevant to the formula $B^*\mathfrak{j}=\mathfrak{j}B$. If $L\in W_+(H_\mathbb{C})$ is parameterized by $B\in \mathbb{B}(N,\bar{N})$, then $L\in \mathcal{L}_+(H_\mathbb{C})$ if and only if for any $x,y\in N$, we shall have
\[\omega_\mathbb{C}(x+Bx,y+By)=\omega_\mathbb{C}(x,By)+\omega_\mathbb{C}(Bx,y)=0.\]
Note that $\omega_\mathbb{C}(u,v)=[u,\bar{v}]$ for all $u,v\in H_\mathbb{C}$. Thus the latter condition is equivalent to
\[[x,\mathfrak{j}By]+[Bx,\mathfrak{j}y]=i(x,\mathfrak{j}By)_+-i(Bx,\mathfrak{j}y)_-=0,\]
which is further equivalent to $B^*\mathfrak{j}=\mathfrak{j}B$.
\end{proof}
If $n$ is finite and an orthonormal basis $\{e_j\}_{j=1}^n$ of $N\in \mathcal{L}_+(H_\mathbb{C})$ is chosen, then $\{\mathfrak{j}e_j\}_{j=1}^n$ is an orthonormal basis of $\bar{N}$. In terms of these, the condition $B^*\mathfrak{j}=\mathfrak{j}B$ is precisely $B=B^t$, where $B^t$ is the transpose of $B$. That's to say, $B$ is a symmetric matrix. If $M=i(Id+B)(Id-B)^{-1}$, then $M$ is a symmetric complex matrix with a positive-definite imaginary part $(M-\bar{M})/2i$. If $\mathcal{L}_+(H_\mathbb{C})$ is identified with the space of such matrices $M$, it is the famous Siegel upper half space of genus $n$. Note that $\textup{dim}W_+(H_\mathbb{C})=n^2$ while $\textup{dim}\mathcal{L}_+(H_\mathbb{C})=\frac{n(n+1)}{2}$. $\mathcal{L}_+(H_\mathbb{C})$ is an irreducible non-compact Hermitian symmetric space of type $III_n$ and $\mathfrak{L}(H_\mathbb{C})$ is its compact dual.

In the above context, we can consider $\mathcal{L}(H_\mathbb{C})$ as well, i.e., the space of Lagrangian subspaces w.r.t. the strong symplectic structure $[\cdot, \cdot]$. Note that $\mathcal{L}(H_\mathbb{C})$ is different from $\mathfrak{L}(H_\mathbb{C})$.
\begin{definition}We say $L\in \mathcal{L}(H_\mathbb{C})$
is real, if $L=\bar{L}$.
\end{definition}
It is necessary that if $L$ is real, then $L\in \mathfrak{L}(H_\mathbb{C})$ and $L$ has to be the complexification of a Lagrangian subspace of $H$ w.r.t. $\omega$.
\begin{proposition}\label{p16}If $N$ as above is fixed, then $L\in \mathcal{L}(H_\mathbb{C})$ is real if and only if it is parameterized by $U\in \mathbb{U}(N,\bar{N})$ satisfying $\mathfrak{j}U^{-1}=U\mathfrak{j}$.
\end{proposition}
\begin{proof}If $L$ is of the form $\{x+Ux\in H_\mathbb{C}|x\in N\}$ for $U\in\mathbb{U}(N,\bar{N})$, then $L=\bar{L}$ if and only if $\mathfrak{j}x=U\mathfrak{j}Ux$ for any $x\in N$. The claim then follows.
\end{proof}
Denote the subset of real elements in $\mathfrak{L}(H_\mathbb{C})$ by $\mathfrak{L}_r(H_\mathbb{C})$.
\begin{proposition}\label{p17}Let $\overline{\mathcal{L}_\pm(H_\mathbb{C})}$ be the closure of $\mathcal{L}_\pm(H_\mathbb{C})$ in $\mathfrak{L}(H_\mathbb{C})$. Then
\[\overline{\mathcal{L}_+(H_\mathbb{C})}\cap\overline{\mathcal{L}_-(H_\mathbb{C})}=\mathfrak{L}_r(H_\mathbb{C}).\]
\end{proposition}
\begin{proof}By Thm.~\ref{thm15}, Prop.~\ref{p16} and Prop.\ref{p11}, this is clear.
\end{proof}

A real symmetric operator $T$ is defined in the same manner as in the complex case. Actually except for statements on spectra, the basic definitions take the same form as in the complex case and we won't bother to spell out all the details.

Let $T$ be a real symmetric operator defined in $H$. Then $\mathcal{B}_T:=D(T^*)/D(T)$ is again a real Hilbert space and it is equipped with the strong symplectic structure
\[\omega_T([x],[y])=(T^*x,y)-(x,T^*y).\]
The symplectic structure $\omega_T$ on $\mathcal{B}_T$ can be viewed as the reduced version of the standard strong symplectic structure on $H\oplus_\bot H$. Let $2n$ be the dimension of $\mathcal{B}_T$. We call $n$ the deficiency index of $T$. $T$ can be complexified in the obvious way to be a symmetric operator $\mathrm{T}$ on $H_\mathbb{C}$ with deficiency indices $(n,n)$. We say $T$ is simple, if $\mathrm{T}$ is so\footnote{It can be seen that $T$ is simple if and only if $T$ cannot be written as a nontrivial orthogonal direct sum of a real self-adjoint operator and another real symmetric operator.}. We only consider simple real symmetric operators.

The space $\mathcal{B}_\mathrm{T}$ can be naturally identified with $\mathcal{B}_T\otimes \mathbb{C}$. As before, for $\lambda\in \mathbb{C}_+\cup \mathbb{C}_-$, $W_T(\lambda)$ is defined to be the image of $\textup{ker}(\mathrm{T}-\lambda)$ in $\mathcal{B}_\mathrm{T}$ (and thus in $\mathcal{B}_T\otimes \mathbb{C}$). We also say $W_T(\lambda)$ is the Weyl map (curve) of $T$. Since $\textup{ker}(\mathrm{T}-\bar{\lambda})=\mathfrak{j}(\textup{ker}(\mathrm{T}-\lambda))$, we have $W_T(\bar{\lambda})=\mathfrak{j}W_T(\lambda)$, where we have used the same $\mathfrak{j}$ to denote the real structure on $\mathcal{B}_\mathrm{T}$.
\begin{theorem}For any $\lambda\in \mathbb{C}_+$, $W_T(\lambda)\in \mathcal{L}_+(\mathcal{B}_T\otimes \mathbb{C})$ and $\mathcal{B}_T\otimes \mathbb{C}=W_T(\lambda)\oplus W_T(\bar{\lambda})$.
\end{theorem}
\begin{proof}For $\lambda\in \mathbb{C}_+$, since by Prop.~\ref{p6} $W_T(\lambda)\in W_+(\mathcal{B}_T\otimes \mathbb{C})$, we only need to prove it is isotropic. Note that for any $x,y\in \textup{ker}(\mathrm{T}-\lambda)$,
\begin{eqnarray*}\omega_{T\mathbb{C}}([x],[y])=(\mathrm{T}^*x,\mathfrak{j}y)_{H_\mathbb{C}}-(x,\mathrm{T}^*\mathfrak{j}y)_{H_\mathbb{C}}
=\lambda(x,\mathfrak{j}y)_{H_\mathbb{C}}-(x,\bar{\lambda}\mathfrak{j}y)_{H_\mathbb{C}}=0.\end{eqnarray*}
The result follows.
\end{proof}
This theorem shows that the single-branched Weyl curve of $\mathrm{T}$ lives in the Siegel upper half space of genus $n$, which is an analytic submanifold of dimension $n(n+1)/2$ in the $n^2$-dimensional complex manifold $W_+(H_\mathbb{C})$. If furthermore $\mathrm{T}$ is entire, then $W_T(\lambda)$ should be an entire curve in $\mathfrak{L}(\mathcal{B}_T\otimes \mathbb{C})\subset Gr(n,2n)$. In this sense, $\mathrm{T}$ as a symmetric operator in $H_\mathbb{C}$ is always algebraically degenerate. This speciality of $\mathrm{T}$ is never noted in the literature, to the best of our knowledge.
\begin{definition}
A simple real symmetric operator $T$ in a real Hilbert space $H$ is called entire if its complexification $\mathrm{T}$ is entire. An entire real symmetric operator $T$ is algebraically non-degenerate if its Weyl curve $W_T(\lambda)$ is Zarisky dense in $\mathfrak{L}(\mathcal{B}_T\otimes \mathbb{C})$. $T$ is called weakly algebraically non-degenerate if each $y\in \mathfrak{L}(\mathcal{B}_T\otimes \mathbb{C})$ as a boundary condition for $\mathrm{T}$ is non-degenerate.
\end{definition}
Obviously, $\mathfrak{L}(\mathcal{B}_T\otimes \mathbb{C})$ should be the suitable ambient space when one is to investigate the spectral theory of an algebraically non-degenerate real entire operator from the viewpoint of value distribution theory. We won't pursue this aspect in this paper any further.

If $\Phi$ is a symplectic isomorphism between $(\mathcal{B}_T, \omega_T(\cdot,\cdot))$ and the standard real strong symplectic Hilbert space $G\oplus_\bot G$ for a real Hilbert space $G$ of dimension $n$, there is an induced map $(\Gamma_0, \Gamma_1):D(T^*)\rightarrow G\oplus_\bot G$ such that for all $x,y\in D(T^*)$
\[(T^*x,y)-(x,T^*y)=(\Gamma_1x,\Gamma_0y)-(\Gamma_0x,\Gamma_1y).\]
It can be checked easily that $(G,\Gamma_0,\Gamma_1)$ can be complexified to produce a boundary triplet $(G\otimes \mathbb{C}, \Gamma_0, \Gamma_1)$ for $\mathrm{T}$. In this sense, we call $(G\otimes \mathbb{C}, \Gamma_0, \Gamma_1)$ a real boundary triplet for $\mathrm{T}$.

\section{Applications to Sturm-Liouville problems}\label{sturm}

Sturm-Liouville problems have a long history of nearly 200 years, and in early 1930s M. Stone first treated this topic in the formalism of unbounded symmetric operators. However, in this section we only consider the special case $L=-\frac{d^2}{dx^2}+q(x)$ for a real-valued smooth ($C^\infty$) potential function $q(x)$ on the interval $[0,\pi]$. The regularity of $q(x)$ can be relaxed, but our basic concern is the influence of different boundary conditions on the spectral theoretic content. The purpose of this section is two-fold: on one side we shall put this classical topic in the new light of the paper and see something new our theory can convey. On the other side, we use this material to demonstrate some concepts introduced in previous sections. Our basic references are \cite[Chap.~6]{behrndt2020boundary} \cite[Chap.~15]{schmudgen2012unbounded} \cite{levitan2018inverse} and \cite{freiling2001inverse}.

By $L^2(0, \pi)$ we mean the space of complex-valued square integrable functions on $[0,\pi]$. For $L=-\frac{d^2}{dx^2}+q(x)$ on $[0, \pi]$, we can define a symmetric operator $T$ by choosing its domain $D(T)$ to be
\[\{f\in L^2(0,\pi)|f, f'\in AC[0,\pi],\, f(0)=f(\pi)=f'(0)=f'(\pi)=0,\, Lf\in L^2(0,\pi)\},\]
where $AC[0,\pi]$ denotes the space of absolutely continuous functions on $[0,\pi]$. $D(T)$ has $C_0^\infty(0,\pi)$ (the space of smooth functions with compact support in $(0,\pi)$) as a subspace and actually $T$ is the closure of $T|_{C_0^\infty(0,\pi)}$ in $L^2(0,\pi)$. Then it's a basic result that
\[D(T^*)=\{f\in L^2(0,\pi)|f, f'\in AC[0,\pi],\, Lf\in L^2(0,\pi)\}.\]
$T$ is simple \cite[Coro.~6.3.5]{behrndt2020boundary} and has deficiency indices $(2,2)$. It should be pointed out that $T$ is in fact the complexification of a real symmetric operator (just obtained by replacing all the above relevant spaces with their real counterparts). In the following, if we want to emphasize the underlying $q$, we shall add a subscript $q$ to the corresponding operators, e.g., $T_q$.

\subsection{Various boundary conditions}

The following boundary triplet is natural (and real in the sense of \S~\ref{real}):
\[G=\mathbb{C}^2, \quad \Gamma_0 y=\left(
                               \begin{array}{c}
                                 y(0) \\
                                 y(\pi) \\
                               \end{array}
                             \right), \quad
                             \Gamma_1 y=\left(
                               \begin{array}{c}
                                 y'(0) \\
                                 -y'(\pi) \\
                               \end{array}
                             \right).
\]
Then all self-adjoint boundary conditions can be written in the following form:
\[\Gamma_1y-i\Gamma_0y=U(\Gamma_1y+i\Gamma_0y),\]
where $U\in \mathbb{U}(2)$. In this scheme, the matrices corresponding to the traditional so-called separated self-adjoint boundary conditions take the form $\left(
 \begin{array}{cc}
 e^{i\alpha} & 0 \\
 0 & e^{i\beta} \\
\end{array}
\right)
$ for $\alpha, \beta\in \mathbb{R}$. It can be checked easily that real self-adjoint extensions are parameterized by $U\in \mathbb{U}(2)$ such that $U=U^t$; in particular, $T_0$ is the Dirichlet extension and $T_1$ the Neumann extension.

Let $s(x, \lambda)$ (resp. $c(x,\lambda)$) be the unique solution of the equation $$T^*y:=-y''+q(x)y=\lambda y$$ with the initial value condition
$\left\{
\begin{array}{ll}
y(0)=0,\\
y'(0)=1.
\end{array}
\right.$ (resp. $\left\{
\begin{array}{ll}
y(0)=1,\\
y'(0)=0.
\end{array}
\right.$). In terms of these solutions, it can be checked that the Weyl function associated to $(G,\Gamma_0,\Gamma_1)$ is
\begin{equation}M(\lambda)=-\frac{1}{s(\pi,\lambda)}\left(
               \begin{array}{cc}
                 c(\pi, \lambda) & -1 \\
                 -1 & s'(\pi,\lambda) \\
               \end{array}
             \right),
\label{eq9}\end{equation}
and the contractive Weyl function is
\[B(\lambda)=\frac{1}{(\mathfrak{c}-i\mathfrak{s})(\mathfrak{s}'-i\mathfrak{s})-1}\left(
                                                                                              \begin{array}{cc}
                                                                                                (\mathfrak{c}+i\mathfrak{s})(\mathfrak{s}'-i\mathfrak{s})-1 & i2\mathfrak{s} \\
                                                                                                i2\mathfrak{s} & (\mathfrak{c}-i\mathfrak{s})(\mathfrak{s}'+i\mathfrak{s})-1 \\
                                                                                              \end{array}
                                                                                            \right),
\] where $\mathfrak{c}=c(\pi, \lambda)$, $\mathfrak{s}=s(\pi, \lambda)$ and $\mathfrak{s}'=s'(\pi, \lambda)$. Note that both $M(\lambda)$ and $B(\lambda)$ are symmetric complex matrices, just as expected from the discussion in \S~\ref{real}.
\begin{proposition}The symmetric operator $T$ is entire and of Weyl order $1/2$.
\end{proposition}
\begin{proof}Since $M(\lambda)$ is analytic away from $\sigma(T_0)$ and has points in $\sigma(T_0)$ as its simple poles, to prove $T$ is entire, we only have to check that $B(\lambda)$ is analytic at each point in $\sigma(T_0)$. Note that there is an identity $\mathfrak{c}\mathfrak{s}'-\mathfrak{c}'\mathfrak{s}\equiv 1$ due to the Liouville formula $c(x,\lambda)s'(x,\lambda)-c'(x,\lambda)s(x,\lambda)\equiv 1$. Then
\[(\mathfrak{c}-i\mathfrak{s})(\mathfrak{s}'-i\mathfrak{s})-1=\mathfrak{s}(\mathfrak{c}'-\mathfrak{s}-i\mathfrak{c}-i\mathfrak{s}').\]
For $\lambda\in \mathbb{R}$, $\mathfrak{c}'-\mathfrak{s}-i\mathfrak{c}-i\mathfrak{s}'$ cannot be zero. Otherwise, we must have $\mathfrak{c}'=\mathfrak{s}$ and $\mathfrak{c}=-\mathfrak{s}'$ for $\lambda\in \mathbb{R}$ and consequently
\[1=\mathfrak{c}\mathfrak{s}'-\mathfrak{c}'\mathfrak{s}=-\mathfrak{c}^2-(\mathfrak{c}')^2,\]
a contradiction! That's to say, real zeros of $(\mathfrak{c}-i\mathfrak{s})(\mathfrak{s}'-i\mathfrak{s})-1$ are precisely those of $\mathfrak{s}$. By using the same identity we can show
\[(\mathfrak{c}\pm i\mathfrak{s})(\mathfrak{s}'\mp i\mathfrak{s})-1\]
also have the factor $\mathfrak{s}$. This implies that $B(\lambda)$ is analytic at $\lambda\in \sigma(T_0)$.

That $T$ is of Weyl order $1/2$ is the consequence of the well-known asymptotic behavior of eigenvalues of $T_0$ and Thm.~\ref{thm16}.
\end{proof}
\begin{corollary}\label{lem4}For any $\lambda\in \mathbb{R}$, $\frac{\mathfrak{s}}{(\mathfrak{c}-i\mathfrak{s})(\mathfrak{s}'-i\mathfrak{s})-1}\neq 0$.
\end{corollary}
\begin{proof}This is obvious from the proof of the above proposition.
\end{proof}
\begin{proposition}For a separated self-adjoint boundary condition, the (geometric or analytic) multiplicity of any eigenvalue is 1.
\end{proposition}
\begin{proof}As before, let $U=\left( \begin{array}{cc}                                                                                                                                             e^{i\alpha} & 0 \\
0 & e^{i\beta} \\
 \end{array}\right)$
 be the matrix parameterizing the separated self-adjoint boundary condition. $\lambda$ is an eigenvalue for the boundary value problem if and only if
 \[\det(U-B(\lambda))=0;\]
in particular, $\lambda$ is of multiplicity 2 if and only if $B(\lambda)=U$. However, according to our previous corollary, $B(\lambda)$ for $\lambda\in \mathbb{R}$ is never diagonal. Due to Thm.~\ref{thm8}, this proves the claim.
\end{proof}
\begin{proposition}Only for real self-adjoint boundary conditions can the corresponding Sturm-Liouville problems have eigenvalues with multiplicity 2.
\end{proposition}
\begin{proof}By Prop.~\ref{p17}, the Weyl curve $W_T(\lambda)$ represents a \emph{real} boundary condition at each point $\lambda\in \mathbb{R}$. Therefore for a self-adjoint boundary condition $U$ which is not real, we can never have $U=B(\lambda)$. The claim follows.
\end{proof}
\emph{Remark}. Traditionally, the above two propositions are proved by checking that the derivative of the determinant at $\lambda$ is non-vanishing. Our proof here is much more structural.

For a general entire operator with deficiency indices $(2,2)$, its Weyl curve lies in $Gr(2,4)$, which is a 4-dimensional complex manifold, but for our $T$, $W_T(\lambda)$ lies in $\mathfrak{L}(\mathcal{B}_T)$, which is a 3-dimensional complex manifold. In this sense, $T$ is algebraically degenerate. Furthermore, even with this in mind, the corresponding real entire operator can still be algebraically degenerate because $T$ may have a nontrivial symmetry.

For $f\in L^2(0,\pi)$, define $(Uf)(x)=f(\pi-x)$. Then $U$ is a unitary operator on $L^2(0,\pi)$. By this $U$, $T_q$ is transformed into another Sturm-Lioville operator $T_{\tilde{q}}$ with $\tilde{q}(x)=q(\pi-x)$. Certainly $T_q$ and $T_{\tilde{q}}$ ought to have the same Weyl class.

If $q(x)=q(\pi-x)$ for any $x\in[0, \pi]$, then $U$ is obviously a symmetry of $T_q$. We note that
\[\Gamma_0Uy=\left(
               \begin{array}{c}
                 y(\pi) \\
                 y(0) \\
               \end{array}
             \right)=\left(
                       \begin{array}{cc}
                         0 & 1 \\
                       1 & 0 \\
                       \end{array}
                     \right)\left(
                              \begin{array}{c}
                                y(0) \\
                                y(\pi) \\
                              \end{array}
                            \right)=\left(
                       \begin{array}{cc}
                         0 & 1 \\
                       1 & 0 \\
                       \end{array}
                     \right)\Gamma_0y,
\]
and similarly,
\[\Gamma_1Uy=\left(
                       \begin{array}{cc}
                         0 & 1 \\
                       1 & 0 \\
                       \end{array}
                     \right)\Gamma_1y.\]
Let $M(\lambda)$ be the associated Weyl function. Then for $\lambda\notin\mathbb{R}$ and $y\in\textup{ker}(T^*-\lambda)$, we have
                   \[\Gamma_1Uy=\left(
                             \begin{array}{cc}
                             0 & 1 \\
                               1 & 0 \\
                           \end{array}
                        \right)
                \Gamma_1y=\left(
                     \begin{array}{cc}
                        0 & 1 \\
                       1 & 0 \\
                      \end{array}
                   \right)M(\lambda)\Gamma_0y.\]
On the other side, $Uy\in \textup{ker}(T^*-\lambda)$ and thus
               \[\Gamma_1Uy=M(\lambda)\Gamma_0Uy=M(\lambda)
               \left(
               \begin{array}{cc}
                        0 & 1 \\
                       1 & 0 \\
                       \end{array}
                    \right)
                    \Gamma_0y.\]
We ought to have
\[\left(
\begin{array}{cc}
                        0 & 1 \\
                                              1 & 0 \\
                       \end{array}
                     \right) M(\lambda)=M(\lambda)\left(\begin{array}{cc}
                         0 & 1 \\
                       1 & 0 \\
                       \end{array}
                     \right).\]
This identity poses the additional requirement that $\mathfrak{c}\equiv \mathfrak{s}'$, implying that $W_T(\lambda)$ should lie in a subvariety of codimension 1 in $\mathfrak{L}(\mathcal{B}_T)$. This is actually not strange: the existence of the nontrivial symmetry $U$ means that $T$ is reducible and can be decomposed into the sum of two symmetric operators with deficiency indices $(1,1)$. These operators are just the restrictions of $T$ on the subspaces of $D(T_0)$ consisting of functions such that $Uf\equiv \pm f$. Therefore, $W_T(\lambda)$ has to lie in the quadratic surface $\mathbb{CP}^1\times \mathbb{CP}^1$. This is an example where a symmetry results in algebraical degeneracy. It is interesting to know whether this is the only way the real symmetric operator acquires its degeneracy.

There are two transversal Picard exceptional boundary conditions $y(0)=y'(0)=0$ and $y(\pi)=y'(\pi)=0$: by Cauchy's existence and uniqueness theorem for initial value problems of ordinary differential equations, the whole $\mathbb{C}$ is the resolvent set for both these two boundary conditions. These boundary conditions are just the reason for the well-known fact that for a generic boundary condition the eigenvalues are precisely zeros of an entire function. The viewpoint of value distribution theory implies that though no eigenvalue exists for these two boundary conditions, the Weyl curve approaches the corresponding Schubert hyperplanes sufficiently often.

Let us show degenerate boundary conditions can really occur. Assume $q\equiv 0$ and consider the boundary conditions\\
\begin{center}
(I)\, $\left \{
\begin{array}{ll}y(0)+y(\pi)=0,\\ y'(0)=y'(\pi).
\end{array}
 \right.$
 and
(II)\, $\left \{
\begin{array}{ll}y(0)=y(\pi),\\ y'(0)+y'(\pi)=0.
\end{array}
 \right.$
\end{center}

\begin{proposition}The above two boundary conditions are degenerate w.r.t. $T_q$ with $q\equiv 0$.
\end{proposition}
\begin{proof}We only consider the first boundary condition and the following argument applies to the second as well. Let $y=a \cos \omega x +b \frac{\sin \omega x}{\omega} $ be the solution for the boundary value problem, where $\omega=\sqrt{\lambda}$ (any branch of $\sqrt{\lambda}$ is OK) and $a, b\in \mathbb{C}$ are constants to be determined. $a, b$ have to satisfy the equations
\begin{equation*}\left \{
\begin{array}{ll}
 (1+\cos(\omega\pi))a+\frac{\sin(\omega\pi)}{\omega}b=0, \\
  \omega\sin(\omega\pi)a+(1-\cos(\omega\pi))b=0.
\end{array}
 \right.
 \end{equation*}
Obviously the determinant of the coefficient matrix is zero for any $\lambda\in \mathbb{C}$. The corresponding operator thus has any $\lambda\in \mathbb{C}$ as its eigenvalue whose geometric multiplicity is 1.
\end{proof}
Thus $T_q$ with $q\equiv 0$ provides an example of an entire operator which is weakly algebraically degenerate. This degeneracy also introduces many other Picard boundary conditions. Let us consider the following boundary conditions
\begin{center}(III) $\left \{
\begin{array}{ll}zy(0)+y(\pi)=0,\\ zy'(0)=y'(\pi).
\end{array}
 \right.$
 \end{center}
 where $z\neq \pm 1$ is a complex constant. It can be checked easily that all these boundary conditions are Picard exceptional and have no eigenvalues at all. It's interesting to point out that the two degenerate boundary conditions (I) (II) and these Picard boundary conditions (III) (with $z=\infty$ included) form a $\mathbb{CP}^1$ in $Gr(2,4)$.

The above weak algebraic degeneracy comes from the fact that $\mathfrak{c}$, $\mathfrak{c}'$, $\mathfrak{s}$ and $\mathfrak{s}'$, as entire functions, are linearly dependent because $\mathfrak{c}\equiv\mathfrak{s}'$. In fact, we have
\begin{proposition}\label{p19}$T_q$ is weakly algebraically non-degenerate if and only if $\mathfrak{c}\not\equiv\mathfrak{s}'$.
\end{proposition}
\begin{proof}We identify $\mathcal{B}_T$ with $\mathbb{C}^4$, simply by taking $(f(0),f'(0),f(\pi),f'(\pi))$ for $[f]\in \mathcal{B}_T$. Let $\{e_i\}$ ($i=1,\cdots, 4$) denote the standard basis of $\mathbb{C}^4$ and $e_{ij}$ ($i\neq j$) the wedge product $e_i\wedge e_j$. Since $\{c(x,\lambda), s(x,\lambda)\}$ is a basis of $\textup{ker}(T^*_q-\lambda)$, $W_T(\lambda)\in Gr(2,4)$ can be represented by the $2\times 4$ matrix
\[\left(
    \begin{array}{cccc}
      1 & 0 & \mathfrak{c} & \mathfrak{c}' \\
      0 & 1 & \mathfrak{s} & \mathfrak{s}' \\
    \end{array}
  \right).
\]
Then $\det W_T(\lambda)$ has a basis
\[e(\lambda):=e_{12}+\mathfrak{s} e_{13}+\mathfrak{s}'e_{14}-\mathfrak{c}e_{23}-\mathfrak{c}'e_{24}+e_{34},\]
where the identity $\mathfrak{c}\mathfrak{s}'-\mathfrak{c}'\mathfrak{s}\equiv1$ was used. Note that for $|\lambda|$ large enough, we have the following asymptotic formulae (see for example \cite[\S~1.1]{freiling2001inverse})
\begin{eqnarray*}\mathfrak{c}&=&\cos (\omega\pi)+c_1\frac{\sin(\omega\pi)}{2\omega}\\&+&(q(\pi)-q(0)-c_2)\frac{\cos(\omega\pi)}{4\omega^2}+O(\frac{e^{|\tau|\pi}}{|\omega|^3}),\end{eqnarray*}
\[\mathfrak{c}'=-\omega\sin (\omega\pi)+O(e^{|\tau|\pi}),\quad
\mathfrak{s}=\frac{\sin(\omega\pi)}{\omega}+O(\frac{e^{|\tau|\pi}}{|\omega|^2}),\]
\begin{eqnarray*}\mathfrak{s}'&=&\cos (\omega\pi)+c_1\frac{\sin(\omega\pi)}{2\omega}\\
&-&(q(\pi)-q(0)+c_2)\frac{\cos(\omega\pi)}{4\omega^2}+O(\frac{e^{|\tau|\pi}}{|\omega|^3}),\end{eqnarray*}
where $\tau=\Im \omega$ and
\[c_1=\int_0^\pi q(x)dx,\quad c_2=\int_0^\pi [q(x)\int_0^xq(\sigma)d\sigma]dx.\]

A generic boundary condition $b\in Gr(2,4)$ now can be represented by a $2\times 4$ matrix $(A|B)$ of rank 2, where $A, B\in M_{2\times 2}$.

(1) If $\det A=0$, then the boundary condition $b$ can be written in either of the following five forms:
\[\left(
    \begin{array}{cccc}
      1 & b_{11} & 0 & b_{12} \\
      0 & 0 & 1 & b_{22} \\
    \end{array}
  \right),\left(
    \begin{array}{cccc}
      1 & b_{11} & b_{12} & 0 \\
      0 & 0 & 0 & 1 \\
    \end{array} \right),\left(
    \begin{array}{cccc}
      0 & 0 & 1 & b_{12} \\
      0 & 1 & 0 & b_{22} \\
    \end{array}\right),
\]
\[\left(
    \begin{array}{cccc}
      0 & 0 & 0 & 1 \\
     0 & 1 & b_{22} & 0 \\
    \end{array}
  \right), \left(
    \begin{array}{cccc}
      0 & 0 & 1 & 0 \\
     0 & 0 & 0 & 1 \\
    \end{array}
  \right).
\]
Each of them cannot be degenerate w.r.t. $T_q$. Indeed, take
\[(A|B)=\left(
    \begin{array}{cccc}
      1 & b_{11} & 0 & b_{12} \\
      0 & 0 & 1 & b_{22} \\
    \end{array}
  \right)\]
  for example. Then $\det V_b$ has a basis
  \[s_b:=e_{13}+b_{22}e_{14}+b_{11}e_{23}+b_{11}b_{22}e_{24}-b_{12}e_{34}.\]
In terms of the basis $e_1\wedge e_2\wedge e_3\wedge e_4$ of $\det \mathbb{C}^4$,
\[s_b(\sigma(\lambda))=\mathfrak{c}'-b_{11}b_{22}\mathfrak{s}+b_{11}\mathfrak{s}'-b_{22}\mathfrak{c}-b_{12}.\]
For other four boundary conditions, we have \[s_b(\sigma(\lambda))=-\mathfrak{c}-b_{11}\mathfrak{s}+b_{12}, \,-\mathfrak{s}'+b_{12}\mathfrak{s}+b_{22}, \,\mathfrak{s}-b_{22} \,\,\textup{or}\,\, 1.\]
By the asymptotic behavior of $\mathfrak{c}, \mathfrak{c}', \mathfrak{s}, \mathfrak{s}'$, we see $s_b(\sigma(\lambda))\not\equiv 0$ holds in each case.

(2) If $\det A\neq 0$, w.l.g., we can set
\[(A|B)=\left(
          \begin{array}{cccc}
            1 & 0 & b_{11} & b_{12} \\
            0 & 1 & b_{21} & b_{22} \\
          \end{array}
        \right).
\]
Then $\det V_b$ has the basis
\[s_b:=e_{12}+b_{21}e_{13}+b_{22}e_{14}-b_{11}e_{23}-b_{12}e_{24}+(b_{11}b_{22}-b_{12}b_{21})e_{34}.\]
In terms of the basis $e_1\wedge e_2\wedge e_3\wedge e_4$ of $\det \mathbb{C}^4$,
\[s_b(\sigma(\lambda))=(b_{11}b_{22}-b_{12}b_{21})+b_{12}\mathfrak{s}-b_{11}\mathfrak{s}'-b_{22}\mathfrak{c}+b_{21}\mathfrak{c}'+1.\]
If $\mathfrak{s}'\equiv \mathfrak{c}$, we can choose $b_{11}=-b_{22}=1$ and $b_{12}=b_{21}=0$ and get a degenerate boundary condition. Conversely, if $\mathfrak{s}'\not\equiv \mathfrak{c}$ and there is a degenerate boundary condition, then by the above asymptotic formulae of $\mathfrak{c}, \mathfrak{c}', \mathfrak{s}, \mathfrak{s}'$ we should have $b_{21}=0$, $b_{11}+b_{22}=0$ and consequently $b_{11}^2=1$. W.l.g, assume $b_{11}=1$. Then
\[0\equiv s_b(\sigma(\lambda))=b_{12}\mathfrak{s}+(\mathfrak{c}-\mathfrak{s}').\]
By the asymptotic formulae, we see, for $|\lambda|$ large enough, $$\mathfrak{c}-\mathfrak{s}'=O(\frac{e^{|\tau|}}{|\omega|^2}).$$
This shows that $\mathfrak{s}$ and $\mathfrak{c}-\mathfrak{s}'$ are linearly independent, a contradiction!
\end{proof}
\begin{proposition}If $q(\pi)\neq q(0)$, then $T_q$ is weakly algebraically non-degenerate.
\end{proposition}
\begin{proof}
This is because for $|\lambda|$ large enough,
\[\mathfrak{c}-\mathfrak{s}'=(q(\pi)-q(0))\frac{\cos(\omega\pi)}{2\omega^2}+O(\frac{e^{|\tau|\pi}}{|\omega|^3}).\]
\end{proof}
Thus if $q(x)=q(\pi-x)$ for any $x\in [0,\pi]$, then $T_q$ is weakly algebraically degenerate and the boundary conditions (I) and (II) are degenerate for all such $T_q$'s. We can also find that in this case the boundary condition (III) with $z\neq \pm 1$ is Picard exceptional.
\begin{proposition}If $T_q$ is weakly algebraically non-degenerate, then for all $y\in Gr(2,4)$ but the two Picard exceptional boundary conditions, the corresponding boundary value problems each have countably infinite eigenvalues with no finite accumulation point.
\end{proposition}
\begin{proof}The analysis in the proof of Prop.~\ref{p19} shows that $s_b(\sigma(\lambda))$ is a non-constant entire function of order $1/2$ if $y$ is neither of the two Picard exceptional boundary conditions. Since a transcendental meromorphic function of non-integer order cannot have more than one Picard exceptional value \cite[Thm.~1.1, Chap.~4]{gol'dberg2008value}, $\infty$ is the only Picard exceptional value of $s_b(\sigma(\lambda))$. Thus $s_b(\sigma(\lambda))$ has countably infinite zeros.
\end{proof}
It's interesting to know if there are Nevanlinna exceptional boundary conditions that are not Picard exceptional.
\subsection{Inverse problems}
Using so-called spectral data (e.g., spectra, Weyl function) to determine the underlying symmetric operator (and its related extensions) is the central problem of inverse spectral theory. In this direction, the inverse Sturm-Liouville theory is the most well-established. In our viewpoint, the correspondence between unitary equivalence classes of simple symmetric operators and congruence classes of Weyl curves should be the basis of such investigations.

In 1946, G. Borg proved in \cite{Borg1946umkehrung} the following now famous and basic result in inverse Sturm-Liouville theory:
\begin{theorem}Given $L=-\frac{d^2}{dx^2}+q(x)$ on $[0,\pi]$. Let $\lambda_0<\lambda_1<\lambda_2<\cdots$ and $\mu_0<\mu_1<\mu_2<\cdots$ be the eigenvalues of $Ly=\lambda y$ with boundary conditions
\begin{center}
 $\left \{
\begin{array}{ll}y'(0)-hy(0)=0,\\ y'(\pi)+Hy(\pi)=0.
\end{array}
 \right.$
 and
$\left \{
\begin{array}{ll}y'(0)-hy(0)=0,\\ y'(\pi)+H'y(\pi)=0.
\end{array}
 \right.$
\end{center}
respectively, where $h$, $H$ and $H'$ are real constants. Then the sequences $\{\lambda_n\}_0^\infty$ and $\{\mu_n\}_0^\infty$ determine $q$, $h$, $H$ and $H'$ uniquely.
\end{theorem}
This theorem also holds if $H'$ is replaced by $\infty$ (then the boundary condition is $y(\pi)=0$), see \cite[Thm.~1.4.4]{freiling2001inverse}. It is important for the two boundary conditions to share a common boundary condition at the left endpoint. V. Pierce constructed an uncountable family of potentials sharing the same Dirichlet and Neumann spectra in \cite{pierce2002determining}.

We give an interpretation of Borg's result as follows. The common boundary condition $y'(0)-hy(0)=0$ gives a symmetric extension $T_h$ of $T$. The domain $D(T_h)$ of $T_h$ is
\[\{f\in L^2(0,\pi)|f, f'\in AC[0,\pi],\,f'(0)=hf(0),\, f(\pi)=f'(\pi)=0,\, Lf\in L^2(0,\pi)\},\]
and
\[D(T^*_h)=\{f\in L^2(0,\pi)|f, f'\in AC[0,\pi],\, Lf\in L^2(0,\pi),\,f'(0)=hf(0)\}.\]
$T_h$ has deficiency indices $(1,1)$ and the two boundary conditions $y'(\pi)+Hy(\pi)=0$ and $y'(\pi)+H'y(\pi)=0$ produce two self-adjoint extensions $T_h^H$, $T_h^{H'}$ of $T_h$. Then Borg's result precisely means the spectra of $T_h^H$, $T_h^{H'}$ determine $T_h$, $H$ and $H'$ completely. This is an improvement of Prop.~\ref{p15} in the specified setting.

Borg's result has two immediate consequences.
\begin{proposition}The potential $q$ of the Sturm-Liouville operator $-\frac{d^2}{dx^2}+q(x)$ is uniquely determined by $M(\lambda)$ in (\ref{eq9}).
\end{proposition}
\begin{proof}Choose real numbers $h$, $H$ and $H'$ and construct the two boundary conditions as in Borg's result. We can obtain the spectra of these boundary value problems by solving the corresponding equation $\det(U-B(\lambda))=0$ where $B(\lambda)$ is the contractive Weyl function determined by $M(\lambda)$ and $U$ the corresponding unitary matrix in terms of the boundary triplet $(\mathbb{C}^4, \Gamma_0, \Gamma_1)$. Applying Borg's result yields the claim.
\end{proof}
\emph{Remark}. Thus the three entire functions $\mathfrak{s}$, $\mathfrak{s}'$ and $\mathfrak{c}$ completely determine $q$.
\begin{proposition}The Sturm-Liouville operator $T_h$ is uniquely determined by its Weyl class.
\end{proposition}
\begin{proof}It's easy to find that all self-adjoint boundary conditions are of the form $y'(\pi)+Hy(\pi)=0$ with $H\in \mathbb{R}\cup \{\infty\}$, where $H=\infty$ corresponds to $y(\pi)=0$. If $m(\lambda)$ is a Weyl function for $T_h$ (we don't know the associated boundary triplet for we only know the Weyl class.) and $B(\lambda)$ its corresponding contractive Weyl function, then we can solve the equations $B(\lambda)\pm 1=0$ and obtain the spectra of two self-adjoint extensions of $T_h$. Applying Borg's theorem leads to our claim.
\end{proof}
\emph{Remark}. Since in this case the Weyl class is completely determined by the first Chern function $c_1$, thus $c_1$ can be viewed as the spectral data to determine $q$ and $h$. It's interesting to know how to recover them from $c_1$ explicitly.

We formulate a question and a conjecture to conclude this subsection. Borg's theorem is special for the symmetric operator $T$ in that the chosen two boundary conditions are not transversal. Since transversality is an open condition, we would like to know the answer to the following question:

\textbf{Question}: How many spectra are necessary and sufficient to determine $q$ if all the involved self-adjoint boundary conditions are \emph{mutually} transversal?

We don't know if there is a finite answer, but we believe the answer will help to clarify how to formulate inverse spectral problems in the general setting of higher deficiency indices. Our previous two propositions also suggest the following conjecture:

\textbf{Conjecture}. If $T_q$ and $T_{\tilde{q}}$ have the same Weyl class, then either $q(x)=\tilde{q}(x)$ or $q(x)=\tilde{q}(\pi-x)$ for all $x\in [0,\pi]$.

This may be proved by carefully examining the asymptotic behavior of the Weyl functions in the Weyl class as $\lambda\rightarrow \infty$, but we won't pursue it here any further.
\section*{Acknowledgemencts}
The author wants to express his deep gratitude towards Prof. Huang Zhenyou for his valuable discussions and suggestion during the preparation of this work and also for his encouragement through all these years.
\section*{Appendix A: Groups and group actions}
This appendix collects the basic terminology concerning group actions. It is only for the convenience of readers unfamiliar with the material, which is actually standard and can be found in many textbooks on Lie groups and differential geometry. We refer the readers to \cite{nicolaescu2020lectures}.

We say a group $G$ (with identity $e$) acts on a set $M$ if there is a group homomorphism $\rho$ between $G$ and the permutation group $\textup{Per}(M)$. If $m\in M$ and $g\in G$, then $\rho(g)(m)$ is usually denoted by $g\cdot m$. For $m\in M$, the set $G_m:=\{g\in G|g\cdot m=m\}$ is obviously a subgroup of $G$, called the isotropy subgroup at $m$. The $G$-action is called effective, if for any $e\neq g\in G$, there is a certain $m\in M$ such that $g\cdot m\neq m$. The set $O(m):=\{g\cdot m\in M|g\in G\}$ is called a $G$-orbit through $m$. $M$ can be partitioned into different $G$-orbits. The set of $G$-orbits in $M$ is denoted by $M/G$. If $M/G$ is a singleton, we say the $G$-action is transitive.

If $G$ is a topological group and $M$ a topological space, the $G$-action is continuous if the map $G\times M\rightarrow M, (g,m)\mapsto g\cdot m$ is continuous. Then each $\rho(g)$ for $g\in G$ is an auto-homeomorphism of $M$. The set $M/G$ with its quotient topology is called the orbit space. To obtain a good orbit space, generally the $G$-action on $M$ must be under control. If $G$ is further a smooth manifold and its structural maps are smooth, $G$ is called a Lie group. A Lie group $G$ acts smoothly on a manifold $M$ if the map $G\times M\rightarrow M, (g,m)\mapsto g\cdot m$ is smooth. Let $G$ act on two manifolds $M$ and $N$. A smooth map $f:M\rightarrow N$ is called $G$-equivariant, if $f(g\cdot m)=g\cdot f(m)$ for any $g\in G$ and $m\in M$. If $M$ is a complex manifold and each element of $G$ acts on $M$ as an automorphism of $M$, we say $G$ acts holomorphically.

For a Hilbert space $H$, $\mathbb{GL}(H)$ consists of invertible elements in $\mathbb{B}(H)$. $\mathbb{GL}(H)$ is a Banach-Lie group and called the general linear group on $H$. $\mathbb{GL}(H)$ has a normal subgroup $\{c\times Id|c\in \mathbb{C}^*\}$ and the quotient group $\mathbb{PGL}(H)$ is called the projective general linear group. The subgroup $\mathbb{U}(H)\subset \mathbb{GL}(H)$ consists of unitary operators on $H$ and is called the unitary group on $H$. If $\dim H=n\in \mathbb{N}$, $\mathbb{U}(H)$ is often written as $\mathbb{U}(n)$. The projective unitary group $\mathbb{PU}(H)$ can be defined as well.

A homomorphism $\rho$ from a group $G$ to $\mathbb{GL}(H)$ is called a representation of $G$ on $H$.
\section*{Appendix B: Linear relations}
We collect some material on linear relations, which is needed for treating simple symmetric operators that are not densely defined. Our goal is to show how to adapt the relevant notions to the new setting. More technical details can be found in \cite{behrndt2020boundary}.

For a Hilbert space $H$, we equip $\mathbb{H}:=H\oplus_\bot H$ with its second standard strong symplectic structure $[\cdot, \cdot]_{new}$. We use $\pi_1$ and $\pi_2$ to denote the projections from $\mathbb{H}$ to the two copies of $H$. A (closed) subspace $A$ in $\mathbb{H}$ is called a (closed) linear relation in $H$. It is a generalization of the graph of a (closed) operator. All linear relations in this appendix are closed.

 For a linear relation $A$. $\pi_1(A)$ (resp. $\pi_2(A)$) is called the domain (resp. range) of $A$. $\textup{ker}A:=\{x\in H|(x,0)\in A\}$ is the kernel of $A$ and $A_{\textup{mul}}:=\{x\in H|(0,x)\in A\}$ is the multi-valued part of $A$. The inverse of $A$ is the linear relation
$A^{-1}=\{(x,y)\in \mathbb{H}|(y,x)\in A\}$. If $\textup{ker}A=0$, then $A^{-1}$ is actually an operator. A linear relation in $H$ is Fredholm if $\dim\textup{ker}A<+\infty$, $\pi_2(A)$ is closed and $\textup{codim}\pi_2(A)<+\infty$. Furthermore $\dim\textup{ker}A-\textup{codim}\pi_2(A)$ is called the index of $A$.

For $\lambda\in \mathbb{C}$,
\[A-\lambda:=\{(x,y-\lambda x)\in \mathbb{H}|(x,y)\in A\}.\]
If $(A-\lambda)^{-1}$ is the graph of an operator in $\mathbb{B}(H)$, $\lambda$ is called a regular value of $A$. The resolvent set $\rho(A)$ of $A$ is the set of all regular values of $A$. $\sigma(A):=\mathbb{C}\backslash \rho(A)$ is called the spectrum of $A$. $\sigma(A)$ can be decomposed into three disjoint parts:
\[\sigma_p(A)=\{\lambda\in \mathbb{C}|\textup{ker}(A-\lambda)\neq 0\},\]
\[\sigma_c(A)=\{\lambda\in \mathbb{C}|\textup{ker}(A-\lambda)=0,\, \overline{\pi_2(A-\lambda)}=H,\, \lambda\not\in \rho(A)\},\]
\[\sigma_r(A)=\{\lambda\in \mathbb{C}|\textup{ker}(A-\lambda)=0,\, \overline{\pi_2(A-\lambda)}\neq H\}.\]
These are called point spectrum, continuous spectrum and residual spectrum of $A$ respectively.

The symplectic complement $A^*$ of $A$ is called the adjoint relation of $A$. We have \begin{itemize}
                                                                                        \item $H=\textup{ker}A\oplus_\bot \overline{\pi_2(A^*)}$,
                                                                                        \item (Closed Range Theorem) $\pi_2(A)$ is closed if and only if $\pi_2(A^*)$ is closed.
                                                                                      \end{itemize}
For details of these claims, see Prop.~1.3.2 and Thm.~1.3.5 in \cite{behrndt2020boundary}. For a linear relation $A$, \begin{itemize}
                                                                                                        \item $\lambda\in \rho(A)$ if and only if $\bar{\lambda}\in \rho(A^*)$;
                                                                                                        \item $\lambda\in \sigma_r(A)$ if and only if $\textup{ker}(A-\lambda)=0$ and $\lambda\in \sigma_p(A^*)$.
                                                                                                      \end{itemize}
See Prop.~1.3.10 in \cite{behrndt2020boundary}.

If $A$ is isotropic or $A\subset A^*$, $A$ is called a symmetric relation. If $A=A^*$, $A$ is called a self-adjoint relation. If $S$ is a self-adjoint relation such that $A\subset S \subset A^*$, $S$ is called a self-adjoint extension of $A$. If $T$ is a closed simple symmetric operator not densely defined in $H$, its graph $A_T$ is a symmetric relation. Then $A_T^*$ is only a linear relation rather than an operator, but $A_T^*/A_T$ is still a strong symplectic Hilbert space. If the signature $(n_+,n_-)$ is such that $n_\pm\geq 1$, $\Gamma_\pm$ can still be defined on $A_T^*$ and a Green's formula holds.

If $A$ is symmetric, $\dim\text{ker}(A^*-\lambda)$ for $\lambda\in \mathbb{C}_+\cup \mathbb{C}_-$ is locally constant \cite[Thm.~1.2.5]{behrndt2020boundary}. These constants are just the above $n_\pm$. Let $\widehat{\text{ker}(A^*-\lambda)}=\{(x, \lambda x)|x\in \text{ker}(A^*-\lambda)\}$. If $\lambda\in \rho(B)$ for an extension $B$ of $A$. Then the following Krein decomposition holds:
\[A^*=B\oplus \widehat{\text{ker}(A^*-\lambda)}.\]
The Weyl curve of $A$ now can be defined in the same way as in \S~\ref{sec4}. If a symplectic isomorphism $\Phi$ between $A_T^*/A_T$ and $H_+\oplus_\bot H_-$ is chosen, then the extensions $A_\pm$ determined by $\Gamma_\pm \hat{x}=0$ for $\hat{x}\in A_T^*$ are actually maximal accumulative and dissipative relations in $H$ respectively. By Coro.~1.6.5 in \cite{behrndt2020boundary}, $\mathbb{C}_\pm\subset \rho(A_\pm)$. Then Prop.~\ref{p2} and Corol.~\ref{c1} in \S~\ref{sec4} continue to hold even the simple symmetric operator $T$ is not densely defined.

$U\in \mathbb{U}(H)$ is called a symmetry of a linear relation $A$ in $H$, if $U\cdot A:=\{(Ux,Uy)\in \mathbb{H}|(x,y)\in A\}=A$.

As was implied in \S~\ref{sec3}, for a self-adjoint relation $S$ in $H$, there is a closed subspace $K \subset H$ and a self-adjoint operator $A$ in $K$ such that
\[S=\{(x,Ax+y)\in \mathbb{H}| x\in D(A),\, y\in K^\bot \}.\]
Then $\sigma(S)=\sigma(A)$ and $S_{\textup{mul}}=K^\bot$. If furthermore $S$ is a self-adjoint extension of a simple symmetric operator with deficiency indices $(n,n)$ where $n\in \mathbb{N}$, then it's easy to see $\dim S_{\textup{mul}}\leq n$.

\end{document}